\documentclass[a4paper, 11pt]{amsart}
\usepackage{amssymb,amsmath,txfonts,mathrsfs,titletoc,appendix}
\usepackage[alphabetic]{amsrefs}
\usepackage{geometry}
\newtheorem{assu}{Assumption}
\newtheorem{theorem}{Theorem}[section]
\newtheorem{prop}[theorem]{Proposition}
\newtheorem{lemma}[theorem]{Lemma}
\newtheorem{remark}[theorem]{Remark}
\newtheorem{claim}[theorem]{Claim}

\newtheorem{question}[theorem]{Question}
\newtheorem{definition}[theorem]{Definition}
\newtheorem{cor}[theorem]{Corollary}

\numberwithin{equation}{section}

\def\pf{{\it Proof:}~}

\begin{document}

\title[Large time behavior of the heat kernel]{Large time behavior of the heat kernel}
\author{Guoyi Xu}
\address{Mathematical Sciences Center\\Tsinghua University, Beijing\\P. R. China, 100084}
\email{gyxu@math.tsinghua.edu.cn}
\date{\today}

\begin{abstract}
In this paper, we study the large time behavior of the heat kernel on complete Riemannian manifolds with nonnegative Ricci curvature, which was studied by P. Li with additional maximum volume growth assumption. Following Y. Ding's original strategy, by blowing down the metric, using Cheeger and Colding's theory about limit spaces of Gromov-Hausdorff convergence, combining with the Gaussian upper bound of heat kernel on limit spaces, we succeed in reducing the limit behavior of the heat kernel on manifold to the values of heat kernels on tangent cones at infinity of manifold with renormalized measure. As one application, we get the consistent large time limit of heat kernel in more general context, which generalizes the former result of P. Li. Furthermore, by choosing different sequences to blow down the suitable metric, we show the first example manifold whose heat kernel has inconsistent limit behavior, which answers an open question posed by P. Li negatively. 
\end{abstract}


\maketitle

\titlecontents{section}[0em]{}{\hspace{.5em}}{}{\titlerule*[1pc]{.}\contentspage}
\titlecontents{subsection}[1.5em]{}{\hspace{.5em}}{}{\titlerule*[1pc]{.}\contentspage}
\tableofcontents

\section{Introduction}\label{SECTION 01}

On $(M^n, g)$, we consider the fundamental solution $H(x, y, t)$, which solves the heat equation with initial data:
\begin{equation}\nonumber
\left\{
\begin{array}{rl}
\big(\frac{\partial}{\partial t}- \Delta\big) F(x, t) &= 0 \quad \quad \quad on\ M^n\times (0, \infty) \\
F(x, 0)&= f(x)  \quad \quad on\ M^n\\
\end{array} \right.
\end{equation}
by setting
\begin{align}
F(x, t)= \int_{M^n} H(x, y, t) f(y) dy \nonumber
\end{align}

It is well-known that there exists a minimal positive fundamental solution of $(M^n, g)$ (cf. Theorem $12.4$ in \cite{Libook}). In \cite{Dodziuk}, J. Dodziuk showed that if the Ricci curvature $(M^n, g)$ is bounded from below, then the minimal positive fundamental solution of $(M^n, g)$ is the unique positive fundamental solution of $(M^n, g)$. In this case, we say that the unique positive fundamental solution $H(x, y, t)$ is the \textbf{heat kernel} of $(M^n, g)$.

Especially, when $(M^n, g)$ has non-negative Ricci curvature, in \cite{LY}, P. Li and S-T. Yau proved that for all $\epsilon> 0$, there exists constants $C(\epsilon)> 0$, such that
\begin{align}
\frac{C(\epsilon)^{-1}}{V\big(\sqrt{t}\big)}\exp\Big(-\frac{d^2(x, y)}{(4- \epsilon)t}\Big)\leq H(x, y, t)\leq \frac{C(\epsilon)}{V\big(\sqrt{t}\big)}\exp\Big(-\frac{d^2(x, y)}{(4+ \epsilon)t}\Big) \label{1.1}
\end{align}
where the terms $V\big(\sqrt{t}\big)$ and $d(x, y)$ denote the volume of the geodesic ball centered at $y$ of radius $\sqrt{t}$ and the geodesic distance from $x$ to $y$, respectively. 

In particular, there are constants $C_1(n)$ and $C_2(n)$ depending only on dimension $n$ of $M^n$, such that
\begin{align}
C_1(n)\leq \varliminf_{t\rightarrow \infty} V\big(\sqrt{t}\big) H(x, y, t)\leq \varlimsup_{t\rightarrow \infty} V\big(\sqrt{t}\big) H(x, y, t)\leq C_2(n) \label{1.2}
\end{align}

For smooth manifold $M^n$ with non-negative Ricci curvature, Bishop-Gromov volume comparison theorem asserts that the relative volume $\frac{V(r)}{r^n}$ is non-increasing in the radius $r$. As $r\rightarrow \infty$, it converges a non-negative number $\Theta$, which is called asymptotic volume ratio. If $\Theta> 0$, then we say that $M^n$ has maximal volume growth.

In \cite{Li}, P. Li initiated the study of large time behavior of heat kernel on open manifolds with $Rc\geq 0$ and maximal volume growth. Among other things, he proved the following theorem:
\begin{theorem}[P. Li]\label{thm Li}
{If $(M^n, g)$ has $Rc\geq 0$ and maximal volume growth, then 
\begin{align}
\lim_{t\rightarrow \infty} V\big(\sqrt{t}\big) H(x, y, t)= \omega(n) (4\pi)^{-\frac{n}{2}} \label{1.3}
\end{align}
where $\omega(n)$ is the volume of the unit $n$-ball in $\mathbb{R}^n$. 
}
\end{theorem}

The key of the proof is Li-Yau's Harnack inequality established in \cite{LY} and the Bishop-Gromov Volume Comparison Theorem. 

Inspired by the above work, in \cite{CM} T. Colding and W. Minicozzi studied the large scale behavior of the Green's function $G(x, y)$. Among other things, they proved
\begin{theorem}[T. Colding and W. Minicozzi]\label{thm CM}
{If $M^n$, $n\geq 3$ has nonnegative Ricci curvature and maximal volume growth, then for a fixed $x\in M^n$,
\begin{align}
\lim_{d(x, y)\rightarrow \infty} \frac{G(x, y)}{G_{\mathbb{R}^n} (x, y)} = \frac{\omega(n)}{\Theta} \nonumber
\end{align} 
where $G_{\mathbb{R}^n}(x, y)$ is the Green's function on $\mathbb{R}^n$.
}
\end{theorem}

And they also pointed out that the geometric motivation behind of Theorem \ref{thm CM} is the fact: every tangent cone at infinity of a manifold satisfying the assumptions of Theorem \ref{thm CM} is a metric cone, which was shown in \cite{CC}. 

Let us recall that for a complete  noncompact manifold $M^n$ with $Rc\geq 0$, a metric space $M_{\infty}$ is a tangent cone at infinity of $M^n$ if it is a Gromov-Hausdorff limit of a sequence of rescaled manifolds $(M^n, p, t_j^{-2}g)$, where $t_j\rightarrow \infty$. By Gromov's compactness theorem, \cite{Gromov}, any sequence $t_j\rightarrow \infty$, has a subsequence, also denoted as $t_j\rightarrow \infty$, such that the rescaled manifolds $(M^n, p, t_j^{-2}g)$ converge to some $M_{\infty}$ in the Gromov-Hausdorff sense. Example of Perelman (\cite{Pere}) shows that tangent cone at infinity is not unique in general even if the manifold with $Rc\geq 0$ has maximal volume growth and quadratic curvature decay. We refer the reader to \cite{CC1} for more examples including collapsing case. Note tangent cones at infinity of $M^n$ reflect the geometry at infinity of manifoold $M^n$. 

Later on, in \cite{LTW}, in addition to providing another proof of Theorem \ref{thm CM}, P. Li, L. Tam and J. Wang proved the sharp bound of the heat kernel under the assumption in Theorem \ref{thm Li}. Their sharp bound of heat kernel shows that the coefficients $\frac{C(\epsilon)^{-1}}{V\big(\sqrt{t}\big)}$ and $\frac{C(\epsilon)}{V\big(\sqrt{t}\big)}$ in (\ref{1.1}) have some relationship with the asymptotic volume ratio $\Theta$. 

As the asymptotic volume ratio is one quantity reflecting the geometry at infinity of manifolds, combined with the above observation about the Green's function and tangent cones at infinity of manifold, it is reasonable to speculate that Theorem \ref{thm Li} has one proof from the view point of tangent cones at infinity of manifold. In other words, the large time behavior of the heat kernel should have close relationship with the geometry at infinity of manifolds.


In \cite{Ding}, under the maximum volume growth assumption, Y. Ding reduced the study of large scale behavior of the Green's function and large time behavior of the heat kernel, to the analysis on tangent cones at infinity of manifolds, where all tangent cones are metric cones and the Gromov-Hausdorff convergence is non-collapsing. Note the analysis on metric cones had been done by J. Cheeger \cite{Cheeger0} in different context. By the above strategy, Y. Ding provided one alternative proof for Theorem \ref{thm Li} and Theorem \ref{thm CM} in unified way. 

However, as pointed out in \cite{Li}, the answer to the following question was still unknown:
\begin{question}\label{ques 1}
{Does $\displaystyle \lim_{t\rightarrow \infty} V\big(\sqrt{t}\big) H(x, y, t)$ exist generally without the assumption of maximal volume growth? 
}
\end{question}

To study the above question, we firstly set up the setting as the following:

\textbf{Blow Down Setup}: Note that $(M^n, g, \mu)$ is a complete Riemannian manifold with $Rc\geq 0$, where $\mu$ is the volume element determined by the metric $g$. We can define $(M_{i}, y, \rho_i,  \nu_i)$, where $M_i$ is the same differential manifold as $M^n$, $\rho_i$ is the metric defined as $\rho_i= t_i^{-1}g$, $\{t_i\}_{i= 1}^{\infty}$ is an increasing positive sequence whose limit is $\infty$, and $y$ is a fixed point on $M_i= M^n$. $\nu_i$ is a Borel regular measure defined by 
\begin{align}
\nu_i (A)\doteqdot \Big(\int_{B_i(1)} 1 d\mu_{i}\Big)^{-1} \Big(\int_{A} 1 d\mu_{i}\Big)= t_i^{\frac{n}{2}} V(\sqrt{t_i})^{-1}\mu_i(A) \label{3.0}
\end{align}
where $A\subset M_i$, $B_i(1)\doteqdot \{z\in M_i|\ d_{\rho_i}(z, y)\leq 1\}$, and $\mu_{i}$ is the volume element determined by $\rho_i$. Then by Gromov's compactness theorem (see \cite{Gromov}) and Theorem $1.6$ in \cite{CC1}, after passing to a suitable subsequence, we have
$(M_i, y, \rho_i, \nu_i)\stackrel{d_{GH}}{\longrightarrow} (M_{\infty}, y_{\infty},\rho_{\infty}, \nu_{\infty})$ in the measured Gromov-Hausdorff sense, where $\nu_{\infty}$ is the renormalized limit measure defined as in Section $1$ of \cite{CC1}. 

Unless otherwise mentioned, in this paper $(M^n, y, g, \mu)$, $(M_i^n, y,  \rho_i, \nu_i)$ and $(M_{\infty}, y_{\infty}, \rho_{\infty}, \nu_{\infty})$ are as in the above \textbf{Blow Down Setup} and $n\geq 3$.

A main result of this paper is the following:
\begin{theorem}\label{thm 8.13}
{Assume $(M_i, y, \rho_i, \nu_i)\stackrel{d_{GH}}{\longrightarrow} (M_{\infty}, y_{\infty},\rho_{\infty}, \nu_{\infty})$ as in the above \textbf{Blow Down Setup} and $n\geq 3$, then
\begin{align}
\lim_{i\rightarrow \infty} V(\sqrt{t_i})H(x, y, t_i)= p_{\infty}(y_{\infty}, y_{\infty}, 1)\, \label{8.13.1}
\end{align}
where $p_{\infty}$ is the heat kernel on the metric measure space $(M_{\infty}, y_{\infty},\rho_{\infty}, \nu_{\infty})$, and the convergence is point-wise convergence.
}
\end{theorem}

\begin{remark}\label{rem 0.1}
{In fact, after some suitable modification, it is not hard to show that the results of this paper also hold on complete Riemann surface, i.e. the $n= 2$ case. For space reason, we will not discuss the $n=2$ case separately here.
}
\end{remark}

To prove Theorem \ref{thm 8.13}, we follow Y. Ding's strategy loosely. However, by combining K.-T. Sturm's study about heat kernel on metric spaces (see \cite{Sturm1}, \cite{Sturm2}, \cite{Sturm3}, \cite{Sturm4}), with Cheeger-Colding's theory about spaces with Ricci curvature bounded from below (see \cite{CC}, \cite{CC1}, \cite{CC2}, \cite{CC3}, \cite{Cheeger}), we manage to overcome the difficulties caused by collapsing during Gromov-Hausdorff convergence. 


More concretely, in \cite{Ding}, the assumption of maximum volume growth was needed to get the Li-Yau's estimate for the Green's function on tangent cones at infinity of manifolds, then the reduction for the Green's function from manifolds to limit space under Gromov-Hausdorff convergence can be obtained, finally the reduction for the heat kernel as in Theorem \ref{thm 8.13} follows from the integral formula connecting the heat kernel with the Green's function. 

Our approach is kind of direct by avoiding the discussion of the Green's function. Note in Ding's proof, the Li-Yau's estimate for the Green's function on the limit spaces (metric cones) plays the essential role in getting the reduction for the Green's function. To get the reduction for the heat kernel, we need such an estimate for the heat kernel on the general limit spaces (metric measure spaces). Following K.-T. Sturm's method, we proved the general existence result and Gaussian-type upper bounds of heat kernel on $M_{\infty}$, which is enough for our use.

Note on compact domains, the heat kernel has the expansion determined by eigenvalues and eigenfunctions. On the other hand, J. Cheeger and T. Colding \cite{CC3} (also see \cite{Cheeger} for some technical details) had proved that the eigenvalues and eigenfunctions on compact metric measure spaces behave continuously under measured Gromov-Hausdorff convergence, which was originally conjectured by K. Fukaya in \cite{Fukaya}. Combining the suitable modifications of these two facts about heat kernel, eigenvalues and eigenfunctions on bounded domains, we can get the reduction of the heat kernel on bounded domains over complete manifolds, see Theorem \ref{thm 5.8}. 

Then applying the crucial Gaussian-type upper bounds of heat kernel on tangent cones at infinity of manifolds and the family of blowing down manifolds, using the suitable compact exhaustion of these complete blowing down manifolds, we succeed in getting the above reduction generally for the heat kernel on complete manifolds,  from the reduction of the heat kernel on bounded domains over complete manifolds.
Note the role of Gaussian-type upper bounds of heat kernel on tangent cones at infinity of manifolds and on blowing down manifolds, in getting our reduction, is analogous to the role that the uniform integrable function bound of measurable functions plays to guarantee two limit processes commute in Lebesgue's Dominated Convergence Theorem. 

A byproduct of the above general reduction result is, a generalization of the former results of P. Li and Y. Ding about the consistent large time behavior of heat kernel. More concretely, we have the following theorem.
\begin{theorem}\label{thm 11.0}
{Assume that $(M^n, g)$ is a complete manifold with cone structures at infinity, $y$ is some fixed point on $M^n$ and $n\geq 3$. Furthermore assume that for any $r> 0$, any two positive sequence $\{s_i\}$, $\{l_i\}$ with the following property: 
\begin{align}
\lim_{i\rightarrow \infty} s_i= \lim_{i\rightarrow \infty} l_i= \infty\ , \quad \lim_{i\rightarrow \infty} \frac{V_y(\sqrt{s_i}r)}{V_y(\sqrt{s_i})}= h(r)\ , \quad \lim_{i\rightarrow \infty} \frac{V_y(\sqrt{l_i}r)}{V_y(\sqrt{l_i})}= \tilde{h}(r) \label{vol}
\end{align}  
where $h(r)$, $\tilde{h}(r)$ are positive functions, the following equation holds:
\begin{align}
\frac{h''(r)}{h'(r)}= \frac{\tilde{h}''(r)}{\tilde{h}'(r)} \label{equ vol}
\end{align}
Then
\begin{align}
\lim_{t\rightarrow \infty} V_y(\sqrt{t})\cdot H(x, y, t)= p_{\infty}(y_{\infty}, y_{\infty}, 1) \label{11.0.2}
\end{align}
where $p_{\infty}$ is the heat kernel on any tangent cone at infinity of manifold $M^n$ with renormalized measure, and the value of the right hand side is consistent.
}
\end{theorem}

The concept of manifolds with cone structures at infinity will be defined in Section \ref{SECTION 08}. Especially, the manifolds with nonnegative Ricci curvature and maximal volume growth satisfy the assumptions in Theorem \ref{thm 11.0}, in fact $h(r)= \tilde{h}(r)= r^n$ in this case.
 
Furthermore, we construct the first example of manifold with $Rc\geq 0$, where the limit in Question \ref{ques 1} does not exist. More precisely, we have the following theorem.
\begin{theorem}\label{thm 1.3}
{There exists a complete Riemannian manifold $(M^8, g)$ with $Rc\geq 0$, such that on $(M^8, g)$,
\begin{align}
\varliminf_{t\rightarrow \infty}V(\sqrt{t})H(x, y, t)< \varlimsup_{t\rightarrow \infty}V(\sqrt{t})H(x, y, t) \nonumber 
\end{align}
}
\end{theorem}

Following Cheeger and Colding's strategy in Section $8$ of \cite{CC1}, we modify the examples there to construct our example. Note that not every two different tangent cones at infinity of manifold will give different values of $p_{\infty}(y, y ,1)$. The different renormalized measures on tangent cones at infinity of manifold are the key point to result in the inconsistent limit behavior of heat kernel. 

The organization of this paper is as the following. In Section \ref{SECTION 2}, we state some background facts about Gromov-Hausdorff convergence, which are needed for later sections. For this part, we mainly refer to \cite{CC}, \cite{CC1}, \cite{Gromov}. And we also review the results about the first order differentiation, Sobolev spaces and Laplacian operator on metric measure spaces, which were proved in \cite{Cheeger} and \cite{CC3}. 

In Section $3$, we proved a Harnack's convergence theorem in Gromov-Hausdorff topology (Theorem \ref{thm 4.2}), which roughly says that the limit (if it exists) of harmonic functions on manifolds, is a harmonic function on limit spaces under some gradient bounds assumption. Theorem \ref{thm 4.2} was originally due to Y. Ding (see Section $3$ of \cite{Ding}). For reader's convenience, we provide a detailed proof here.

In Section $4$, as in \cite{Ding}, combining with the well-known estimates of eigenvalues and eigenfunctions, the convergence of eigenvalues and eigenfunctions in Gromov-Hausdorff sense follows from the Harnack's convergence theorem proved in Section \ref{SECTION 3}.
 
In Section $5$, the heat equation on metric measure space $M_{\infty}$ is discussed. Using the theory of abstract Cauchy problem developed in \cite{LM}, we get the existence of the solutions of heat equation on $M_{\infty}$ as in \cite{Sturm2}. In addition, some mean value inequality of the heat equation solutions are obtained, whose proof imitates L. Saloff-Coste's argument on smooth manifolds (cf. see \cite{SB}).

In Section $6$, we follow closely the argument of K.-T. Sturm in \cite{Sturm2} (also see \cite{Sturm1}, \cite{Sturm3} and \cite{Sturm4}) and L. Saloff-Coste in \cite{SB} (also see \cite{Saloff}, \cite{Saloff1}) to prove the existence and Gaussian upper bound of heat kernel on metric measure space $(M_{\infty}, \rho_{\infty}, \nu_{\infty})$. We believe that some results in this section are well-known to experts in this field in more general context, but we provide the details here to make our argument self-contained.

In Section $7$, using the results established in the former sections, we manage to reduce the $\displaystyle \lim_{i\rightarrow \infty} V\big(\sqrt{t_i}\big) H(x, y, t_i)$ to the heat kernel value $p_{\infty}(y, y, 1)$ on $(M_{\infty}, \nu_{\infty})$, where $M_{\infty}$ is any tangent cone at infinity of complete manifold $M^n$ with $Rc\geq 0$ and $\nu_{\infty}$ is the renormalized measure on $M_{\infty}$.

In Section $8$, by the general reduction results obtained in Section $7$, the general criterion in Theorem \ref{thm 11.0} is given to determine whether the limit behavior of heat kernel is consistent. This general criterion includes the former related results of P. Li and Y. Ding as a special case. 

In Section $9$, using the generalized Hopf fibration of $\mathbb{S}^7$, we construct the example $(M^8, g)$ by modifying the metric on $\mathbb{R}^8$ step by step. When $M_{\infty}$ have cone structure $dr^2+ f(r)^2 dX$, one key point to get different heat kernel values $p_{\infty}(y, y, 1)$ on $(M_{\infty}, \nu_{\infty})$ is, to assure that (\ref{equ vol}) does not hold for two specially chosen positive sequences whose limits are infinity. The computation involved in the construction of this example is long but straightforward, we give the details for completeness.

Finally in Appendix \ref{App 1}, some $L^p$-convergence results in Gromov-Hausdorff sense are stated, and the proof of the Rellich-type compactness theorem is also provided for reader's convenience.

\section{Preliminaries on Cheeger-Colding's theory}\label{SECTION 2}
In this section we review some background material about Gromov-Hausdorff convergence and analysis on limit spaces, which were established in \cite{Gromov} and \cite{CC1}, \cite{CC2}, \cite{CC3}, \cite{Cheeger}. Especially, the doubling condition and local Poincar\'e inequality on limit spaces are showed. Also the existence of self-adjoint Laplacian operator on limit spaces is established. Those two results are used repeatedly through the whole paper.




Let $\big\{(M_i^n, y_i, \rho_i)\big\}$ be a sequence of pointed Riemannian manifolds, where $y_i\in M_i^n$ and $\rho_i$ is the metric on $M_i^n$. If $\big\{(M_i^n, y_i, \rho_i)\big\}$ converges to $(M_{\infty}, y_{\infty}, \rho_{\infty})$ in the Gromov-Hausdorff sense, we write $\displaystyle (M_i^n, y_i, \rho_i)\stackrel{d_{GH}}{\longrightarrow} (M_{\infty}, y_{\infty}, \rho_{\infty})$. See \cite{Gromov} for the definition and basic facts concerning Gromov-Hausdorff convergence.

Obviously if a sequence of pointed metric spaces converges to a pointed space $(X, p)$ in the Gromov-Hausdorff sense, it also converges to its completion. We will only consider complete metric spaces as Gromov-Hausdorff limits. Then, similarly to the case of ordinary convergence, a Gromov-Hausdorff limit of pointed spaces is essentially unique. For general background on metric space and length space, we refer the reader to \cite{BBI}. 

Let $\displaystyle (X_i, p_i)\stackrel{d_{GH}}{\longrightarrow} (X, p)$ where $X_i$ are length spaces and $X$ is a complete metric space, from Theorem $8.1.9$ in \cite{BBI}, $X$ is a complete length space.

From the above argument, we get that $\big(M_{\infty}, y_{\infty}, \rho_{\infty}, \nu_{\infty}\big)$ is a complete length space.

A metric space is said to be boundedly compact if all closed bounded sets in it are compact. By Exercise $8.1.8$ in \cite{BBI}, $(M_{\infty}, \rho_{\infty})$ is also boundedly compact.

We define the convergence concept for functions on manifolds $\{M_i^n\}$ as the following, it is so called "uniform convergence in Gromov-Hausdorff topology", for simplification, sometimes it is written as "uniform convergence in G-H topology".

\begin{definition}[Uniform Convergence in G-H topology]\label{def 1.1}
{Suppose 
\begin{align}
K_i\subset M_i^n\stackrel{d_{GH}}{\longrightarrow} K_{\infty}\subset M_{\infty} \nonumber
\end{align}
Assume that $\{f_i\}_{i= 1}^{\infty}$ are functions on $M_i^n$, $f_{\infty}$ is a function on $M_{\infty}$. and $\Phi_{i}: K_{\infty}\rightarrow K_i$ are $\epsilon_i$-Gromov-Hausdorff approximations, $\lim_{i\rightarrow \infty}\epsilon_i= 0$. If $f_i\circ \Phi_i$ converge to $f_{\infty}$ uniformly, we say that $f_i\rightarrow f_{\infty}$ uniformly over $K_i\stackrel{d_{GH}}{\longrightarrow} K_{\infty}$.
}
\end{definition}

As in Section $9$ of \cite{Cheeger}, we have the following definition.
\begin{definition}\label{def MGH}
{If $\nu_i$, $\nu_{\infty}$ are Borel regular measures on $M_i^n$, $M_{\infty}$, we say that $(M_i^n, y_i, \rho_i, \nu_i)$ converges to $(M_{\infty}, y_{\infty}, \rho_{\infty}, \nu_{\infty})$ \textbf{in the measured Gromov-Hausdorff sense}, if $(M_i^n, y_i, \rho_i)\stackrel{d_{GH}}{\longrightarrow} (M_{\infty}, y_{\infty}, \rho_{\infty})$, in addition, for any $x_i\rightarrow x_{\infty}$, ($x_i\in M_i^n$, $x_{\infty}\in M_{\infty}$), $r> 0$, we have 
\[\nu_{i}\Big(B_{i}(x_i, r)\Big)\rightarrow \nu_{\infty}\Big(B_{\infty} (x_{\infty}, r)\Big)\]
where $(M_{\infty}, \rho_{\infty})$ is a length space with length metric $\rho_{\infty}$, and  
\begin{align}
B_i(x_{i}, r)= \{z\in M_i^n|\  d_{\rho_i}(z, x_i)\leq r\}\ , \quad B_{\infty}(x_{\infty}, r) = \{z\in M_{\infty}|\  d_{\rho_{\infty}}(z, x_{\infty})\leq r\} \nonumber
\end{align}  
}
\end{definition}

In the rest of this section, we assume that $\{M_i^n\}$ is a sequence of complete noncompact manifolds with non-negative Ricci curvature, $\nu_i$ is the renormalized measure on $M_i^n$ defined as $\nu_i(A)= \frac{\mu_i(A)}{\mu_i(B_i(1))}$, where $\mu_i$ is the volume element determined by $\rho_i$. And $(M_i^n, y_i, \rho_i, \nu_i)$ converges to $(M_{\infty}, y_{\infty}, \rho_{\infty}, \nu_{\infty})$ in the measured Gromov-Hausdorff sense. Note from Theorem $1.6$ in \cite{CC1}, any sequence $(M_i, y_i, \rho_i)$ with $Rc\geq 0$, there is a subsequence, $(M_i^n, y_i, \rho_i, \nu_i)$, convergent to some $(M_{\infty}, y_{\infty}, \rho_{\infty}, \nu_{\infty})$ in the measured Gromov-Hausdorff sense.

Before discussing the analysis on $M_{\infty}$, we firstly consider the general metric measure space $(X, m)$, where $X$ is a metric space and $m$ is a Borel regular measure on $X$. Hence $(M_{\infty}, \rho_{\infty}, \nu_{\infty})$ is a special case of $(X, m)$. Fixed a set $A\subset X$, let $f$ be a function on $A$ with values in the extended real numbers. 

\begin{definition}\label{grad}
{An \textbf{upper gradient}, $g$, for $f$ is an extended real valued Borel function, $g: A\rightarrow [0, \infty]$, such that for all points, $z_1$, $z_2\in A$, and all continuous rectifiable curves, $\gamma: [0, l]\rightarrow A$, parameterized by arc-length, $s$, with $\gamma(0)= z_1$, $\gamma(l)= z_2$, we have
\begin{align}
|f(z_1)- f(z_2)|\leq \int_0^{l}g(\gamma(s)) ds \label{grad 1}
\end{align}
}
\end{definition}

Fix an open set $U\subset X$, and until further notice, write $L^p$ for $L^p(U)$. For $f\in L^p$, we set
\begin{align}
|f|_{1, p}\doteqdot |f|_{L^p}+ \inf_{\{g_i\}}\liminf_{i\rightarrow \infty} |g_i|_{L^p} \label{sobolev}
\end{align}
where the inf is taken over all sequences $\{g_i\}$, for which there exists a sequence, $f_i\stackrel{L^p}{\longrightarrow} f$, such that $g_i$ is an upper gradient for $f_i$, for all $i$.

\begin{definition}\label{SoboSpace}
{For $p\geq 1$, the Sobolev space $W^{1, p}(U)$ is the subspace of $L^p(U)$ consisting of functions, $f$, for which $|f|_{1, p}< \infty$, equipped with the norm $|\cdot|_{1, p}$.                                           
}      
\end{definition}

Let $0\rightarrow W^{1, p}\stackrel{\mathbf{i}}{\rightarrow} L^p$ denote the natural map, $U_{\eta}\subset U$ denote the set of points at distance $\geq \eta$ from $\partial U$. Let $\mathcal{K}(U)$ denote the subset of $W^{1, p}(U)$ consisting of those functions, $f$, for which there exists $\eta> 0$, such that $\mathbf{i}(f)$, the image of $f$, in $L^p(U)$, has a representative with support in $U_{\eta}$. 

\begin{definition}\label{SoboSpace2}
{The Sobolev space $W_{0}^{1, p}(U)\subset W^{1, p}(U)$, is the closure of the space $\mathcal{K}(U)$ in $W^{1, p}(U)$. 
}
\end{definition}

From Definition $2.8$, $2.9$ and Theorem $2.10$ in \cite{Cheeger}, we have the following theorem.
\begin{theorem}[Cheeger]\label{thm mgug}
{For all $1< p< \infty$ and $f\in W^{1, p}(U)$, there exists a unique $g_f\in L^p(U)$ (up to modification on subsets of measure zero) such that
\begin{align}
|f|_{1, p}= |f|_{L^p}+ |g_f|_{L^p} \label{mgug1}
\end{align}
and there exist sequences, $f_i\stackrel{L^p}{\rightarrow} f$, $g_i\stackrel{L^p}{\rightarrow} g_f$, where $g_i$ is an upper gradient for $f_i$, for all $i$.
}
\end{theorem}

$g_f$ is called a \textbf{minimal generalized upper gradient} for $f$, which may depend on the choice of $p$ and $U$.

When $p= 2$, the above Sobolev spaces become Hilbert spaces, we use the following notations:
\begin{align}
H^1\doteqdot W^{1, 2}\ , \quad H_0^1\doteqdot W_0^{1, 2} \nonumber
\end{align}

We define the following properties:


Property ($\mathscr{B}$)(\textbf{the doubling condition}): For all balls $B_{2r}(x)\subset X$, we have 
\begin{align}
m\Big(B_{2r}(x)\Big)\leq 2^{n}\cdot m\Big(B_r(x)\Big) \label{double}
\end{align}

Property ($\mathscr{C}$): There exists a constant $C= C(n)$ such that for all balls $B_{2r}(x)\subset X$, we have 
\begin{align}
\int_{B_r(x)} |f- f_{x, r}|^2 dm\leq Cr^2 \int_{B_{2r}(x)} |g_{f}|^2 dm \label{Poincare}
\end{align}
for all $f\in H^1\big(X, m\big)$, and
\begin{align}
f_{x, r}= \frac{1}{m\big(B_r(x)\big)} \int_{B_r(x)} f dm \label{ave int}
\end{align}



We have the following proposition about $(M_{\infty}, y_{\infty}, \rho_{\infty}, \nu_{\infty})$.
\begin{prop}\label{prop 3.9}
Property ($\mathscr{B}$), ($\mathscr{C}$) hold on $\big(M_{\infty}, \rho_{\infty}, y, \nu_{\infty}\big)$.
\end{prop}

\pf
{It follows from Volume Comparison Theorem that Property ($\mathscr{B}$) holds on $(M_i^n, y_i, \rho_i, \nu_i)$. 

By $Rc\geq 0$ on $M_i^n$, from Theorem $5.6.5$ in \cite{SB}, we have
\begin{align}
\int_{B_i(z, r)} |f- f_{z, r}|^{\frac{3}{2}} d\nu_i\leq C(n)r^{\frac{3}{2}} \int_{B_i(z, r)} |\nabla f|^{\frac{3}{2}} d\nu_i \ ,  \quad f\in H^1(M_i, \nu_i)\label{3.9.0}
\end{align}
Using H\"older inequality, we obtain that 
\begin{align}
\Big(\big|f- f_{z, r}\big|\Big)_{z, r}\leq C(n)r \Big[\big(|\nabla f|^{\frac{3}{2}}\big)_{z, r}\Big]^{\frac{2}{3}}\ ,  \quad f\in H^1(M_i, \nu_i) \label{3.9.1}
\end{align} 

By Theorem $9.6$ in \cite{Cheeger}, we get Property ($\mathscr{B}$) and the following inequality holds on $(M_{\infty}, y_{\infty}, \rho_{\infty}, \nu_{\infty})$:
\begin{align}
\Big(|f- f_{z, r}|\Big)_{z, r}\leq C(n)r\Big[\big(|g|^2\big)_{z, r}\Big]^{\frac{1}{2}}\label{3.9.2}
\end{align}
where $f\in H^1(M_{\infty}, \nu_{\infty})$ and $g$ is any upper gradient for $f$. 

Using Theorem \ref{thm mgug}, there exist sequences, $f_i\stackrel{L^2}{\rightarrow} f$, $g_i\stackrel{L^2}{\rightarrow} g_f$, and $g_i$ is an upper gradient for $f_i$. From (\ref{3.9.2}), we get
\begin{align}
\Big(|f_i- (f_i)_{z, r}|\Big)_{z, r}\leq C(n)r\Big[\big(|g_i|^2\big)_{z, r}\Big]^{\frac{1}{2}}\nonumber
\end{align}
taking $i\rightarrow \infty$ in the above inequality, we have
\begin{align}
\Big(|f- f_{z, r}|\Big)_{z, r}\leq C(n)r\Big[\big(|g_f|^2\big)_{z, r}\Big]^{\frac{1}{2}}\ ,  \quad f\in H^1(M_{\infty}, \nu_{\infty}) \label{3.9.2.0} 
\end{align}

From the argument in the beginning of Section \ref{SECTION 2}, we know that $(M_{\infty}, \rho_{\infty})$ is a complete boundedly compact length space. By Corollary $1$ in \cite{HK}, $B_{\infty}(z, r)$ satisfies the $C(\lambda, M)$ condition (defined in \cite{HK}) for $\lambda= 1$ and some independent constant $M$. Then we can use (\ref{3.9.2.0}) and Theorem $1$ in \cite{HK} to get
\begin{align}
\Big[\big(|f- f_{z, r}|^{2\chi}\big)_{z, r}\Big]^{\frac{1}{2\chi}} \leq \tau r\Big[\big(|g_f|^2\big)_{z, r}\Big]^{\frac{1}{2}}\ ,  \quad f\in H^1(M_{\infty}, \nu_{\infty}) \label{3.9.3}
\end{align}
where $\chi= \chi(n)> 1$, $\tau= \tau(n, \chi)> 0$ are some constants.

By (\ref{3.9.3}) and H\"older inequality, we conclude that 
\begin{align}
\int_{B_{\infty}(z, r)} \big|f- f_{z, r}\big|^2 d\nu_{\infty} \leq C(n)r^2 \int_{B_{\infty}(z, r)} |g_f|^2 d\nu_{\infty}\ ,  \quad f\in H^1(M_{\infty}, \nu_{\infty}) \label{3.9.4}
\end{align}
which implies Property ($\mathscr{C}$) on $(M_{\infty}, y_{\infty}, \rho_{\infty}, \nu_{\infty})$.
}
\qed

We have the following theorem about ``$df$":
\begin{theorem}[\cite{Cheeger}, \cite{CC3}]\label{thm 3.8}
{$f\in H^1(M_{\infty})$ $\Big(H_0^1(M_{\infty})\Big)$, if and only if there exists a sequence of Lipschitz functions (compactly supported Lipschitz functions) $f_i\stackrel{L^2}{\longrightarrow} f$ and $df_i\stackrel{L^2}{\longrightarrow} \omega$ for some $L^2$-section $\omega$ of $T^{*}M_{\infty}$, and $\omega$ is unique. 
}
\end{theorem}

\pf
{By Theorem $4.47$ in \cite{Cheeger} (also see Theorem $6.7$ in \cite{CC3}) and Proposition \ref{prop 3.9} above, we get our conclusion.
}
\qed

\begin{remark}\label{rem df}
{$\omega$ in Theorem \ref{thm 3.8} is called a strong $L^2$ exterior derivative of $f$ in \cite{CC3}, we can define $df\doteqdot \omega$ for $f\in H_0^1(M_{\infty})$, then $df$ is the $L^{\infty}$ section of $T^{*}M_{\infty}$ (the cotangent tensor bundle) determined by $f$, which is called the differential of $f$. From the Theorem above, it is well defined.
}
\end{remark}

We define 
\begin{align}
\mathscr{L}(U)&= \{f|\ f\ is\ Lipschitz\ function\ on\ U\} \nonumber \\
\mathscr{L}_c(U)&= \{f|\ f\ is\ compactly\ supported\ Lipschitz\ function\ on\ U\} \nonumber
\end{align}

From Theorem \ref{thm 3.8} above, we know that $\mathscr{L}_c(U)$ is dense in $H_0^1(U)$. We define $H_0^1(M_{\infty})$ as the closure of $\mathscr{L}_c(M_{\infty})$ in $H^1(M_{\infty})$. 

It is easy to see $C_c(U)$ is dense in $L^2(U)$, from the fact that any compactly supported continuous function can be uniformly approximated by compactly supported Lipschitz functions, we get that $\mathscr{L}_c(U)$ is dense in $L^2(U)$. Then $H_0^1(U)$ is also dense in $L^2(U)$.

Because the operator $d$ is well defined on $\mathscr{L}(M_{\infty})$, we can view the operator $d$ on $L^2(M_{\infty})$ as a densely defined unbounded operator. By Theorem \ref{thm 3.8}, this operator is closable as an operator on $L^2(M_{\infty})$. We have the existence of self-adjoint operator $\Delta_{\infty}$ on $M_{\infty}$ as the following.
\begin{theorem}[\cite{CC3}]\label{thm Laplacian} 
{The bilinear form $\displaystyle \int_{M_{\infty}} <df_1, df_2> d\nu_{\infty}$ is a densely defined, closed symmetric form on $L^2(M_{\infty})$. Hence, there is a unique self-adjoint operator, $\Delta_{\infty}$, (associated to the minimal closure), such that
\begin{align}
\int_{M_{\infty}} |df|^2 d\nu_{\infty}= \int_{M_{\infty}} <(-\Delta_{\infty})^{\frac{1}{2}} f , (-\Delta_{\infty})^{\frac{1}{2}} f> d\nu_{\infty}\ , \quad f\in H_0^1(M_{\infty}) \label{Lap 1}
\end{align}
}
\end{theorem}

\pf
{It follows from Theorem $2.23$ of \cite{Kato}.
}
\qed

\section{Harnack's convergence theorem in the Gromov-Hausdorff sense}\label{SECTION 3}
In this section, we will show that under uniform gradient bound assumption, the uniform limit of solutions, of Poisson equations on a sequence of convergent manifolds (in Gromov-Hausdorff sense), if it exists, will be the solution of Poisson equation on the limit space. The result of this section will only be needed in Section \ref{SECTION 04}.   

Compared with the classical Harnack's convergence theorem (cf. Theorem $2.9$ in \cite{GT}), which says that the limit of monotonic increasing bounded harmonic functions is still harmonic, where monotonicity is used to apply Harnack estimate on harmonic functions. With the uniform gradient bound assumption replacing Harnack estimate, one may think of our theorem (Theorem \ref{thm 4.2}) as Harnack's convergence theorem in the Gromov-Hausdorff sense, which is crucial in the proof of Proposition \ref{prop 3.7}.

On Riemannian manifold $(M_i^n, \rho_i, \nu_i)$, one solves the Poisson equation
\begin{equation}\nonumber
\left\{
\begin{array}{rl}
\Delta_{\rho_i} u= f  \\
u\Big|_{\partial B_i(x_i, r)}= h
\end{array} \right.
\end{equation}
for Lipschitz functions $f$, $h$ on $B_i(x_i, r)\subset M_i^n$. By the Dirichlet's principle, $u$ is the unique minimizer of the functional 
\begin{align}
I(u, \nu_i, x_i, r)= \int_{B_i(x_i, r)} \Big(\frac{1}{2}|\nabla u|^2+ fu \Big) d\nu_i \nonumber
\end{align}
within the space $\mathscr{H}_i\doteqdot h+ H_0^1\Big(B_i(x_i, r)\Big)$.

Similarly, for $(M_{\infty}, \rho_{\infty}, \nu_{\infty})$, by Theorem \ref{thm Laplacian}, the solution of the Poisson equation  
\begin{equation}\nonumber
\left\{
\begin{array}{rl}
\Delta_{\infty} u= f  \\
u\Big|_{\partial B_{\infty}(x_{\infty}, r)}= h
\end{array} \right.
\end{equation}
is the unique minimizer of the functional 
\begin{align}
I(u, \nu_{\infty}, x_{\infty}, r)= \int_{B_{\infty}(x_{\infty}, r)} \Big(\frac{1}{2}|d u|^2+ fu \Big) d\nu_{\infty} \nonumber
\end{align}
within the space $\mathscr{H}_{\infty}\doteqdot h+ H_0^1\Big(B_{\infty}(x_{\infty}, r)\Big)$.

The following theorem was originally proved by Y. Ding. We present a detailed proof here for completeness, which is loosely based on that in \cite{Ding}.

\begin{theorem}\label{thm 4.2}
{Suppose $u_i$, $f_i$ are $C^2$ functions over $B_{i}(x_i, 2r)\subset (M_i^n, y_i, \rho_i, \nu_i)$, where $B_{i}(x_i, 2r)= \{z\in M_i^n|\ d_{\rho_i}(z, x_i)\leq 2r\}$; $\Delta_{\rho_i} u_i= f_i$ on $B_{i}(x_i, 2r)$ and $r$ is some fixed positive constant.
Also assume $u_i\rightarrow u_{\infty}$, $f_i\rightarrow f_{\infty}$ uniformly over the sequence of converging balls $B_{i}(x_i, 2r)\rightarrow B_{\infty}(x_{\infty}, 2r)\subset (M_{\infty}, y_{\infty}, \rho_{\infty}, \nu_{\infty})$, and there exists $L> 0$ such that for any $i$:
\begin{align}
|\nabla u_i(x)|\leq L \ ,  \quad |\nabla f_i(x)|\leq L  \quad for \  x\in B_i(x_i, 2r)  \label{2.1.1}
\end{align}

Then
\begin{align}
\Delta_{\infty} u_{\infty}= f_{\infty}  \quad  on\  B_{\infty}\big(x_{\infty}, r\big) \label{4.2.1}
\end{align}
}
\end{theorem}

\pf
{To prove the theorem, we need the following lemma:
\begin{lemma}\label{lem 2.1}
{Let $u_{\infty}$, $f_{\infty}$ be as in Theorem \ref{thm 4.2}, then we have
\begin{align}
I(u_{\infty}, \nu_{\infty}, x_{\infty}, r)\leq \liminf_{i\rightarrow \infty} I(u_{i}, \nu_i, x_{i}, r) \label{2.1.2}
\end{align}
where 
\begin{align}
I(u_{\infty}, \nu_{\infty}, x_{\infty}, r)&= \int_{B_{\infty}(x_{\infty}, r)} \big(\frac{1}{2}|d u_{\infty}|^2+ f_{\infty}u_{\infty}\big)d\nu_{\infty} \nonumber \\
I(u_i, \nu_i, x_i, r)&= \int_{B_i(x_i, r)} \big(\frac{1}{2}|\nabla u_i|^2+ f_i u_i\big)d\nu_i \nonumber 
\end{align}
}
\end{lemma}

The proof of the Lemma is deferred to the end of this section.  We assume that  Lemma \ref{lem 2.1} holds, and prove the theorem by contradiction. Assume $\Delta_{\infty} u_{\infty}= f_{\infty}$ is not true over $B_{\infty}(x, s)\subset \subset B_{\infty}\big(x_{\infty}, r\big)$. 

By solving the Dirichlet problem on $B_{\infty}(x, s)$ (see Theorem $7.8$ and Remark $7.11$ in \cite{Cheeger}), we can find $\tilde{u}_{\infty}$ with the same boundary value as $u_{\infty}$ over $\partial B_{\infty}(x, s)$ and 
\begin{align}
I(\tilde{u}_{\infty}, \nu_{\infty}, x, s)< I(u_{\infty}, \nu_{\infty}, x, s) - 2\delta \label{4.2.2}
\end{align}
where $\delta> 0$ is some constant.

By Lemma \ref{lem 2.1}, assume that $x^{(i)}\rightarrow x$, then there exists $i_1> 0$, for $i> i_1$,
\begin{align}
I(u_{\infty}, \nu_{\infty}, x, s)\leq I(u_i, \nu_i, x^{(i)}, s)+ \delta \label{4.2.3}
\end{align}

By Lemma $10.7$ in \cite{Cheeger}, we can find a sequence of Lipschitz functions $\tilde{u}_i: B_i(x^{(i)}, s)\rightarrow \mathbb{R}$, such that $\tilde{u}_i$ converges uniformly to $\tilde{u}_{\infty}$ and 
\begin{align}
\displaystyle \varlimsup_{i\rightarrow \infty}  \int_{B_i(x^{(i)}, s)} |\nabla \tilde{u}_i|^2 d\nu_i \leq \int_{B_{\infty}(x, s)} |d \tilde{u}_{\infty}|^2 d\nu_{\infty} \nonumber
\end{align}
Hence there exists $i_2> 0$, for $i> i_2$, 
\begin{align}
I(\tilde{u}_i,  \nu_i, x^{(i)}, s)< I(\tilde{u}_{\infty}, \nu_{\infty}, x, s)+ \frac{1}{2}\delta \label{4.2.4}
\end{align}

By (\ref{4.2.2}), (\ref{4.2.3}) and (\ref{4.2.4}), we get that for $i> i_0$, where $i_0= \max\{i_ 1, i_ 2\}$,
\begin{align}
I(\tilde{u}_i, \nu_i, x^{(i)}, s)< I(u_i, \nu_i,  x^{(i)}, s)- \frac{1}{2}\delta \label{4.2.5}
\end{align}

When $i> i_0$, solve the following Dirichlet problem:
\begin{equation}\nonumber
\left\{
\begin{array}{rl}
\Delta \hat{u}_i= f_i &\quad \quad on\ B_i(x^{(i)}, s)\\
\hat{u}_i= \tilde{u}_i &\quad \quad on \ \partial B_i(x^{(i)}, s) \\
\end{array} \right.
\end{equation}
then by Dirichlet principle and (\ref{4.2.5}), we get that
\begin{align}
I(\hat{u}_i, \nu_i, x^{(i)}, s)\leq I(\tilde{u}_i, \nu_i, x^{(i)}, s)< I(u_i, \nu_i, x^{(i)}, s)- \frac{1}{2}\delta \label{4.2.6}
\end{align}

Note in fact we have
\begin{equation}\nonumber
\left\{
\begin{array}{rl}
\Delta (\hat{u}_i- u_i)= 0 &\quad \quad on\ B_i(x^{(i)}, s)\\
(\hat{u}_i- u_i)= (\tilde{u}_i- u_i) &\quad \quad on \ \partial B_i(x^{(i)}, s) \\
\end{array} \right.
\end{equation}
and 
\begin{align}
\lim_{i\rightarrow \infty} \sup_{\partial B_i(x^{(i)}, s)} |\tilde{u}_i- u_i|= \sup_{\partial B_{\infty}(x, s)} |\tilde{u}_{\infty}- u_{\infty}|= 0 \nonumber 
\end{align}

By maximum principle, we get
\begin{align}
\lim_{i\rightarrow \infty} \sup_{z\in B_i(x^{(i)}, s)} |\big(\hat{u}_i- u_i\big)(z)|\leq \lim_{i\rightarrow \infty} \sup_{z\in \partial B_i(x^{(i)}, s)} |\big(\tilde{u}_i- u_i\big)(z)|= 0 \label{4.2.7}
\end{align}

From (\ref{4.2.6}) and (\ref{4.2.7}), there exists $i_3> 0$, such that for $i> i_3$, 
\begin{align}
\frac{1}{2}\int_{B_i(x^{(i)}, s)} |\nabla \hat{u}_i|^2 d\nu_i< \frac{1}{2} \int_{B_i(x^{(i)}, s)} |\nabla u_i|^2 d\nu_i- \frac{1}{4}\delta \nonumber
\end{align}

By $|\nabla u_i|\leq L$ in (\ref{2.1.1}) and volume convergence of $B_{i}(x^{(i)}, s)$, there exists $i_4> 0$ and $s_1\in (0, s)$, such that for $i> i_4$,
\begin{align}
\int_{B_i(x^{(i)}, s)\backslash B_i(x^{(i)}, s_1)} |\nabla u_i|^2 d\nu_i< \frac{1}{100}\delta \nonumber
\end{align}
hence for $i> i_4$, we have
\begin{align}
\int_{B_i(x^{(i)}, s)} |\nabla \hat{u}_i|^2 d\nu_i< \int_{B_i(x^{(i)}, s_1)} |\nabla u_i|^2 d\nu_i- \frac{1}{4}\delta \label{4.2.8}
\end{align}

On $B_{i}(x^{(i)}, s_1)\subset\subset B_i(x^{(i)}, s)$, from Cheng-Yau's gradient estimate (also see Lemma \ref{lem 3.6} later), we get 
\begin{align}
\sup_{B_i(x^{(i)}, s_1)}|\nabla \hat{u}_i- \nabla u_i|\leq \frac{C(n)}{s- s_1}\sup_{B_{i}(x^{(i)}, s)} |\hat{u}_i- u_i| \label{4.2.9}
\end{align}

From (\ref{4.2.7}), (\ref{4.2.9}) and $|\nabla u_i|\leq L$, there exists $i_5> 0$, for $i> i_5$, 
\begin{align}
\int_{B_{i}(x^{(i)}, s_1)} |\nabla u_i|^2 d\nu_i- \int_{B_{i}(x^{(i)}, s_1)} |\nabla \hat{u}_i|^2 \leq \frac{1}{100}\delta \label{4.2.10}
\end{align}

From (\ref{4.2.8}) and (\ref{4.2.10}), we get
\begin{align}
\int_{B_{i}(x^{(i)}, s)\backslash B_i(x^{(i)}, s_1)} |\nabla \hat{u}_i|^2< -\frac{1}{8}\delta \nonumber
\end{align}

That is contradiction, the theorem is proved.
}
\qed

\bigskip
{\it \textbf{Proof of Lemma \ref{lem 2.1}}:}~
{Recall the Bochner formula:
\begin{align}
\frac{1}{2}\Delta \Big(|\nabla u_i|^2\Big)= \Big|\nabla^2 u_i\Big|^2+ <\nabla \Delta u_i, \nabla u_i>+ Rc(\nabla u_i, \nabla u_i) \label{2.1.3}
\end{align}

Multiply by a cut-off function $\phi$ with $supp (\phi)\subset B_i(x_i, 2r)$, $\phi|_{B_i(x_i, \frac{3}{2}r)}= 1$, $|\Delta \phi|\leq C(n, r)$, $\frac{|\nabla \phi|^2}{\phi}\leq C(n, r)$ (see Theorem $6.33$ of \cite{CC}):
\begin{align}
\frac{1}{2}\phi \Delta \Big(|\nabla u_i|^2\Big)= \phi \Big|\nabla^2 u_i\Big|^2+ \phi Rc(\nabla u_i, \nabla u_i)+ \phi <\nabla \Delta u_i, \nabla u_i> \label{2.1.4}
\end{align}

Integration by parts, using $Rc\geq 0$, we get
\begin{align}
&\frac{1}{2}\int_{B_i(x_i, 2r)} |\nabla u_i|^2 \Delta \phi d\nu_i \geq \int_{B_i(x_i, 2r)} \Big[\phi \Big|\nabla^2 u_i\Big|^2- \phi |\Delta u_i|^2- \Delta u_i \Big(\nabla \phi \cdot \nabla u_i\Big)\Big] d\nu_i \nonumber \\
&\quad \quad \geq \int_{B_i(x_i, 2r)} \Big[\phi \Big|\nabla^2 u_i\Big|^2- \frac{3}{2} \phi |\Delta u_i|^2- \frac{|\nabla \phi|^2}{2\phi} \big|\nabla u_i\big|^2 \Big]d\nu_i 
\nonumber \\
&\quad \quad \geq \int_{B_i(x_i, 2r)} \Big[\phi \Big|\nabla^2 u_i\Big|^2- C(n, r) |\nabla u_i|^2- \frac{3}{2}\phi |f_i|^2 \Big] d\nu_i \nonumber
\end{align}

Hence when $i$ is big enough,
\begin{align}
&\int_{B_i(x_i, 2r)} \phi \Big|\nabla^2 u_i \Big|^2 d\nu_i\leq C(n, r)\int_{B_i(x_i, 2r)} |\nabla u_i|^2 \Big(|\Delta \phi|+ 1\Big) d\nu_i  \nonumber \\
&\quad \quad \quad \quad \quad \quad \quad \quad \quad  + \frac{3}{2}\int_{B_{\infty}(x_{\infty}, 2r)} |f_{\infty}|^2 d\nu_{\infty}+ 1\nonumber \\
&\quad \leq C(n, r)L\cdot \Big[\nu_{\infty} \big(B_{\infty}(x_{\infty}, 2r)\big)+ 1\Big] + \frac{3}{2}\int_{B_{\infty}(x_{\infty}, 2r)} |f_{\infty}|^2 d\nu_{\infty}+ 1 \nonumber
\end{align}
We get a uniform upper bound of $\int_{B_i(x_i, \frac{3}{2}r)} \Big|\nabla^2 u_i \Big|^2 d\nu_i$.

By Theorem \ref{thm 2.3} in the Appendix, we can get that some subsequence of $|\nabla u_i|$ converges to a function $\Gamma$ on $B_{\infty}(x_{\infty}, r)$ in $L^2\Big(B_{\infty}(x_{\infty}, r), \nu_{\infty}\Big)$, from (\ref{2.1.1}) we also know that $\Gamma\in L^{\infty}\Big(B_{\infty}(x_{\infty}, r), \nu_{\infty}\Big)$. By Lusin's theorem for general topological spaces with measure and $\Gamma\in L^2\Big(B_{\infty}(x_{\infty}, r), \mu_{\infty}\Big)$, for any $\epsilon> 0$, there exists $K_{\epsilon}\subset\subset B_{\infty}\big(x_{\infty}, r\big)$ and $\nu_{\infty}\Big(B_{\infty}(x_{\infty}, r)\backslash K_{\epsilon}\Big)< \epsilon$, $\Gamma$ is continuous on $K_{\epsilon}$, note $K_{\epsilon}$ is $\nu_{\infty}$-measurable.

Note $\nu_{\infty}$ satisfies the doubling condition, which implies the Vitali Covering Theorem \Big(see Chapter $2$ of \cite{Ma}\Big), hence the Lebesgue Differentiation Theorem holds for measure $\nu_{\infty}$. Then 
\begin{align}
\lim_{s\rightarrow 0}\frac{\nu_{\infty}\Big(B_{\infty}(x, s)\cap K_{\epsilon}\Big)}{\nu_{\infty}\Big(B_{\infty}(x, s)\Big)}= 1 \,  \quad \quad \quad \nu_{\infty}\ a.e.\ x\in K_{\epsilon} \label{2.1.5} 
\end{align}

For $x\in K_{\epsilon}$ satisfying (\ref{2.1.5}), we will show
\begin{align}
\Big|d u_{\infty} (x)\Big|\leq \Gamma (x) \label{2.1.6}
\end{align}

Finally for $\displaystyle x\in \mathop{\cup}_{i= 1}^{\infty} K_{2^{-i}}$, (\ref{2.1.6}) is valid. Hence for $\nu_{\infty}$ a.e. $x\in B_{\infty}(x_{\infty}, r)$, (\ref{2.1.6}) is valid, which implies (\ref{2.1.2}) holds.

To prove (\ref{2.1.6}), it is enough to prove that for any $\delta> 0$, there exists $1> \epsilon(\delta)> 0$, when $d_{\rho_{\infty}}(y, x)< \epsilon(\delta)$, the following holds:
\begin{align}
\Big|u_{\infty}(x)- u_{\infty}(y)\Big|\leq d_{\rho_{\infty}}(y, x)\Big(\Gamma(x)+ 7\delta\Big) \label{2.1.7}
\end{align}

By contradiction. Then there is $1> \delta_{0}> 0$, $\{y_i\}_{i= 1}^{\infty}$, $y_i\in B_{\infty}(x_{\infty}, r)$, such that $d_{\rho_{\infty}}(y_i, x)= \ell_{i}\rightarrow 0$, and
\begin{align}
\Big|u_{\infty}(x)- u_{\infty}(y_i)\Big|> d_{\rho_{\infty}}(y_i, x)\Big[\Gamma(x)+ 7\delta_{0} \Big] \label{2.1.8}
\end{align}

Then for $z\in B_{\infty}\Big(x, \frac{\ell_i \delta_0}{L}\Big)$, $y\in B_{\infty}\Big(y_i, \frac{\ell_i \delta_0}{L}\Big)$, we have
\begin{align}
\Big|u_{\infty}(z)- u_{\infty}(y)\Big|&\geq \Big|u_{\infty}(x)- u_{\infty}(y_i)\Big|- \Big|u_{\infty}(z)- u_{\infty}(x)\Big|- \Big|u_{\infty}(y_i)- u_{\infty}(y)\Big| \nonumber\\
&> \ell_i \Big[\Gamma(x)+ 7\delta_0\Big]- L\cdot d_{\rho_{\infty}}(y_i, y)- L\cdot d_{\rho_{\infty}}(z, x) \nonumber\\
&\geq \ell_{i} \Big[\Gamma(x)+ 5\delta_0\Big] \label{2.1.9}
\end{align}

Pick $\tilde{x}_j$, $y_{j, i}\in M_j^{n}$, $\tilde{x}_j\rightarrow x$, $y_{j, i}\rightarrow y_i$, and $d(\tilde{x}_j, y_{j, i})= d(x, y_i)$. When $j$ is big enough, for all $z_j\in B_{j}\Big(\tilde{x}_j, \frac{\ell_i \delta_0}{L}\Big)$, $\tilde{y}_j\in B_j\Big(y_{j, i}, \frac{\ell_i \delta_0}{L}\Big)$ and all minimal geodesic $\gamma_j$ connecting $z_j$, $\tilde{y}_j$, by (\ref{2.1.9}), we have
\begin{align}
\int_{\gamma_j} |\nabla u_j| d\rho_j \geq \ell_i \Big[\Gamma(x)+ 4\delta_0 \Big] \label{2.1.10} 
\end{align}

Since $|\nabla u_j|\leq L$, a simple computation shows along every $\gamma_j$,
\begin{align}
|\nabla u_j|> \Gamma (x)+ 2\delta_0 \label{2.1.11}
\end{align}
on a subset of $\gamma_j$, which has $1$-dim Hausdorff measure at least 
$\frac{2\delta_0\ell_i} {L}$.

By $Rc\geq 0$ and Theorem $2.11$ in \cite{CC}, we get that the global segment inequality holds on $\Big(M_j^n, \rho_j, y, \nu_j \Big)$:
\begin{align}
\int_{A_1\times A_2} \Big(\int_0^{d_{\rho_j}(p, q)} e\big(\gamma_{p, q}(s)\big) ds\Big) d p d q\leq C(n) D \Big[\nu_j(A_1)+ \nu_j(A_2)\Big]\cdot \Big(\int_{W} e d\nu_{j}\Big)  \label{2.1.12}
\end{align}
where $e$ is any nonnegative integrable function on $W\subset M_j^n$, and $\gamma_{p, q}$ is a minimal geodesic from $p$ to $q$, 
\begin{align}
D\doteqdot \max_{p\in A_1, q\in A_2} d_{\rho_j}(p, q)\ , \quad A_1, A_2\subset M_j^n \ , \quad \mathop{\cup}_{p, q}\gamma_{p, q}\subset W \nonumber
\end{align}

Choose $A_1= B_j\Big(\tilde{x}_j, \frac{\ell_i \delta_0}{L}\Big)$, $A_2= B_j\Big(y_{j,i}, \frac{\ell_i \delta_0}{L}\Big)$ and $e= \chi_{E_j^i}$ in (\ref{2.1.12}), where 
\begin{align}
E_j^i\doteqdot \Big\{z|\ z\in B_j\Big(\tilde{x}_j, \ell_i\big(1+ \frac{\delta_0}{L}\big)\Big),\ |\nabla u_j (z)|> \Gamma(x)+ 2\delta_0 \Big\} \nonumber
\end{align}
then we get
\begin{align}
&\nu_j(E_j^i)\cdot C(n)\Big[1+ \frac{\delta_0}{L}\Big]\ell_i\Big[\nu_j\Big(B_j(\tilde{x}_j, \frac{\ell_i \delta_0}{L})\Big)+ \nu_j\Big(B_j(y_{j, i}, \frac{\ell_i \delta_0}{L})\Big)\Big] \nonumber \\
&\geq \frac{2\delta_0 \ell_i}{L}\cdot \nu_j\Big(B_j(\tilde{x}_j, \frac{\ell_i \delta_0}{L})\Big)\cdot \nu_j\Big(B_j(y_{j, i}, \frac{\ell_i \delta_0}{L})\Big) \nonumber
\end{align}
Using the Bishop-Gromov volume comparison theorem, we get that for any $i$, if $j$ big enough, 
\begin{align}
\frac{\nu_j(E_j^i)}{\nu_j\Big(B_j\big(\tilde{x}_j, \ell_i (1+ \frac{\delta_0}{L})\big)\Big)}\geq C(\delta_0, L, \Gamma(x), n) \ ,  \label{2.1.13}
\end{align}

From (\ref{2.1.13}), we obtain that there exists 
\begin{align}
\mathscr{C}_i\subset \mathscr{B}_i\doteqdot B_{\infty}(x, \ell_i\big(1+ \frac{\delta_0}{L}\big)) \nonumber
\end{align}
such that $\nu_{\infty}(\mathscr{C}_i)\geq \delta_1 \nu_{\infty}(\mathscr{B}_i)$, where $\delta_1= \frac{1}{2}C(\delta_0, L, \Gamma(x), n)$, and 
\begin{align}
F_j^i\subset E_j^i\ , \quad F_j^i \stackrel{d_{GH}}{\longrightarrow} \mathscr{C}_i\ as\ j\rightarrow \infty \nonumber 
\end{align}

For fixed $i$, we further assume $\varphi_j: F_j^i\rightarrow \mathscr{C}_i$ is a measure approximation and an $\epsilon_j$-Gromov-Hausdorff approximation for some $\epsilon_j\rightarrow 0$.

Let $\tau_1= \frac{\delta_1}{10}\nu_{\infty}(\mathscr{B}_i)$, $\tau_2= \frac{\delta_1\delta_0^2}{40}\nu_{\infty}(\mathscr{B}_i)$.

Let $h_j= |\nabla u_j|$, note that $h_j$ converges to $\Gamma$ in $L^2$ on $B_{\infty}(x_{\infty}, r)$. By Definition \ref{def Lp}, on $\mathscr{C}_i\subset B_{\infty}(x_{\infty}, r)$ there exists $h_{\infty}^{(k)}: \mathscr{C}_i\rightarrow \mathbb{R}$, such that
\begin{align}
\lim_{k\rightarrow \infty} \int_{\mathscr{C}_i} |h_{\infty}^{(k)}- \Gamma|^2 d\nu_{\infty}= 0 \label{2.13.1} 
\end{align}
and 
\begin{align}
\lim_{k\rightarrow \infty}\varlimsup_{j\rightarrow \infty} \int_{F_j^i} |h_j- h_{\infty}^{(k)}\circ\varphi_j|^2 d\nu_j= 0 \label{2.13.2}
\end{align}

For $\tau_1$, from (\ref{2.13.1}) and Egoroff's Theorem, there exists $A\subset \mathscr{C}_i$, such that $\nu_{\infty}(A)< \tau_1$, and on $\mathscr{C}_i- A$, $h_{\infty}^{(k)}\rightarrow \Gamma$ uniformly.

Note there exists $C_0> 0$, such that $\nu_j(B_j\big(\tilde{x}_j, \ell_i (1+ \frac{\delta_0}{L})\big))\leq C_0$ for any $i, j$. And there exists $k_1> 0$, if $k> k_1$, 
\begin{align}
|h_{\infty}^{(k)}- \Gamma|\leq \sqrt{\frac{\tau_2}{C_0}} \quad\quad \quad on \ \mathscr{C}_i- A \label{2.13.3}
\end{align}

For $\tau_2> 0$, from (\ref{2.13.2}), there exists $k_2> k_1> 0$, if $k\geq k_2$, 
\begin{align}
\varlimsup_{j\rightarrow \infty} \int_{F_j^i} |h_j- h_{\infty}^{(k)}\circ \varphi_j|^2 d\nu_j\leq \frac{\tau_2}{2} \nonumber
\end{align}
hence, there exists $j_1> 0$, if $j> j_1$, then
\begin{align}
\int_{F_j^i} |h_j- h_{\infty}^{(k_2)}\circ \varphi_j|^2 d\nu_j< \tau_2 \label{2.13.4}
\end{align}

Let $Q_j^i= F_j^i- \varphi_j^{-1}(A)$, then when $j> j_1$,
\begin{align}
\int_{Q_j^i} |h_j- \Gamma\circ\varphi_j|^2 d\nu_j&\leq 2\Big[\int_{Q_j^i} |h_{\infty}^{(k_2)}\circ \varphi_j- \Gamma\circ \varphi_j|^2+ \int_{Q_j^i} |h_j- h_{\infty}^{(k_2)}\circ\varphi_j|^2 \Big] \nonumber \\
&\leq  4\tau_2 \label{2.13.5}
\end{align}
the last inequality above follows from (\ref{2.13.3}) and (\ref{2.13.4}).

Define 
\begin{align}
\mathscr{W}_{\infty}= \{z|\ \Gamma(z)\leq \Gamma(x)+ \delta_0\ , \ z\in \mathscr{C}_i- A\} \nonumber
\end{align}
and 
\begin{align}
\mathscr{W}_j= \varphi_j^{-1}(\mathscr{W}_{\infty})\subset F_j^i- \varphi_j^{-1}(A)= Q_j^i \nonumber
\end{align}
hence on $\mathscr{W}_j$, $h_j(z)> \Gamma(x)+ 2\delta_0$, and $\big(\Gamma\circ \varphi_j\big)(z)\leq \Gamma(x)+ \delta_0$, we get 
\begin{align}
\int_{\mathscr{W}_j} |h_j- \Gamma\circ\varphi_j|^2\geq \int_{\mathscr{W}_j} \delta_0^2= \delta_0^2 \nu_j(\mathscr{W}_j) \label{2.13.6}
\end{align}

From (\ref{2.13.5}) and (\ref{2.13.6}),
\begin{align}
\nu_j(\mathscr{W}_j)\leq \frac{4\tau_2}{\delta_0^2}= \frac{\delta_1}{10}\nu_{\infty}(\mathscr{B}_i) \nonumber
\end{align}

Hence 
\begin{align}
\nu_{\infty}(\mathscr{W}_{\infty})= \lim_{j\rightarrow \infty} \nu_j(\mathscr{W}_j)\leq \frac{\delta_1}{10}\nu_{\infty}(\mathscr{B}_i) \nonumber
\end{align}

Define $\mathscr{A}_i= \{z\in \mathscr{B}_i|\ \Gamma(z)> \Gamma(x)+ \delta_0 \}$, note that on $\mathscr{C}_i- A- \mathscr{W}_{\infty}\subset \mathscr{B}_i$, $\Gamma(z)> \Gamma(x)+ \delta_0$, hence 
\begin{align}
\nu_{\infty}(\mathscr{A}_i)\geq \nu_{\infty}(\mathscr{C}_i)- \nu_{\infty}(A)- \nu_{\infty}(\mathscr{W}_{\infty})\geq \frac{4}{5}\delta_1\nu_{\infty}(\mathscr{B}_i) \label{2.14}
\end{align}


Note $\delta_1= \frac{1}{2}C(\delta_0, L, \Gamma(x), n)$, we get 
\begin{align}
\frac{\nu_{\infty}(\mathscr{A}_i)}{\nu_{\infty}\big(\mathscr{B}_i\big)}\geq C(\delta_0, L, \Gamma(x), n)> 0 \label{2.1.14}
\end{align}
where $C(\delta_0, L, \Gamma(x), n)$ in different lines may be different.

Now we have 
\begin{align}
0&< C(\delta_0, L, \Gamma(x), n)\leq \frac{\nu_{\infty}(\mathscr{A}_i)}{\nu_{\infty}\big(\mathscr{B}_i\big)}= \frac{\nu_{\infty}(\mathscr{A}_i\cap K_{\epsilon})+ \nu_{\infty}(\mathscr{A}_i\backslash K_{\epsilon})}{\nu_{\infty}(\mathscr{B}_i)} \nonumber \\
&\leq \frac{\nu_{\infty}(\mathscr{B}_i\backslash K_{\epsilon})}{\nu_{\infty}(\mathscr{B}_i)} + \frac{\nu_{\infty}(\mathscr{A}_i\cap K_{\epsilon})}{\nu_{\infty}(\mathscr{B}_i)}= (I)_i+ (II)_i \label{2.1.15}
\end{align}

From (\ref{2.1.5}) and the choice of $x$, we get $\lim_{i\rightarrow \infty} (I)_i= 0$. Because $\Gamma$ is continuous on $K_{\epsilon}$, it is easy to see that $(II)_i= 0$ when $i$ is big enough. We take $i\rightarrow \infty$ in (\ref{2.1.15}), it is contradiction. Hence (\ref{2.1.7}) holds for any $\delta> 0$, (\ref{2.1.6}) holds $\nu_{\infty}$ a.e. $B_{\infty}(x_\infty, r)$. We are done.
}
\qed

\section{The convergence of eigenfunctions in the Gromov-Hausdorff sense}\label{SECTION 04}

In this section, we will show that the eigenvalues, eigenfunctions on the convergent sequence of manifolds converge (subsequentially) to eigenvalues, eigenfunctions on limit space under Gromov-Hausdorff convergence. The main tools are eigenvalue and eigenfunction estimates obtained by P. Li, S-Y. Cheng, S-T. Yau and Harnack's convergence theorem in the Gromov-Hausdorff sense (Theorem \ref{thm 4.2}). 

Write $\lambda_{j, i}^{(R)}$ for the $j$-th Dirichlet eigenvalue over $B_i(R)\subset (M_i, y, \rho_i, \nu_i)$, $\phi_{j, i}^{(R)}$ is the corresponding eigenfunction satisfying the following:
\begin{equation}\label{3.4}
\left\{
\begin{array}{rl}
\Delta_{\rho_i} \phi_{j, i}^{(R)}&= \lambda_{j, i}^{(R)}\phi_{j, i}^{(R)} \quad \quad on \ B_i(R) \\
\phi_{j, i}^{(R)}(x)&= 0 \quad \quad \quad \quad on \ \partial B_i(R) 
\end{array} \right.
\end{equation}
and $\int_{B_i(R)} \phi_{j, i}^{(R)}\cdot \phi_{k, i}^{(R)} d\nu_i = \delta_{jk}$, where $\Delta_{\rho_i}$ is the Laplace operator with respect to the metric $\rho_i$.

From Theorem $3.1$ in \cite{Saloff}, for any $f\in H_0^{1}(B_i(R))$, we get  
\begin{align}
\Big[\int_{B_i(R)} |f|^{\frac{2n}{n- 2}} d\mu_i \Big]^{\frac{n- 2}{n}} \leq C(n)\frac{R^2}{\mu_i\Big(B_i(R)\Big)^{\frac{2}{n}}}\cdot \Big[\int_{B_i(R)} \big(|\nabla f|^2+ R^{-2} f^2 \big) d\mu_i\Big] \label{3.5.2}
\end{align}

Using Corollary $1.1$ in \cite{LS}, 
\begin{align}
R^{-2}\int_{B_i(R)} f^2 d\mu_i \leq C(n) \int_{B_i(R)} |\nabla f|^2 d\mu_i \label{3.5.3}
\end{align}

By (\ref{3.5.2}) and (\ref{3.5.3}), we have
\begin{align}
\int_{B_i(R)} |\nabla f|^2 d\mu_i &\geq C(n)\mu_i\Big(B_i(R)\Big)^{\frac{2}{n}} R^{-2}\cdot \Big[\int_{B_i(R)} |f|^{\frac{2n}{n- 2}} d\mu_i \Big]^{\frac{n- 2}{n}} \nonumber \\
&= C_{\mathscr{SD}} \Big[\int_{B_i(R)} |f|^{\frac{2n}{n- 2}} d\mu_i \Big]^{\frac{n- 2}{n}} \label{3.5.4}
\end{align}
where 
\begin{align}
C_{\mathscr{SD}}\doteqdot C(n)\mu_i\Big(B_i(R)\Big)^{\frac{2}{n}} R^{-2} \label{def C(SD)}
\end{align}

\begin{lemma}\label{lem 3.3}
{There exists a constant $C(n)$ such that
\begin{align}
C(n)^{-1}\cdot R^{-2}\cdot j^{\frac{1}{n}}  \leq \lambda_{j, i}^{(R)}\leq C(n)\cdot R^{-2}\cdot j^2 \label{3.5}
\end{align}
}
\end{lemma}

\pf
{Define $\displaystyle C_1(n)\doteqdot \sum_{\ell= 0}^{\infty}\frac{1}{2\beta^{\ell}- 1}$, where $\displaystyle \beta= \frac{n}{n- 2}$. Then we have $\frac{n}{4}\leq C_1(n)\leq \frac{n}{2}$.
By the argument of ($10.9$) in \cite{Libook} (also see \cite{LiSobolev}), we get the lower bound of $\lambda_{j, i}^{(R)}$ as the following:
\begin{align}
\lambda_{j, i}^{(R)}\geq C(n) j^{\frac{1}{2C_1(n)}} C_{\mathscr{SD}} \cdot \mu_i\big(B_i(R)\big)^{-\frac{2}{n}} \label{3.3.1}
\end{align}
combining with the definition of $C_{\mathscr{SD}}$ in (\ref{def C(SD)}), we have 
\begin{align}
\lambda_{j, i}^{(R)}\geq C(n) j^{\frac{1}{2C_1(n)}} R^{-2}\geq C(n)\cdot R^{-2}\cdot j^{\frac{1}{n}}  \label{3.3.2}
\end{align}
By the similar argument of Theorem $2$ on page $105$ of \cite{SY} (also see \cite{Cheng}), we get the upper bound of $\lambda_{j, i}^{(R)}$.
}
\qed

The following lemma is standard, for completeness, we provide the proof following the argument of Theorem $10.1$ in \cite{Libook}.
\begin{lemma}\label{lem 3.5}
{If $R> 2$, we have 
\begin{align}
\Vert \phi_{j, i}^{(R)}\Vert_{L^{\infty}(\nu_i)} \leq C(n, R)j^{\frac{n}{2}} \label{3.5.1}
\end{align}
where $\Vert \cdot \Vert_{L^{k}(\nu_i)}$ denotes the $L^k$ norm with respect to the measure $\nu_i$.
}
\end{lemma}

\pf
{We observe that for a $C^{\infty}$ function $u$, from Lemma $7.6$ and Lemma $7.7$ in \cite{GT}, $|\nabla u|^2= |\nabla |u||^2$ for $\mu_i$ a.e. $x$. The identities 
\begin{align}
\Delta (u^2)= 2u\Delta u+ 2|\nabla u|^2 \nonumber
\end{align}
and
\begin{align}
\Delta (|u|^2)= 2|u|\Delta |u|+ 2|\nabla |u||^2 \nonumber
\end{align}
imply $u\Delta u= |u|\Delta |u|$ a.e. Hence we have
\begin{align}
|\phi_{j, i}^{(R)}|\Delta_{\rho_i}|\phi_{j, i}^{(R)}|= \phi_{j, i}^{(R)} \Delta_{\rho_i} \phi_{j, i}^{(R)}= -\lambda_{j, i}^{(R)} |\phi_{j, i}^{(R)}|^2 \label{3.5.0}
\end{align}
For any constant $k\geq 2$, by (\ref{3.5.0}), (\ref{3.5.4}) and integration by parts,
\begin{align}
\int_{B_i(R)} \Big|\phi_{j, i}^{(R)}\Big|^k d\mu_i &= -\frac{1}{\lambda_{j, i}^{(R)}} \int_{B_i(R)} \Big|\phi_{j, i}^{(R)}\Big|^{k -1} \cdot \Delta_{\rho_i} \Big|\phi_{j, i}^{(R)}\Big| d\mu_i \nonumber \\
&= \frac{4(k- 1)}{\lambda_{j, i}^{(R)}\cdot k^2} \int_{B_{i}(R)} \Big|\nabla \big(|\phi_{j, i}^{(R)}|^{\frac{k}{2}}\big) \Big|^2 d\mu_i \nonumber \\
&\geq \frac{2 C_{\mathscr{SD}}}{k\lambda_{j, i}^{(R)}} \Big(\int_{B_i(R)} \Big|\phi_{j, i}^{(R)}\Big|^{\frac{kn}{n- 2}} d\mu_{i}\Big)^{\frac{n- 2}{n}} \nonumber
\end{align}
Denote $\beta= \frac{n}{n- 2}$, then for all $k\geq 2$, 
\begin{align}
\Vert \phi_{j, i}^{(R)}\Vert_{L^{k}(\mu_i)}\geq \Big(\frac{2 C_{\mathscr{SD}}}{k\lambda_{j, i}^{(R)}}\Big)^{\frac{1}{k}} \Vert \phi_{j, i}^{(R)}\Vert_{L^{k\beta}(\mu_i)} \nonumber 
\end{align}

Setting $k= 2\beta^s$ for $s= 0, 1, 2, \dots$, we have
\begin{align}
\Vert \phi_{j, i}^{(R)}\Vert_{L^{2\beta^{s+ 1}}(\mu_i)}\leq \Big(\frac{\beta^j \lambda_{j, i}^{(R)}}{C_{\mathscr{SD}}}\Big)^{\frac{1}{2\beta^s}} \cdot \Vert \phi_{j, i}^{(R)}\Vert_{L^{2\beta^s}(\mu_i)} \nonumber 
\end{align}
Iterating this estimate and using 
\begin{align}
\Vert \phi_{j, i}^{(R)}\Vert_{L^2(\mu_i)}= t_i^{-\frac{n}{4}} V(\sqrt{t_i})^{\frac{1}{2}} \Vert \phi_{j, i}^{(R)}\Vert_{L^2(\nu_i)}= t_i^{-\frac{n}{4}} V(\sqrt{t_i})^{\frac{1}{2}} \nonumber 
\end{align}
we conclude that
\begin{align}
\Vert \phi_{j, i}^{(R)}\Vert_{L^{2\beta^{s+ 1}}(\mu_i)}\leq \Big[\prod_{l= 0}^{s}\Big(\frac{\beta^l \lambda_{j, i}^{(R)}}{C_{\mathscr{SD}}}\Big)^{\frac{1}{2\beta^l}} \Big]\cdot  t_i^{-\frac{n}{4}} V(\sqrt{t_i})^{\frac{1}{2}} \nonumber
\end{align}

Let $s\rightarrow \infty$ and applying the fact that 
\begin{align}
\Vert \phi_{j, i}^{(R)}\Vert_{L^{\infty}(\nu_i)}= \Vert \phi_{j, i}^{(R)}\Vert_{L^{\infty}(\mu_i)}= \lim_{p\rightarrow \infty} \Vert \phi_{j, i}^{(R)}\Vert_{L^{p}(\mu_i)} \nonumber
\end{align}

We obtain
\begin{align}
\Vert \phi_{j, i}^{(R)}\Vert_{L^{\infty}(\nu_i)} &\leq \Big(C_{\mathscr{SD}}\Big)^{-\frac{n}{4}} \cdot \Big(\lambda_{j, i}^{(R)}\Big)^{\frac{n}{4}} \cdot C(n) t_i^{-\frac{n}{4}} V(\sqrt{t_i})^{\frac{1}{2}} \nonumber \\
&= C(n)\Big[\frac{V(\sqrt{t_i}) R^n}{V(\sqrt{t_i}R)}\Big]^{\frac{1}{2}}\cdot \Big(\lambda_{j, i}^{(R)}\Big)^{\frac{n}{4}} \nonumber \\
&\leq C(n) R^{\frac{n}{2}} \Big(\lambda_{j, i}^{(R)}\Big)^{\frac{n}{4}} \label{3.5.5}
\end{align}

Combining with Lemma \ref{lem 3.3}, we get 
\begin{align}
\Vert \phi_{j, i}^{(R)}\Vert_{L^{\infty}(\nu_i)} \leq C(n, R)j^{\frac{n}{2}} \nonumber
\end{align}
}
\qed

Note that the volume element $\nu_i$ of $(M_i, y, \rho_i, \nu_i)$ is not determined by the metric $\rho_i$, the heat kernel of $(M_i, y, \rho_i, \nu_i)$ is 
\begin{align}
H_i(x, y, s)= t_i^{\frac{n}{2}}\mu_i\Big(B_i(1)\Big)\cdot H(x, y, t_i s)= V(\sqrt{t_i})\cdot H(x, y, t_i s) \label{Hi}
\end{align}
where $H(x, y, s)$ is the heat kernel of $(M^n, y, g, \mu)$, $\mu$ is the volume element determined by $g$, $V(\sqrt{t_i})\doteqdot \int_{B(t_i)} 1 d\mu$, and $B(t_i)\doteqdot \{z\in M^n|\ d_{g}(z, y)\leq t_i\}$. Note $H_i(x, y, s)$ is different from the heat kernel $\tilde{H}_i(x, y, s)$ of $(M_i, y, \rho_i, \mu_i)$, which is $t_i^{\frac{n}{2}} H(x, y, t_i s)$.

Hence we have
\begin{align}
\lim_{i\rightarrow \infty} V(\sqrt{t_i})H(x, y, t_i)= \lim_{i\rightarrow \infty} H_{i}(x, y, 1) \label{3.1}
\end{align}
and by (\ref{1.1}),
\begin{align}
H_{i}(x, y, t)= V(\sqrt{t_i})H(x, y, t_it)\leq C(n) V(\sqrt{t_i})V(\sqrt{t_it})^{-1}e^{-\frac{d^2_{g}(x, y)}{5t_it}} \label{3.1.1}
\end{align}

Let us denote by $H_R(x, y, t)$ the Dirichlet heat kernel on the metric ball
\begin{align}
B(R)= \{z\in M^n|\ d_{g}(z, y)\leq R\}\subset (M^n, g, \mu) \nonumber
\end{align}
where $R> 0$ is a constant, and put $H_R= 0$ outside of $B(R)$. Similarly, we denote by $H_{R, i}(x, y, t)$ the Dirichlet heat kernel on $B_i(R)\subset (M_i, y, \rho_i, \nu_i)$.

From Lemma \ref{lem 3.3} and Lemma \ref{lem 3.5}, using similar argument in the proof of Theorem $10.1$ in \cite{Libook}, it is easy to get the following eigenfunction expansion of $H_{R, i}(x, y, t)$:
\begin{align}
H_{R, i}(x, y, t)= \sum_{j= 1}^{\infty}e^{-\lambda_{j, i}^{(R)}t} \phi_{j, i}^{(R)}(x) \phi_{j, i}^{(R)}(y) \label{H(R, i)}
\end{align}

\begin{lemma}\label{lem 3.4}
{For any $N> 0$, there exists a function $\epsilon (N, R, \delta)$ such that for any fixed $R> 2$, $\displaystyle \lim_{\delta\rightarrow 0}\epsilon (N, R, \delta)= 0$. And for $j$ satisfying $\lambda_{j, i}^{(R)}< N$, we have
\begin{align}
\int_{A_i(R- \delta, R)} \Big|\phi_{j, i}^{(R)}\Big|^2 d\nu_i \leq \epsilon (n, N, R, \delta) \quad for \ 0< \delta\leq 1 \label{3.6}
\end{align}
where $A_i(R- \delta, R)\doteqdot \{z\in M_i|\ R-\delta\leq d_{\rho_i}(z, y)\leq R\}$.
}
\end{lemma}

\pf
{Using (\ref{H(R, i)}) and (\ref{3.1.1}), we get
\begin{align}
\int_{A_i(R- \delta, R)} \Big|\phi_{j, i}^{(R)}(x)\Big|^2 d\nu_i(x) &\leq \int_{A_i} e^{\lambda_{j, i}^{(R)}} H_{R, i} (x, x, 1) d\nu_i(x) \nonumber \\
&\leq e^{\lambda_{j, i}^{(R)}}\int_{A_i} H_{i}(x, x, 1) d\nu_i(x) \nonumber \\
&\leq C(n) e^{N} \int_{A_i} e^{-\frac{1}{5t_i}} d\nu_i(x) \nonumber \\
&\leq C(n, N) \frac{\mu(A_i)}{V(\sqrt{t_i})} \leq C(n, N) \Big[R^n- (R- \delta)^n\Big] \nonumber 
\end{align}
in the last inequality above, we used the Bishop-Gromov inequality. Our conclusion is proved.
}
\qed

The following lemma follows from a standard argument of Cheng-Yau in \cite{CY}, which is needed in the proof of Proposition \ref{prop 3.7}.
\begin{lemma}\label{lem 3.6}
{Assume that $(M^n, g)$ is a complete manifold with $Rc\geq 0$, if $\Delta u= -\lambda u$ on $B_p(2r)\subset M^n$ and $\lambda\geq 0$, then we have
\begin{align}
|\nabla u|(x)\leq C(n)[r^{-1}+ \lambda] \cdot \sup_{x\in B_p(2r)}|u(x)|\ , \quad x\in B_p(r) \nonumber
\end{align} 
where $B_p(r)= \{z\in M^n |\ d_g(z, p)\leq r\}$.
}
\end{lemma}

\pf
{Let $\displaystyle \mathscr{M}= \sup_{x\in B_p(2r)}|u(x)|$, $f(x)= u(x)+ \mathscr{M}$, without loss of generality, assume $\mathscr{M}> 0$. It is easy to get
$\Delta f= -\lambda f+ \lambda \mathscr{M}$ on $B_{p}(2r)$, and $f\geq 0$. 

Apply Theorem $6$ in \cite{CY} to $f(x)$, we get 
\begin{align}
|\nabla f(x)|\leq C(n)[r^{-1}+ \lambda] \cdot [f(x)+ \mathscr{M}]\ , \quad x\in B_{p}(r) \label{3.6.1}
\end{align}
By the definition of $f(x)$ and $\mathscr{M}$, our conclusion follows from (\ref{3.6.1}).
}
\qed

\begin{prop}\label{prop 3.7}
{For fixed $j$, $k> 0$, assume (for a subsequence of the eigenvalues) $\lambda_{j, i}^{(R)}\rightarrow \lambda_{j, \infty}^{(R)}$, $\lambda_{k, i}^{(R)}\rightarrow \lambda_{k, \infty}^{(R)}$ as $i\rightarrow \infty$. Then there is a subsequence (denoted also by $\phi_{j, i}^{(R)}$, $\phi_{k, i}^{(R)}$) that converges uniformly on compact subsets of $\mathring{B}_{\infty}(R)$, and also in $L^2\Big(B_{\infty}(R)\Big)$, to two compactly supported Lipschitz functions $\phi_{j, \infty}^{(R)}$, $\phi_{k, \infty}^{(R)}$ on $B_{\infty}(R)$, where $B_{\infty}(R)= \{z\in M_{\infty}|\ d_{\rho_{\infty}}(z, y)\leq R\}$, $\mathring{B}_{\infty}(R)$ denotes the interior of $B_{\infty}(R)$. Moreover, 
\begin{align}
&\Delta_{\infty} \phi_{j, \infty}^{(R)}= \lambda_{j, \infty}^{(R)}\phi_{j, \infty}^{(R)}\ , \quad \Delta_{\infty} \phi_{k, \infty}^{(R)}= \lambda_{k, \infty}^{(R)}\phi_{k, \infty}^{(R)}\ , \label{3.7.1} \\ 
&\int_{B_{\infty}(R)} \phi_{j, \infty}^{(R)}\phi_{k, \infty}^{(R)} d\nu= \delta_{j, k} \label{3.7.2}
\end{align} 
}
\end{prop}

\pf
{Locally uniform convergence follows from Lemma \ref{lem 3.5} and \ref{lem 3.6}. The $L^2$ convergence and (\ref{3.7.2}) are implied by locally uniform convergence and Lemma \ref{lem 3.4}. Finally, (\ref{3.7.1}) follows from Theorem \ref{thm 4.2} and Lemma \ref{lem 3.6}.
}
\qed

\section{Solutions of the heat equations on metric measure spaces}\label{SECTION 5}

In this section, on metric measure spaces, we will show the existence of the solution of the heat equations and the parabolic mean value inequality. For smooth manifolds, all these results are well-known. On metric measure spaces, our setup is closely related with the discussion in \cite{Sturm2}. 

Assume $U\subset M_{\infty}$, and $U$ is open. We will be concerned with the following Banach spaces.
\begin{itemize}
\item $L^2\big((0, T); H_0^1(U)\big)$ is the Hilbert space consisting of functions $u(x, t)$, measurable on $(0, T)$ with range in $H_0^1(U)$ (for the Lebesgue measure $dt$ on $(0, T)$), for any $t\in (0, T)$, $u(\cdot, t)\in H_0^1(U)$ and the norm of the space is $\displaystyle \Big(\int_0^T \big|u(\cdot, t)\big|_{H_0^1(U)}^2 dt\Big)^{\frac{1}{2}}$
\item $H^1\big((0, T); H_0^1(U)^{*}\big)$ is the Sobolev space of functions $u$, where $H_0^1(U)^{*}$ is the dual space of $H_0^1(U)$, and $u\in L^2\big((0, T); H_0^1(U)^{*}\big)$, and it has distributional time derivative $\frac{\partial}{\partial t}u\in L^2\big((0, T); H_0^1(U)^{*}\big)$ equipped with the norm \[\displaystyle \Big(\int_0^T \big|u(\cdot, t)\big|_{H_0^1(U)^{*}}^2+ \big|\frac{\partial}{\partial t}u(\cdot, t)\big|_{H_0^1(U)^{*}}^2 dt\Big)^{\frac{1}{2}}\].
\item $\displaystyle \mathscr{F}\big((0, T)\times U\big)\doteqdot L^2\big((0, T); H_0^1(U)\big)\cap H^1\big((0, T); H_0^1(U)^{*}\big)$. We mention the following important result from \cite{RR}:
\begin{align}
\mathscr{F}\big((0, T)\times U\big)\subset C\big([0, T], L^2(U)\big) \nonumber
\end{align}
\item Similarly, $\displaystyle \mathscr{G}\big((0, T)\times U\big)\doteqdot L^2\big((0, T); H^1(U)\big)\cap H^1\big((0, T); H^1(U)^{*}\big)$.
\end{itemize}

\begin{definition}\label{def solution}
{A function $u$ is called a \textbf{Dirichlet solution} of the heat equation on $(0, T)\times U$:
\begin{align}
\frac{\partial}{\partial t}u= \Delta_{\infty} u \quad on\ (0, T)\times U \label{solution 1}
\end{align}
iff $u\in \mathscr{F}\big((0, T)\times U\big)$, and for all $\phi\in \mathscr{F}\big((0, T)\times U\big)$:
\begin{align}
\int_0^T \int_U <du, d\phi> d\nu_{\infty} dt+ \int_0^T\int_{U} \frac{\partial u}{\partial t}\cdot \phi d\nu_{\infty}dt= 0 \label{solution 2}
\end{align}
}
\end{definition}

\begin{remark}\label{rem subso}
{For $u\in \mathscr{G}\big((0, T)\times U\big)$, we say that 
\begin{align}
\Big(\frac{\partial}{\partial t}- \Delta_{\infty}\Big) u= (\leq) 0 \quad \quad on\ (0, T)\times U \nonumber
\end{align}
if for almost all $t\in (0, T)$ except a subset of $(0, T)$ with Lebesgue measure $0$, 
\begin{align}
\int_U <du, d\phi> d\nu_{\infty}+ \int_{U} \frac{\partial u}{\partial t}\cdot \phi d\nu_{\infty}= (\leq) 0 \nonumber
\end{align} 
holds for all non-negative $\phi\in H_0^1(U)$. Such $u$ is also called a \textbf{solution (subsolution) of the heat equation} on $(0, T)\times U$.
}
\end{remark}

\begin{definition}\label{def ini sol}
{Given a function $f\in L^2(U)$, the function $u$ is called a \textbf{Dirichlet solution of the initial value problem} on $[0, T)\times U$:
\begin{equation}\label{ini sol 1}
\left\{
\begin{array}{rl}
\frac{\partial}{\partial t}u= \Delta_{\infty} u &\quad on \ (0, T)\times U \\
u(\cdot , 0)= f(\cdot) &\quad on \ U \\
\end{array} \right.
\end{equation}
iff $u$ is a Dirichlet solution of (\ref{solution 1}) and $\lim_{t\rightarrow 0}\int_U |u(x, t)- f(x)|^2 d\nu_{\infty}= 0$.
}
\end{definition}




\begin{prop}\label{prop exist}
{For every $f\in L^2(U)$, there exists a unique Dirichlet solution $u\in \displaystyle \mathscr{F}\big((0, T)\times U\big)$ of the initial value problem (\ref{ini sol 1}).
}
\end{prop}

\pf
{It follows from Theorem $4.1$ and Remark $4.3$ in Chapter $3$ of \cite{LM}.
}
\qed

For the solutions of heat equations on $M_{\infty}$, we have the following mean value inequalities.
\begin{theorem}\label{thm 6.8}
{If $\frac{\partial u}{\partial t}- \Delta_{\infty} u= 0$ in $Q_1$, then for any $0< \delta< 1$, we have
\begin{align}
\sup_{z\in Q_{\delta}} u^2(z)& \leq \frac{C(n)}{(1- \delta)^{n+ 2} r^2 \nu_{\infty}(B)} \int_{Q_1} u^2 d\nu_{\infty}dt \label{mean value inq} \\
\sup_{z\in Q_{\delta}} u(z)& \leq \frac{C(n)}{(1- \delta)^{n+ 2} r^2 \nu_{\infty}(B)} \int_{Q_1} u d\nu_{\infty}dt \label{L1 MVI}
\end{align}
where $B= B_{\infty}(x, r)$, $s> r^2 > 0$, $\tau> 0$ is a fixed positive constant, and  
\begin{align}
Q_1 \doteqdot (s- r^2, s)\times B_{\infty}(x, r)\ , \quad Q_{\delta}\doteqdot (s- \delta r^2, s)\times B_{\infty}(x, \delta r) \nonumber
\end{align}
}
\end{theorem}

\begin{remark}\label{rem MVI}
{The parabolic mean value inequality on smooth manifold were firstly proved in \cite{LT}, however the proof there used the upper bound of heat kernel, which is the target we want to prove. The conclusion on metric measure spaces was essentially obtained in \cite{SB}, although the context there are smooth manifolds. The following argument is just slight modification of the original argument there, hence it is sketchy. For the complete details, we refer the reader to that book.
}
\end{remark}

\pf
{Firstly, from the argument of Lemma $5.3.2$, Lemma $5.2.5$ in \cite{SB} and Proposition \ref{prop 3.9}, we can get the following Dirichlet Poincar\'{e} Inequality:

There exists positive constant $C(n)> 0$, such that for any $B= B_{\infty}(x, r)\subset M_{\infty}$,
\begin{align}
|f|_{L^2}\leq C(n) r|g_f|_{L^2} \ , \quad f\in H_0^1(B) \label{Local Poincare}
\end{align}

Secondly, from the argument of Theorem $5.3.3$ in \cite{SB}, Proposition \ref{prop 3.9} and (\ref{Local Poincare}), we can obtain Local Sobolev Inequality as the following:

There exists $C(n)> 0$, such that for any $B= B_{\infty}(x ,r)\subset M_{\infty}$, we have 
\begin{align}
\Big(\int_B |f|^{\frac{2n}{n- 2}} d\nu_{\infty}\Big)^{\frac{n- 2}{n}}\leq C(n)\frac{r^2}{\nu_{\infty}(B)^{\frac{2}{n}}} \Big(\int_B |g_f|^2 d\nu_{\infty}\Big)\ , \quad f\in H_0^1(B) \label{sobo inq}
\end{align}

Next, employing (\ref{sobo inq}), we can use almost exactly the same argument of Theorem $5.2.9$ in \cite{SB} to get the following two inequalities: 

If $\frac{\partial u}{\partial t}- \Delta_{\infty} u\leq 0$ in $Q_1$ and $u\geq 0$, then for any $0< \delta< 1$, (\ref{mean value inq}) and (\ref{L1 MVI}) hold.

Finally, for any $\epsilon> 0$, it is easy to show that $v\doteqdot \sqrt{u^2+ \epsilon}$ is the solution of the heat equation, which was defined in Remark \ref{rem subso}, and $v\geq 0$. By the above argument,
\begin{align}
\sup_{z\in Q_{\delta}} (u^2+ \epsilon)(z)\leq \frac{C(n)}{(1- \delta)^{n+ 2} r^2 \nu_{\infty}(B)} \int_{Q_1} (u^2+ \epsilon) d\nu_{\infty} dt\label{6.8.1}
\end{align}
Let $\epsilon\rightarrow 0$ in (\ref{6.8.1}), we get (\ref{mean value inq}).

Similar argument yields (\ref{L1 MVI}).
}
\qed





\section{The existence and Gaussian upper bound of heat kernel on limit spaces}\label{SECTION 6}

In this section we will prove the existence of heat kernel on limit spaces under Gromov-Hausdorff convergence, and establish Gaussian upper bound of heat kernel. 

To prove the existence of heat kernel on limit spaces, we are inspired by the method of K.-T. Sturm in \cite{Sturm2}. Firstly, from Proposition \ref{prop exist}, there exists a uniquely determined operator:
\begin{align}
T: L^2\big(M_{\infty}\big)\rightarrow \mathscr{F}\big((0, T)\times M_{\infty}\big)\label{T}
\end{align}
with the property that for every $f\in L^2(M_{\infty})$, the unique Dirichlet solution of (\ref{ini sol 1}) ($U= M_{\infty}$ there) is given by $u(x, t)= [Tf](x, t)$. 

We also define $\big[T_t f\big](x)= [Tf](x, t)$ for every $t\in (0, T)$, then 
\begin{align}
T_{t}:\ L^2\big(M_{\infty}\big)\rightarrow L^2\big(M_{\infty}\big) \label{T(t)}
\end{align}

\begin{lemma}\label{lem 6.9}
{There exists $C(n)> 0$ such that for any $t\in (0, 8R^2)$,  
\begin{align}
\sup_{x\in B_{\infty}(R)} \big|\big(T_{t}f\big) (x)\big| \leq C(n) \Big(\frac{R}{\sqrt{t}}\Big)^{n+ 2} \nu_{\infty}\big(B_{\infty}(R)\big)^{-\frac{1}{2}} |f|_{L^2(M_{\infty})}\ , \quad \forall f\in L^2\big(M_{\infty}\big) \nonumber
\end{align}
where $R> 0$ is any positive constant.
}
\end{lemma}

\pf
{We will apply Theorem \ref{thm 6.8} on $T_t f$ for given $t\in (0, 4R^2)$. Let $r= 2R$, $\delta= 1- \frac{t}{10 R^2}$, $s= (2R)^2+ \frac{1}{2}t$, $\tau= 1$ in (\ref{mean value inq}), note that $t\in (s- \delta r^2, s)$, then we get
\begin{align}
\sup_{x\in B_{\infty}(R)} |(T_t f)(x)|&\leq \sup_{Q_{\delta}} |(T_t f)(x)| \nonumber \\
&\leq C(n) \Big(\frac{1}{1- \delta}\Big)^{\frac{n+ 2}{2}} \Big(\frac{1}{(2R)^2 \nu_{\infty}\big(B_{\infty}(2R)\big)} \int_{Q_1} |T_t f|^2\Big)^{\frac{1}{2}} \nonumber \\ 
&\leq C(n) \Big(\frac{R}{\sqrt{t}}\Big)^{n+ 2} \nu_{\infty}\big(B_{\infty}(R)\big)^{-\frac{1}{2}} |f|_{L^2(M_{\infty})} \nonumber
\end{align}
in the last inequality, we used that
\begin{align}
\int_{M_{\infty}} |T_t f|^2 d\nu_{\infty} \leq \int_{M_{\infty}} |f|^2 d\nu_{\infty}\ , \quad \forall t> 0 \nonumber
\end{align}
which follows from (\ref{solution 2}).
}
\qed

We also have the following parabolic maximum principle on $M_{\infty}$ (for the proof, see Proposition $4.11$ in \cite{GH2}).
\begin{lemma}[\cite{GH2}]\label{lem max prin}
{Assume $h$ is a solution of the heat equation on $(0, T+1)\times B_{\infty}(z, R)$, and
\begin{align}
\lim_{t\rightarrow 0} \int_{B_{\infty}(z, R)} h^2(x, t) d\nu_{\infty}(x)= 0\ , \quad h|_{\partial B_{\infty}(z, R)\times (0, T]}\leq 0 \label{3.12.1}
\end{align}
for any $f(x)\in L^2\big(B_{\infty}(z, R)\big)$. Then $h\leq 0$ on $(0, T]\times B_{\infty}(z, R)$.
}
\end{lemma}

The following result is one modification of classical result in functional analysis, which was due to J-X Hu and Grigor'yan (see Lemma $3.3$ in \cite{GH}).
\begin{lemma}[\cite{GH}]\label{lem 7.2}
{Let $K: L^2(Y)\rightarrow L^{\infty}(X)$ be a bounded linear operator, with the norm bounded by $C$, that is, for any $f\in L^2(Y)$, 
\begin{align}
\sup_{X} |Kf|\leq C|f|_2 \label{7.2.1}
\end{align}
There exists a mapping $k: X\rightarrow L^2(Y)$ such that, for all $f\in L^2(Y)$, and almost all $x\in X$,
\begin{align}
Kf(x)= (k(x), f) \label{7.2.2}
\end{align} 
Moreover, for all $x\in X$, $||k(x)||_{L^2(Y)}\leq C$. Furthermore, there is a function $k(x, y)$ that is jointly measurable in $(x, y)\in M\times M$, such that, for almost all $x\in X$, $k(x, \cdot)= k(x)$ almost everywhere on $Y$.
}
\end{lemma}

Now we can prove the existence of the heat kernel with respect to the Dirichlet boundary condition on $M_{\infty}$.

\begin{theorem}\label{thm 6.10}
{There exists a nonnegative measurable function 
\[p_{\infty}: M_{\infty}\times M_{\infty}\times \mathbb{R}^{+} \rightarrow [0, \infty]\] with the following properties:
\begin{enumerate}
\item On $[0, \infty)\times M_{\infty}$, the function 
\begin{align}
u(x, t)= \int_{M_{\infty}} p_{\infty}(x, z, t)f(z) d\nu_{\infty}(z) \nonumber
\end{align}
is a solution of (\ref{ini sol 1}), where $f\in L^2\big(M_{\infty}\big)$.
\item 
For any fixed $w\in M_{\infty}$, any $T> 0$, 
\[p_{\infty}(x, w, t)\in L^2\big((0, T); H_0^1(M_{\infty})\big)\cap H^1\big((0, T); H_0^1(M_{\infty})^{*}\big)\]
is a Dirichlet solution of the heat equation (defined as in Definition \ref{def solution}).
\end{enumerate}
}
\end{theorem}

\begin{remark}\label{rem heat}
{Such $p_{\infty}$ is called the \textbf{heat kernel of $(M_{\infty}, \rho_{\infty}, \nu_{\infty})$}.
}
\end{remark}

\pf
{By Lemma \ref{lem 6.9} and Lemma \ref{lem 7.2}, there exists $p_{\infty}(x, z, t)$, which is jointly measurable in $(x, z)\in M_{\infty}\times M_{\infty}$, such that 
\[T_{t}(f)(x)= \int_{M_{\infty}} p_{\infty}(x, z, t)f(z) d\nu_{\infty}(z) \]
From Lemma \ref{lem max prin}, we get that if $f\geq 0$, $T_{t}(f)\geq 0$. It follows from Lemma $3.2$ in \cite{GH}, $p_{\infty}(x, z, t)\geq 0$. Then $p_{\infty}\geq 0$ and the conclusion in $(1)$ above are proved.

For any $f\in L^2\big(B_{\infty}(R)\big)$, from the uniqueness of solution in Proposition \ref{prop exist} and the definition of $T$, $T_{t}$, we get
\begin{align}
\big[T_{t+s}f\big](z)&= \big[T f\big](z, t+s)= T\big[(T f)(\cdot, s)\big](z, t) \nonumber \\
&= T_{t}\big[(T f)(\cdot, s)\big](z)= \int_{M_{\infty}} p_{\infty}(z, x, t)\cdot \big[T f\big](x, s) d\nu_{\infty}(x) \nonumber \\
&= \int_{M_{\infty}} p_{\infty}(z, x, t)\Big(\int_{M_{\infty}} p_{\infty}(x, w, s)f(w) d\nu_{\infty}(w)\Big) d\nu_{\infty}(x) \nonumber \\
&= \int_{M_{\infty}} \Big(\int_{M_{\infty}} p_{\infty}(z, x, t)p_{\infty}(x, w, s) d\nu_{\infty}(x)\Big)\cdot f(w) d\nu_{\infty}(w) \nonumber
\end{align}
Hence we have
\begin{align}
p_{\infty}(z, w, t+s)&= \int_{M_{\infty}} p_{\infty}(z, x, t)p_{\infty}(x, w, s) d\nu_{\infty}(w) \nonumber\\
&= \Big[T_{t}\big(p_{\infty}(\cdot, w, s)\big)\Big](z)= \Big[T\big(p_{\infty}(\cdot, w, s)\big)\Big](z, t) \nonumber
\end{align} 
By the definition of $T$, and $t+s$ can be chosen as any positive number, we get that $p_{\infty}(z, w, t)$ is a Dirichlet solution of the heat equation on $(0, \infty)\times M_{\infty}$. 
}
\qed


And we have the following theorem about the upper bound of $p_{\infty}(x, y, t)$. We will follow the method developed by E. B. Davies on smooth manifolds (see \cite{Davies}, also \cite{SB}), our proof is just slight modification of the proof given in \cite{SB}, and it is presented here for completeness and reader's convenience.
\begin{theorem}\label{thm 6.11}
{Assume that $p_{\infty}(x, y, t)$ is the heat kernel of $\big(M_{\infty}, y, \rho_{\infty}, \nu_{\infty}\big)$, then
\begin{align}
p_{\infty}(x, y, t)\leq C(n) \nu_{\infty}\Big(B_{\infty}(y, \sqrt{t})\Big)^{-1} e^{-\frac{1}{5t} \rho_{\infty}^2(x, y)} \label{6.11.1}
\end{align}
where $C(n)$ is the positive constant depending only on $n$.
}
\end{theorem}

We firstly need to prove a lemma. 
\begin{lemma}\label{lem 6.12}
{For any function $\phi\in H_0^1({M_\infty})$ with $|g_{\phi}|\leq 1$ and any $\alpha\in \mathbb{R}$, we define the operator $H_{t}^{\alpha, \phi}$ as the following:
\begin{align}
H_t^{\alpha, \phi} f(x)\doteqdot e^{-\alpha\phi(x)}\int_{M_{\infty}} p_{\infty}(x, y, t)e^{\alpha\phi(y)} f(y) d\nu_{\infty}(y) \label{Ht}\ , \quad \quad f\in L^2(M_{\infty})
\end{align}
Then as an operator from $L^2(M_{\infty})$ to $L^{2}(M_{\infty})$, $H_t^{\alpha, \phi}$ satisfies $||H_t^{\alpha, \phi} ||\leq e^{\alpha^2 t}$. 
}
\end{lemma}

\pf
{For any $f\in L^2(M_{\infty})$, set $u(t)= |H_t^{\alpha, \phi} f|_{L^2}^2$, then
\begin{align}
u'(t)&= 2\int_{M_\infty} \frac{\partial}{\partial t} \big(H_t^{\alpha, \phi} f\big)\cdot H_t^{\alpha, \phi} f \nonumber \\
&= \int_{M_{\infty}} e^{-\alpha \phi(x)} \Delta \Big(e^{\alpha \phi(x)} H_t^{\alpha , \phi} f(x) \Big) H_t^{\alpha , \phi} f(x) d\nu_{\infty}(x) \nonumber \\
&= -2\int_{M_{\infty}} <d \Big(e^{\alpha \phi(x)}H_t^{\alpha, \phi} f(x)\Big), d\Big(e^{-\alpha \phi(x)}H_t^{\alpha, \phi}f(x)\Big)> \nonumber \\
&= 2\Big[\alpha^2 \int_{M_{\infty}} |d\phi|^2 |H_t^{\alpha, \phi} f|^2- \int_{M_{\infty}} |H_t^{\alpha, \phi} f|^2 \Big]\leq 2\alpha^2 u(t) \nonumber
\end{align}
Hence $u(t)\leq e^{2\alpha^2 t} u(0)$, note $u(0)= \Big|f \Big|_{L^2}^2$, we get
\begin{align}
\Big|H_t^{\alpha, \phi}f\Big|_{L^2}^2\leq e^{2\alpha^2 t} \Big|f\Big|_{L^2}^2 \nonumber
\end{align}
The conclusion follows from the above inequality.
}
\qed

\bigskip
{\it \textbf{Proof of Theorem \ref{thm 6.11}}:}~
{Fix $x$, $y\in M_{\infty}$, and $r_1$, $r_2> 0$. Let $\chi_1$ (respectively $\chi_2$) be the function equal to $1$ on $B_1= B_{\infty}(x, r_1)$ (respectively $B_2= B_{\infty}(y, r_2)$) and equal to $0$ otherwise. Then
\begin{align}
&\int_{B_1}\int_{B_2} p_{\infty}(\xi, \zeta, t)e^{-\alpha(\phi(\xi)- \phi(\zeta))} d\zeta d\xi = \int_{M_{\infty}} \chi_1(\xi) \big(H_t^{\alpha, \phi} \chi_2\big)(\xi) d\xi \nonumber \\
&\quad \quad \quad \leq || H_t^{\alpha, \phi}||\cdot ||\chi_1 ||_{L^2}\cdot || \chi_2 ||_{L^2}\leq e^{\alpha^2 t} \nu_{\infty}(B_1)^{\frac{1}{2}} \nu_{\infty} (B_2)^{\frac{1}{2}} \nonumber 
\end{align}

Using $|d\phi|\leq 1$, we get 
\begin{align}
&\int_{B_1}\int_{B_2} p_{\infty}(\xi, \zeta, t) d\zeta d\xi\leq \int_{B_1}\int_{B_2} p_{\infty} e^{-\alpha\big[(\phi(\xi)- \phi(x))- (\phi(\zeta)- \phi(y))\big]} \cdot e^{|\alpha|(r_1+ r_2)} \nonumber \\
&\quad \quad \quad \leq\Big[\nu_{\infty}(B_1)\nu_{\infty}(B_2)\Big]^{\frac{1}{2}} \exp\{\alpha^2 t+ \alpha [\phi(x)- \phi(y)]+ |\alpha|(r_1+ r_2)\} \nonumber 
\end{align}

As $p_{\infty}(x, \cdot, t)$ is a Dirichlet solution of heat equation in $(0, \infty)\times M_{\infty}$, assume $t\geq \frac{1}{4}r_2^2$ and applying Theorem \ref{thm 6.8}, we obtain
\begin{align}
p_{\infty}(\xi, y, t)\leq \frac{C(n)}{r_2^2 \nu_{\infty}(B_2)} \int_{t- \frac{1}{4}r_2^2}^t \int_{B_2} p_{\infty}(\xi, \zeta, s)d\zeta ds \nonumber
\end{align}

Thus
\begin{align}
\int_{B_1} p_{\infty}(\xi, y, t) d\xi\leq \frac{C(n)\nu_{\infty}(B_1)^{\frac{1}{2}}}{\nu_{\infty}(B_2)^{\frac{1}{2}}}\cdot \exp\Big\{\alpha^2 t+ \alpha[\phi(x)- \phi(y)]+ |\alpha| (r_1+ r_2)\Big\} \nonumber
\end{align}

Assume $t\geq \frac{1}{4}(r_1^2+ r_2^2)$, by Theorem \ref{thm 6.8} again, combining with the above inequality, we get
\begin{align}
p_{\infty}(x, y, t)&\leq \frac{C(n)}{r_1^2\nu_{\infty}(B_1)}\int_{t- \frac{1}{4}r_1^2}^t \int_{B_1} p_{\infty}(\xi, y, s) d\xi ds \nonumber \\
&\leq \frac{C(n)}{\Big[\nu_{\infty}(B_1)\nu_{\infty}(B_2)\Big]^{\frac{1}{2}}} \exp\{\alpha^2 t+ \alpha[\phi(x)- \phi(y)]+ |\alpha|(r_1+ r_2)\} \nonumber 
\end{align}

Taking $\alpha= \frac{\phi(y)- \phi(x)}{2t}$, $r_1= r_2= \frac{t}{\sqrt{t}+ \rho_{\infty}(x, y)}$, we obtain
\begin{align}
p_{\infty}\leq \frac{C(n)}{\Big[\nu_{\infty}(B_1)\nu_{\infty}(B_2)\Big]^{\frac{1}{2}}}\cdot \exp\{-\frac{\big(\phi(x)- \phi(y)\big)^2}{4t}+ \frac{|\phi(x)- \phi(y)|}{\sqrt{t}+ \rho_{\infty}(x, y)}\} \label{6.11.3}
\end{align}

Choosing $\phi(\cdot)= \rho_{\infty}(x, \cdot)$ in (\ref{6.11.3}) gives
\begin{align}
p_{\infty}(x, y, t)&\leq \frac{C(n)}{\Big[\nu_{\infty}(B_1)\nu_{\infty}(B_2)\Big]^{\frac{1}{2}}} \exp\{-\frac{\rho_{\infty}^2(x, y)}{4t}\} \nonumber \\
&\leq \frac{C(n)\Big(1+ \frac{\rho_{\infty}(x, y)}{\sqrt{t}}\Big)^{\frac{n}{2}}}{\Big[\nu_{\infty}\big(B_{\infty}(x, \sqrt{t})\big)\nu_{\infty}\big(B_{\infty}(y, \sqrt{t})\big)\Big]^{\frac{1}{2}}} \exp\{-\frac{\rho_{\infty}^2(x, y)}{4t}\} \nonumber \\
&\leq \frac{C(n)}{\nu_{\infty}\big(B_{\infty}(y, \sqrt{t})\big)} \exp\{-\frac{\rho_{\infty}^2(x, y)}{5t}\} \nonumber 
\end{align}
The conclusion is proved.
}
\qed

\begin{cor}\label{cor 3.10}
{For positive constant $T> 0$, there exists a positive constant $\epsilon (n, T, R)$ with $\displaystyle \lim_{R\rightarrow \infty} \epsilon (n, T, R)= 0$ such that for $t\in (0, T]$: 
\begin{align}
\int_{M_{\infty}\backslash B_{\infty}(R)} p_{\infty}(x, y, t) d\nu_{\infty}(x) \leq \epsilon (n, T, R) \label{c 3.10.1}
\end{align}
}
\end{cor}

\pf
{From (\ref{6.11.1}) and Property ($\mathscr{B}$) on $(M_{\infty}, \rho_{\infty}, \nu_{\infty})$, we get 
\begin{align}
\int_{M_{\infty}\backslash B_{\infty}(R)} p_{\infty}(x, y, t) d\nu_{\infty}(x)&\leq C(n)\int_{M_{\infty}\backslash B_{\infty}(R)} \nu_{\infty}\Big(B_{\infty}(\sqrt{t})\Big)^{-1} e^{-\frac{\rho_{\infty}^2(x, y)}{5t}} d\nu_{\infty} \nonumber \\
&\leq \frac{C(n)}{\nu_{\infty}\big(B_{\infty}(\sqrt{t})\big)} \sum_{k= 0}^{\infty} \int_{B_{\infty}(2^{k+1}R)\backslash B_{\infty}(2^k R)} e^{-\frac{\rho_{\infty}^2(x, y)}{5t}} d\nu_{\infty} \nonumber \\
&\leq \frac{C(n)}{\nu_{\infty}\big(B_{\infty}(\sqrt{t})\big)} \sum_{k= 0}^{\infty} e^{-\frac{(2^k R)^2}{5t}}\cdot (2^k R)^n \nonumber \\
&\leq C(n, T)\sum_{k= 0}^{\infty} e^{-\frac{1}{5}\Big(\frac{2^k R}{\sqrt{t}}\Big)^2}\cdot \Big(\frac{2^k R}{\sqrt{t}}\Big)^n \nonumber 
\end{align}

Without loss of generality, we can assume $R\geq \sqrt{T}$. Then from the above,
\begin{align}
\int_{M_{\infty}\backslash B_{\infty}(R)} p_{\infty}(x, y, t) d\nu_{\infty}(x)\leq C(n, T)\int_{\frac{R}{\sqrt{T}}}^{\infty} e^{-\frac{1}{5}s^2}s^n ds \leq \epsilon (n, T, R) \nonumber
\end{align}

}
\qed

\section{The convergence of heat kernels in the Gromov-Hausdorff sense}\label{SECTION 07}

In this section, we will prove one main theorem of this paper, Theorem \ref{thm 8.13}. The eigenfunction expansion of heat kernel and Proposition \ref{prop 3.7} provides the bridge between local Dirichlet heat kernels on bounded regions of $M_i$ and $M_{\infty}$. Combined with Gaussian upper bounds of heat kernels on $M_i$, $M_{\infty}$, maximum principle leads to the convergence of local Dirichlet heat kernel to global Dirichlet heat kernel on $M_i$, $M_{\infty}$. From all these, the hear kernels' convergence in the Gromov-Hausdorff sense is proved.


\begin{lemma}\label{lem 3.1}
{For positive constant $T> 0$, there exists $\epsilon (n, T, R)> 0$ with $\lim_{R\rightarrow \infty} \epsilon (n, T, R)= 0$, such that for $t\in (0, T]$: 
\begin{align}
\int_{M_i\backslash B_i(R)} H_i(x, y, t) d\nu_i(x) \leq \epsilon (n, T, R) \label{3.2}
\end{align} 
}
\end{lemma}

\pf
{Without loss of generality, assume $R\geq \sqrt{T}$, then from $\nu_i= \frac{1}{V(\sqrt{t_i})} \mu$ and (\ref{3.1.1}), we get
\begin{align}
\int_{M_i\backslash B_i(R)} H_{i}(x, y, t) d\nu_i(x) &\leq C(n)\int_{M^n\backslash B(\sqrt{t_i}R)} V(\sqrt{t_it})^{-1} e^{-\frac{d^2}{5t_it}} d\mu \nonumber \\
&= \frac{C(n)}{V(\sqrt{t_i t})} \int_{\sqrt{t_i}R}^{\infty} e^{-\frac{r^2}{5t_it}}A(r)r^{n- 1} dr\nonumber \\
&\leq C(n) \frac{A(\sqrt{t_i T}) \cdot \big(t_i t\big)^{\frac{n}{2}}}{V(\sqrt{t_it})}\cdot \Big(\int_{\frac{R}{\sqrt{t}}}^{\infty} e^{-\frac{1}{5}s^2} s^{n- 1} ds\Big) \nonumber \\
&\leq C(n) \frac{A(\sqrt{t_i T})}{A(\sqrt{t_i t})} \cdot \Big(\int_{\frac{R}{\sqrt{T}}}^{\infty} e^{-\frac{1}{5}s^2} s^{n- 1} ds\Big) \nonumber\\
&\leq C(n)\cdot \int_{\frac{R}{\sqrt{T}}}^{\infty} e^{-\frac{1}{5}s^2} s^{n- 1} ds \leq \epsilon \Big(n, T, R\Big) \nonumber
\end{align}
where $A(r)r^{n- 1}$ in the first equality is the surface area element of $\partial B(r)$, in the second inequality above we used the fact that $A(r)$ is non-increasing (from Bishop-Gromov inequality) and $R\geq \sqrt{T}$; the third inequality from the end follows from the fact $V(\sqrt{t_i t})\geq \frac{1}{n}A(\sqrt{t_i t}) \big(t_i t\big)^{\frac{n}{2}}$.
}
\qed

\begin{prop}\label{prop 3.2}
{\begin{align}
\lim_{R\rightarrow \infty} H_{R, i}(\cdot , y, t)= H_{i}(\cdot , y, t) \label{3.3}
\end{align}
The convergence is uniform on $x\in M_i^n$, $i= 1, 2, \dots$, and uniform in $L^1(\nu_i)$ on any finite time interval $(0, T]$.
}
\end{prop}

\pf
{Assume $R\geq 1$, put 
\begin{align}
M_{R, i}\doteqdot \sup\{H_i(x, y, t)|\ x\in \partial B_i(R),\ 0< t\leq T \} \label{3.2.1}
\end{align}
By (\ref{3.1.1}) and Volume Comparison Theorem, we have
\begin{align}
M_{R, i}&\leq \sup_{0< t\leq T} C(n) V(\sqrt{t_i})V(\sqrt{t_i t})^{-1} e^{-\frac{R^2}{5t}} \nonumber \\
&\leq C(n)\cdot \max\{e^{-\frac{R^2}{5T}},\ \sup_{0< t\leq 1} t^{-\frac{n}{2}} e^{-\frac{R^2}{5t}} \} \label{3.2.2}\\
&\leq C(n)\max\{e^{-\frac{R^2}{5T}},\ R^{-n}\} \label{3.2.3}
\end{align}

By the maximum principle, when $x\in B_i(R)$,
\begin{align}
H_{i}(x, y, t)- M_{R, i}\leq H_{R, i}(x, y, t)\leq H_i(x, y, t) \label{3.2.4}
\end{align}

From (\ref{3.2.3}) and (\ref{3.2.4}), we get $\displaystyle \lim_{R\rightarrow \infty} H_{R, i}(\cdot, y, t)= H_i(\cdot, y, t)$ uniformly on $(0, T]\times B_i(R)$, $i= 1, 2, 3, \cdots$. Combining with (\ref{3.1.1}), we get that the convergence is uniform on $(0, T]\times M_i^n$ and $i= 1, 2, 3, \cdots$.

From (\ref{3.2.2}) and Volume Comparison Theorem, we get
\begin{align}
\lim_{R\rightarrow \infty} M_{R, i}\nu_i\Big(B_i(R)\Big)&\leq \lim_{R\rightarrow \infty}C(n)R^{n} \cdot \max\{e^{\frac{R^2}{5T}},\ \sup_{0< t\leq 1} t^{-\frac{n}{2}} e^{-\frac{R^2}{5t}}\} \nonumber \\
&\leq \lim_{R\rightarrow \infty} C(n) \max\{R^n e^{-\frac{R^2}{5T}},\ \sup_{s\geq R^2} s^{\frac{n}{2}} e^{-\frac{s}{5}}\}= 0 \label{3.2.5}
\end{align}

Combining (\ref{3.2.4}), (\ref{3.2.5}) with Lemma \ref{lem 3.1}, we have
\begin{align}
\Vert H_{R, i}(\cdot ,y, t)- H_i(\cdot , y, t) \Vert_{L^1(\nu_i)}\leq \epsilon (n, T, R) \label{3.2.6}
\end{align}
and $\displaystyle \mathop{\lim}_{R\rightarrow \infty} \epsilon (n, T, R)= 0$.
}
\qed

By Lemma \ref{lem 3.3}, Lemma \ref{lem 3.5} and Proposition \ref{prop 3.7}, we can assume, after passing to a subsequence of $\{i\}_{i= 1}^{\infty}$, that for every $j$, eigenvalue and eigenfunction converge:
\begin{align}
\lim_{i\rightarrow \infty}\lambda_{j, i}^{(R)}= \lambda_{j, \infty}^{(R)}\ , \quad \lim_{i\rightarrow \infty}\phi_{j, i}^{(R)}= \phi_{j, \infty}^{(R)} \label{3.8}
\end{align}

\begin{theorem}\label{thm 5.8}
{\begin{align}
H_{R, \infty}(x, y, t)\doteqdot \sum_{j= 1}^{\infty} e^{-\lambda_{j, \infty}^{(R)}t}\phi_{j, \infty}^{(R)}(x)\phi_{j, \infty}^{(R)}(y) \label{3.8.1}
\end{align}
is well defined on $B_{\infty}(R)\times B_{\infty}(R)\times (0, \infty)$, where $R> 2$. And 
\begin{align}
\lim_{i\rightarrow \infty}H_{R, i}(\cdot, y, t)= H_{R, \infty}(\cdot, y, t) \label{3.8.2}
\end{align}
where the convergence is in $L^2\big(B_{\infty}(R)\big)$, and is also locally uniform on $\mathring{B}_{\infty}(R)$. Furthermore, $H_{R, \infty}(\cdot, y, t)$ is locally Lipschitz on $\mathring{B}_{\infty}(R)$.
}
\end{theorem}

\begin{remark}\label{rem 3.9}
$H_{R, \infty}$ may depend on the choice of subsequence of $\{M_i^n\}$.
\end{remark}

\pf
{By (\ref{3.5.5}) and (\ref{3.8}), we get
\begin{align}
\Vert \phi_{j, \infty}^{(R)}\Vert_{L^{\infty}(\nu_{\infty})}\leq C(n) R^{\frac{n}{2}} \Big(\lambda_{j, \infty}^{(R)}\Big)^{\frac{n}{4}} \label{3.8.3}
\end{align}

Using (\ref{3.8.3}), when $t\in [t_0, \infty)$, $t_0> 0$ is any positive constant, we can obtain
\begin{align}
\Vert e^{-\lambda_{j, \infty}^{(R)}t} \phi_{j, \infty}^{(R)} (x) \phi_{j, \infty}^{(R)} (y) \Vert_{L^{\infty}(\nu_{\infty})} &\leq e^{-\lambda_{j, \infty}^{(R)} t} \Vert \phi_{j, \infty}^{(R)} \Vert_{L^{\infty}(\nu_{\infty})}^2  \nonumber \\
&\leq C(n, R) e^{-\lambda_{j, \infty}^{(R)}t} \Big(\lambda_{j, \infty}^{(R)}\Big)^{\frac{n}{2}} \leq C(n, R, t_0)e^{-\frac{\lambda_{j, \infty}^{(R)}t}{2}} \label{3.8.4}
\end{align}

Applying (\ref{3.8}) and Lemma \ref{lem 3.3}, we conclude that 
\begin{align}
\sum_{j= 1}^{\infty} \Big| e^{-\lambda_{j, \infty}^{(R)}t} \phi_{j, \infty}^{(R)} (x) \phi_{j, \infty}^{(R)} (y) \Big|\leq C(n, R, t_0) \sum_{j= 1}^{\infty} e^{-\big[C(n, R)j^{\frac{1}{n}} t\big]} \nonumber
\end{align}
which clearly converges uniformly on $B_{\infty}(R)\times B_{\infty}(R)\times [t_0, \infty)$ for any $t_0> 0$. Hence the kernel $H_{R, \infty}(x, y, t)$ is well defined and locally Lipschitz on $\mathring{B}_{\infty}(R)$.

Similar as (\ref{3.8.4}), it is easy to see 
\begin{align}
\Vert e^{-\lambda_{j, i}^{(R)}t} \phi_{j, i}^{(R)} (x) \phi_{j, i}^{(R)} (y) \Vert_{L^{\infty}(\nu_{i})}\leq C(n, R, t_0)e^{-\big[C(n, R)j^{\frac{1}{n}} t\big]} \label{3.8.5}
\end{align}
when $t\in [t_0, \infty)$. Then (\ref{3.8.2}) follows from (\ref{H(R, i)}), (\ref{3.8.1}),(\ref{3.8.4}), (\ref{3.8.5}) and Proposition \ref{prop 3.7}.
}
\qed

Fix one increasing sequence $R_{k}\rightarrow \infty$, by a diagonal argument, we can choose one subsequence of $\{M_i^n\}$, also denoted as $\{M_i^n\}$, such that for each $k$, $H_{R_k, i}\rightarrow H_{R_k, \infty}$ in $L^2\big(B_{\infty}(R_k)\big)$ and also locally uniform on $\mathring{B}_{\infty}(R_k)$.

On $M_i^n$, for $R_j< R_k$, we have 
\begin{align}
H_{R_j, i}(x, y, t)\leq H_{R_k, i}(x, y, t)\leq H_i(x, y, t)\leq \frac{C(n)}{\nu_i\big(\sqrt{t}\big)} e^{-\frac{d_{\rho_i}^2(x, y)}{5t}} \label{3.8.6}
\end{align}  
where $\nu_i(\sqrt{t})\doteqdot \nu_i\big(B_i(\sqrt{t})\big)$. Taking $i\rightarrow \infty$ in (\ref{3.8.6}), we get
\begin{align}
0\leq H_{R_j, \infty}(x, y, t)\leq H_{R_k, \infty}(x, y, t)\leq \frac{C(n)}{\nu_{\infty}\big(\sqrt{t}\big)} e^{-\frac{d_{\rho_{\infty}}^2(x, y)}{5t}} \label{3.8.7}
\end{align}  
where $\nu_{\infty}(\sqrt{t})= \nu_{\infty}\big(B_{\infty}(\sqrt{t})\big)$.
Thus we can get that the non-decreasing sequence $H_{R_j, \infty}$ converges pointwise to some function $H_{\infty}$:
\begin{align}
H_{\infty}(x, y, t)= \lim_{k\rightarrow \infty}H_{R_k, \infty}(x, y, t)= \lim_{k\rightarrow \infty} \lim_{i\rightarrow \infty} H_{R_k, i}(x_i, y, t) \label{3.8.8}
\end{align}
for some subsequence of $\{M_i^n\}_{i= 1}^{\infty}$, $\{R_k\}_{k= 1}^{\infty}$ and any $x_i\rightarrow x$.

\begin{prop}\label{prop 3.11}
{$H_{R, \infty}$ is a Dirichlet solution of the heat equation
\begin{equation}\label{3.11.1}
\left\{
\begin{array}{rl}
\Big(\frac{\partial}{\partial t}- \Delta_{\infty}\Big) H_{R, \infty}&= 0\\
\displaystyle \lim_{t\rightarrow 0}H_{R, \infty}(x, y, t)&= \delta_{y}(x)
\end{array} \right.
\end{equation}
on $B_{\infty}(R)\subset \Big(M_{\infty}, y, \rho_{\infty}, \nu_{\infty}\Big)$.
}
\end{prop}

\pf
{By Lemma \ref{lem 3.3}, \ref{lem 3.5}, \ref{lem 3.6} and Proposition \ref{prop 3.7}, we have
\begin{align}
\lim_{k\rightarrow \infty}\sum_{j= k}^{\infty} \Big|e^{-\lambda_j t} \big[d \phi_{j, \infty}(x)\big] \phi_{j, \infty}(y)\Big|= 0\ , \quad \quad x\in \mathring{B}_{\infty}(R)\label{3.11.2}
\end{align}
Hence $H_{R, \infty}$ is a Dirichlet solution of the heat equation by directly checking that (\ref{solution 2}) holds for it. 

From (\ref{3.8.2}), (\ref{3.8.7}) and the definition of $H_{R, i}$, using the similar argument as in Lemma \ref{lem 3.1}, we get 
\begin{align}
\lim_{t\rightarrow 0} \int_{M_{\infty}} H_{R, \infty}(x, y, t)f(x) d\nu_{\infty}(x)= f(y) \label{3.11.3}
\end{align}
where $f$ is any Lipschitz function with compact support on $M_{\infty}$.
}
\qed

\begin{prop}\label{prop 3.13}
{\begin{align}
\lim_{R\rightarrow \infty} H_{R, \infty}(\cdot , y, t)= p_{\infty}(\cdot , y, t) \label{3.13.1}
\end{align}
The convergence is uniform on $x\in M_{\infty}$, and uniform in $L^1(\nu_{\infty})$ on any finite time interval $(0, T]$.
}
\end{prop}

\pf
{Assume $R\geq 1$, put 
\begin{align}
M_{R, \infty}\doteqdot \sup\{p_{\infty}(x, y, t)|\ x\in \partial B_{\infty}(R),\ 0< t\leq T \} \label{3.13.2}
\end{align}
By (\ref{6.11.1}) and Property ($\mathscr{B}$) on $M_{\infty}$(from Proposition \ref{prop 3.9}), we have
\begin{align}
M_{R, \infty}&\leq \sup_{0< t\leq T} C(n) \nu_{\infty}\Big(B_{\infty}(\sqrt{t})\Big)^{-1} e^{-\frac{R^2}{5t}} \nonumber \\
&\leq C\cdot \max\Big\{e^{-\frac{R^2}{5T}},\ \sup_{0< t\leq 1} \nu_{\infty}\big(B_{\infty}(\sqrt{t})\big) e^{-\frac{R^2}{5t}} \Big\} \nonumber\\
&\leq C\cdot \max\{e^{-\frac{R^2}{5T}},\ \sup_{0< t\leq 1} t^{-\frac{n}{2}} e^{-\frac{R^2}{5t}} \} \label{3.13.3}\\
&\leq C(n)\max\{e^{-\frac{R^2}{5T}},\ R^{-n}\} \label{3.13.4}
\end{align}

From Proposition \ref{prop 3.11} and comparison inequalities for heat kernels on metric measure spaces (see Proposition $4.1$ in \cite{GHL}), we get 
\begin{align}
p_{\infty}(x, y, t)- M_{R, \infty}\leq H_{R, \infty}(x, y, t)\leq p_{\infty}(x, y, t) \label{3.13.5}
\end{align}

From (\ref{3.13.4}) and (\ref{3.13.5}), $\displaystyle \lim_{R\rightarrow \infty} H_{R, \infty}(\cdot, y, t)= p_{\infty}(\cdot, y, t)$ uniformly on $B_{\infty}(R)$. Combining with (\ref{6.11.1}), the convergence is uniform on $(0, T]\times M_{\infty}$.

From (\ref{3.13.3}) and Property ($\mathscr{B}$), note $\nu_{\infty}\big(B_{\infty}(1)\big)= 1$, we get
\begin{align}
\lim_{R\rightarrow \infty} M_{R, \infty}\nu_{\infty}\Big(B_{\infty}(R)\Big)&\leq \lim_{R\rightarrow \infty}C(n)R^{n} \cdot \max\{e^{\frac{R^2}{5T}},\ \sup_{0< t\leq 1} t^{-\frac{n}{2}} e^{-\frac{R^2}{5t}}\} \nonumber \\
&\leq \lim_{R\rightarrow \infty} C(n) \max\{R^n e^{-\frac{R^2}{5T}},\ \sup_{s\geq R^2} s^{\frac{n}{2}} e^{-\frac{s}{5}}\}= 0 \label{3.13.6}
\end{align}

Combining (\ref{3.13.6}) with Corollary \ref{cor 3.10}, we have
\begin{align}
\Vert H_{R, \infty}(\cdot ,y, t)- p_{\infty}(\cdot , y, t) \Vert_{L^1(\nu_{\infty})}\leq \epsilon (n, T, R) \label{3.13.7}
\end{align}
and $\displaystyle \mathop{\lim}_{R\rightarrow \infty} \epsilon (n, T, R)= 0$.
}
\qed

\begin{prop}\label{prop 3.14}
{Assume $x_i\rightarrow x$ as $(M_i, y, \rho_i, \nu_i)\stackrel{d_{GH}}{\longrightarrow} (M_{\infty}, y_{\infty},\rho_{\infty}, \nu_{\infty})$, then
\begin{align}
\lim_{i\rightarrow \infty} H_{i}(x_i, y, t)= p_{\infty}(x, y, t)\, \quad\quad t\in (0, \infty) \label{3.14.1}
\end{align}
The convergence is locally uniform on $M_{\infty}$.
}
\end{prop}

\begin{remark}\label{rem 9.8}
{$H_{\infty}$ in (\ref{3.8.8}) is equal to $p_{\infty}$ in (\ref{3.14.1}).
}
\end{remark}

\pf
{For any sequence $(M_i, y, \rho_i, \nu_i)\stackrel{d_{GH}}{\longrightarrow} (M_{\infty}, y_{\infty},\rho_{\infty}, \nu_{\infty})$, we can get a subsequence of $\{M_i^n\}$ as before, also denoted as $\{M_i^n\}$, such that, there exists increasing sequence $R_k\rightarrow \infty$, and 
\begin{align}
\lim_{i\rightarrow \infty} H_{R_k, i}(\cdot, y, t)= H_{R_k, \infty}(\cdot, y, t)\, \quad \quad k= 1, 2, 3, \cdots \nonumber
\end{align}
where the convergence is as in Theorem \ref{thm 5.8}. 

Then
\begin{align}
&|H_i(x_i, y, t)- p_{\infty}(x, y, t)|\leq \Big(|H_i(x_i, y, t)- H_{R_k, i}(x_i, y, t)|\nonumber \\
&\quad + |H_{R_k, \infty}(x, y, t)- p_{\infty}(x, y, t)|\Big) + |H_{R_k, i}(x_i, y, t)- H_{R_k, \infty}(x, y, t)| \label{3.14.2}
\end{align}

For any $\epsilon> 0$, from Proposition \ref{prop 3.2} and \ref{prop 3.13}, we get the first two terms on the right side of  (\ref{3.14.2}) will be less than $\frac{1}{3}\epsilon$ when $k$ is big enough. Now fixed $k$ such that $x\in \mathring{B}_{\infty}(R_k)$ and 
\[\Big(|H_i(x_i, y, t)- H_{R_k, i}(x_i, y, t)|+ |H_{R_k, \infty}(x, y, t)- p_{\infty}(x, y, t)|\Big)< \frac{2}{3}\epsilon \]
Using Theorem \ref{thm 5.8}, if $i$ is big enough (which may depend on $k$ we chose above), then we get 
\[|H_{R_k, i}(x_i, y, t)- H_{R_k, \infty}(x, y, t)|< \frac{1}{3}\epsilon\]

By the above argument, we get that for such subsequence of $\{M_i^n\}$, 
\begin{align}
\lim_{i\rightarrow \infty}H_i(x_i, y, t)= p_{\infty}(x, y, t) \nonumber
\end{align}

However, any subsequence of $\{M_i^n\}$ must contain a subsequence whose limit is also $p_{\infty}$ by the above argument. Hence, in fact we prove that for the original sequence $\{M_i^n\}$, (\ref{3.14.1}) holds.
}
\qed

{\it \textbf{Proof of Theorem \ref{thm 8.13}}}:~
{From (\ref{Hi}), (\ref{3.14.1}) and $x\rightarrow y_{\infty}$ as $i\rightarrow \infty$ for any $x\in M_i$.
}
\qed

\section{Analysis on manifolds with cone structures at infinity}\label{SECTION 08}

In this section we will discuss large time behavior of the heat kernel on manifolds with cone structures at infinity (see Definition below), and prove Theorem \ref{thm 11.0}.

\begin{definition}\label{def cone}
{Assume that $(M^n, g)$ is a complete manifold with $Rc\geq 0$, $y$ is some fixed point on $M^n$, and for any $t_i\rightarrow \infty$, define $(M_i, y, \rho_i, \nu_i)$ as in \textbf{Blow Down Setup}, such that 
\begin{align}
(M_i, y, \rho_i, \nu_i)\stackrel{d_{GH}}{\longrightarrow} (M_{\infty}, y_{\infty}, \rho_{\infty}, \nu_{\infty}) \label{11.0.1} 
\end{align}
If $(M_{\infty}, y_{\infty}, \rho_{\infty})$ (may be different for different choice of $\{t_i\}$) always has the cone structure, i.e. 
\begin{align}
\rho_{\infty}= dr^2+ l(r)^2 dX \label{cone}
\end{align}
where $X$ is some compact metric space, $l(r)> 0$ is some function of $r$, then we say that $M^n$ is a \textbf{manifold with cone structures at infinity}.
}
\end{definition}

{\it \textbf{Proof of Theorem \ref{thm 11.0}}}:~
{Assume that $s_i\rightarrow \infty$, blowing down the metric $g$ by $s_i^{-1}$ instead of $t_i^{-1}$, define $(M_i, y, \rho_i, \nu_i)$ as in \textbf{Blow Down Setup}, and the following holds:
\begin{align}
(M_i, y, \rho_i, \nu_i)\stackrel{d_{GH}}{\longrightarrow} (M_{\infty}, y_{\infty}, \rho_{\infty}, \nu_{\infty}) \nonumber
\end{align}
From (\ref{3.0}) and (\ref{vol}), it is easy to get $\nu_{\infty}(B_{\infty}(y_{\infty}, r))= h(r)$. By the assumption that $M^n$ is a complete manifold with cone structures at infinity, we 
get that the heat kernel $p_{\infty}$ on $(M_{\infty}, y_{\infty}, \rho_{\infty}, \nu_{\infty})$, only depends on $r= \rho (x, y_{\infty})$ and $t$, denoted as $p_{\infty}(r, t)$. 

It is easy to get 
\begin{align}
\Delta p_{\infty}(r, t)= \frac{\partial^2 p_{\infty}}{\partial r^2}+ \Big(\frac{h''(r)}{h'(r)}\Big)\cdot \frac{\partial p_{\infty}}{\partial r} \nonumber
\end{align}

Hence $p_{\infty}(r, t)$ is the unique positive solution of 
\begin{equation}\nonumber
\left\{
\begin{array}{rl}
\frac{\partial p_{\infty}}{\partial t}&= (p_{\infty})_{rr}+ \big(\frac{h''}{h'}\big) (p_{\infty})_r  \\
\displaystyle\lim_{t\rightarrow 0}p_{\infty}(r, t)&= \delta_{y_{\infty}}(x) \\
\end{array} \right.
\end{equation}
From the above, it is easy to see that $p_{\infty}(r, t)$ is uniquely determined by $(\frac{h''}{h'})(r)$. The conclusion follows from Theorem \ref{thm 8.13}, the above argument and (\ref{equ vol}).
}
\qed

\begin{remark}\label{rem 11.1}
{Note (\ref{equ vol}) is equivalent to the assumption that $\frac{h'(r)}{\tilde{h}'(r)}$ is a constant independent of $r$. Although the tangent cones at infinity of manifold $M^n$ may be different metric measure spaces for different choices of $s_i$, $p_{\infty}$ only depends on the function $h'(r)$ when the tangent cone at infinity $(M_{\infty}, y_{\infty}, \rho_{\infty})$ has the cone structure as in (\ref{cone}).
}
\end{remark}

\begin{theorem}\label{thm 11.1}
{Assume that $(M^n, g)$ is a complete manifold with nonnegative sectional curvature, $n\geq 3$, $y$ is some fixed point on $M^n$, and for any $r> 0$,
\begin{align}
\lim_{s\rightarrow \infty}\frac{V_y(sr)}{V_y(s)}= h(r) \label{11.1.0}
\end{align}
where $h(r)> 0$ is some positive function. 

Then there exists a unique $(M_{\infty}, y_{\infty}, \rho_{\infty}, \nu_{\infty})$, where $\nu_{\infty}(B_{\infty}(y_{\infty}, r))= h(r)$, such that for any $t_i\rightarrow \infty$, define $(M_i, y, \rho_i, \nu_i)$ as in \textbf{Blow Down Setup}, we have
\begin{align}
(M_i, y, \rho_i, \nu_i)&\stackrel{d_{GH}}{\longrightarrow} (M_{\infty}, y_{\infty}, \rho_{\infty}, \nu_{\infty}) \label{11.1.1} \\
\lim_{t\rightarrow \infty} V_y(\sqrt{t})\cdot H(x, y, t)&= p_{\infty}(y_{\infty}, y_{\infty}, 1) \label{11.1.2}
\end{align}
where $p_{\infty}$ is the heat kernel on $(M_{\infty}, y_{\infty}, \rho_{\infty}, \nu_{\infty})$.
}
\end{theorem}

\pf
{Because $M^n$ has non-negative sectional curvature, from Theorem $I. 26$ in \cite{RFTA}, we know that the tangent cone at infinity $(M_{\infty}, y_{\infty}, \rho_{\infty})$ is the unique metric cone. Hence $M^n$ is a manifold with cone structures at infinity and (\ref{11.1.1}) is obtained.

From the assumption (\ref{11.1.0}) and the above argument, we can apply Theorem \ref{thm 11.0}, (\ref{11.1.2}) is obtained.
}
\qed

As an application of the above theorem, we have the following interesting result about nonnegatively curved manifolds with asymptotic polynomial volume growth.

\begin{cor}\label{cor 11.2}
{Assume that $(M^n, g)$ is a complete manifold with nonnegative sectional curvature, $n\geq 3$ and it has asymptotic polynomial volume growth, i.e.
\begin{align}
\lim_{r\rightarrow \infty} \frac{V(r)}{r^k}= C_0 \nonumber
\end{align}
where $k\geq 1$ and $C_0> 0$ are constants. Then (\ref{11.1.2}) holds.
}
\end{cor} 

\pf
{The proof follows directly from Theorem \ref{thm 11.1}.
}
\qed




\section{Example with $\varliminf_{t\rightarrow \infty}V(\sqrt{t})H(x, y, t)< \varlimsup_{t\rightarrow \infty}V(\sqrt{t})H(x, y, t)$}\label{SECTION 09}

In this section we will construct the first example, which is a complete manifold with nonnegative Ricci curvature and $\varliminf_{t\rightarrow \infty}V(\sqrt{t})H(x, y, t)< \varlimsup_{t\rightarrow \infty}V(\sqrt{t})H(x, y, t)$. 

From Theorem \ref{thm 8.13}, the example should have different tangent cones at infinity of the manifold with renormalized measure. Furthermore, from Theorem \ref{thm 11.0} and its proof, if two tangent cones at infinity of $(M^n, g)$ have the cone structure as defined in Definition \ref{def cone}, \textbf{only different renormalized measure} will result in the inconsistent limit behavior of heat kernel. Note in this context, if there exists $r> 0$, such that $\nu_{\infty}\big(B_{\infty}(r)\big)\neq \tilde{\nu}_{\infty}\big(\tilde{B}_{\infty}(r)\big)$, where $B_{\infty}(r)\subset M_{\infty}= C(X)$ and $\tilde{B}_{\infty}(r)\subset \tilde{M}_{\infty}= C(\tilde{X})$ are two balls with the same radius $r$ in different metric tangent cones $C(X)$, $C(\tilde{X})$; we say that the renormalized measures $\nu_{\infty}$, $\tilde{\nu}_{\infty}$ are different.

Hence, the different structure of tangent cones at infinity alone can not guarantee the inequality the inconsistent limit behavior of heat kernel. As mentioned in the introduction of this paper, Perelman (\cite{Pere}) had constructed the manifold with $Rc\geq 0$, maximal volume growth and quadratic curvature decay, where the tangent cone at infinity is not unique. However it is not hard to see that the renormalized measure on those different tangent cones (in fact, metric cones) are the same, so will not lead to inconsistent limit behavior of heat kernel on such manifolds.

In fact, from Theorem \ref{thm Li}, the example manifold must be collapsing case. The construction of the following example is inspired by the related discussion in Section $8$ of \cite{CC1}. However, we need to do some suitable modifications to assure the different renormalized measure on different tangent cones at infinity of manifold.

Let us start from the generalized Hopf fibration of $\mathbb{S}^7$ as the following:
\begin{align}
\mathbb{S}^3\longrightarrow \mathbb{S}^7\stackrel{\pi}{\longrightarrow} \mathbb{S}^4\ , \quad g^{\mathbb{S}^7}= k_1+ k_2 \nonumber
\end{align}
where $\mathbb{S}^3$, $\mathbb{S}^7$, $\mathbb{S}^4$ carry the metrics $g^{\mathbb{S}^3}$, $g^{\mathbb{S}^7}$, $\frac{1}{4}g^{\mathbb{S}^4}$; $\pi$ is a Riemannian submersion with totally geodesic fibers and $k_1= g^{\mathbb{S}^3}$, $k_2= \pi^{*}\big(\frac{1}{4}g^{\mathbb{S}^4}\big)$; $g^{\mathbb{S}^n}$ denotes the canonical metric of curvature $\equiv 1$ on $\mathbb{S}^n$.

Define $\tilde{g}= f^2 k_1+ h^2 k_2$, then the following formulas are well-known (for example, see Section $2$ in \cite{BKN}):
\begin{align}
Rc(\tilde{g})|_{k_1}= \Big(\frac{2}{f^2}+ \frac{4f^2}{h^4}\Big)I\ , \quad Rc(\tilde{g})|_{k_2}= \frac{6(2h^2- f^2)}{h^4}I \label{8.0}
\end{align}
Other mixed $Rc(\tilde{g})= 0$.

Then for metric $g= dr^2+ f^2(r)k_1+ h^2(r)k_2$ on $M^8$, which is diffeomorphic to $\mathbb{R}^8$, from $(8.13)$ in \cite{CC1} and (\ref{8.0}), we have 
\begin{align}
Rc(g)|_{k_1}&= \frac{2\big(1-(f')^2\big)}{f^2}- \frac{f''}{f}+ \frac{4f^2}{h^4}- 4\frac{f'h'}{fh} \label{8.1} \\
Rc(g)|_{k_2}&= \frac{6(2h^2- f^2)}{h^4}- \frac{h''}{h}- 3\frac{(h')^2}{h^2}- 3\frac{f'h'}{fh} \label{8.2} \\
Rc(g)(\vec{n}, \vec{n})&= -\Big[3\frac{f''}{f}+ 4\frac{h''}{h}\Big] \label{8.3} 
\end{align}

Our construction will be broken into four steps in subsections \ref{subsection 8.1}-\ref{subsection 8.4} separately, we will verify that our example $(M^8, g)$ has the property $\displaystyle \varliminf_{t\rightarrow \infty}H(x, y, t)< \varlimsup_{t\rightarrow \infty}H(x, y, t)$ in subsection \ref{subsection 8.4}.

\subsection{Step $(\mathbf{I})$}\label{subsection 8.1}
$\newline$Initial approximation $\bar{f}$, $\bar{h}$ to the functions $f$, $h$ will be constructed inductively at this stage. These approximations have jump discontinuities at the points $b_i$; see (\ref{tau1}), (\ref{tau2}), (\ref{delta1}), (\ref{delta2}). However, the left- and right-hand limits of the first derivatives do agree at all $b_i$, $i\geq 1$, see (\ref{1st der 1}) and (\ref{1st der 2}).

We can define $\bar{f}(r)$ as the following:
\begin{equation}\label{f1}
\bar{f}(r)= \left\{
\begin{array}{rl}
\beta_{2i}b_{2i+ 1}^{-\omega_{2i}} r^{1- \eta_1} \quad \quad r\in (b_{2i}, b_{2i+ 1}] \\
\beta_{2i+ 1}b_{2i+ 2}^{\omega_{2i+ 1}} r^{1- \eta_2} \quad \quad  r\in (b_{2i+ 1}, b_{2i+ 2}] \\
\end{array} \right.
\end{equation}
where for $i= 0, 1, 2, \cdots$ we have:
\begin{assu}\label{Assu 1}
\begin{align}
&1- \epsilon_0> \frac{1- \eta_2}{1- \eta_1}\geq \frac{99}{100}\ , \quad 1> \eta_2> \eta_1> \frac{1}{2}(1+ \epsilon_0) \label{8.1.0} \\
&\beta_{2i+ 2}> \beta_{2i}> 0\ , \quad  \beta_{2i+ 1}> \beta_{2i+ 3}> 0\ , \quad \beta_0\geq \frac{99}{100}\ , \quad \beta_1\leq \frac{1}{100}  \label{8.1.1}\\
&\lim_{i\rightarrow \infty}\beta_{2i}= 1\ , \quad \lim_{i\rightarrow \infty}\beta_{2i+ 1}= 0 \label{8.1.1.5}
\end{align}
\end{assu}
$\beta_i$, $\omega_i$, $\epsilon_0$ are positive constants to be determined later and satisfy
\begin{assu}\label{Assu 2}
{$\displaystyle
\lim_{i\rightarrow \infty}\omega_i= 0\ , \quad \frac{\eta_2- \eta_1}{100}> \omega_0> \omega_1> \omega_2> \cdots $
}
\end{assu}
We have the following equations for $i= 0, 1, 2, \cdots$
\begin{assu}\label{Assu 3}
\begin{align}
\frac{1-\eta_1}{1- \eta_2}\cdot \frac{\beta_{2i}}{\beta_{2i+ 1}} &= \frac{b_{2i+ 2}^{\omega_{2i+ 1}}}{b_{2i+ 1}^{\eta_2- \eta_1- \omega_{2i}}}  \label{8.1.3} \\
\frac{1-\eta_2}{1- \eta_1}\cdot \frac{\beta_{2i+ 1}}{\beta_{2i+ 2}} &= \frac{b_{2i+ 2}^{\eta_2- \eta_1- \omega_{2i+ 1}}}{b_{2i+ 3}^{\omega_{2i+ 2}}} \label{8.1.4}
\end{align}
\end{assu}
which implies that for $i\geq 1$, 
\begin{align}
\bar{f}'(b_{i})= \lim_{r\rightarrow b_{i}^{+}}\bar{f}'(r)= \lim_{r\rightarrow b_{i}^{-}}\bar{f}'(r) \label{1st der 1}
\end{align}

We define $\bar{h}(r)$ in the following way:
\begin{equation}\label{h1}
\bar{h}(r)= \left\{
\begin{array}{rl}
\alpha_{2i} b_{2i+ 1}^{-\epsilon_{2i}} r^{1+ \epsilon_{2i}} \quad \quad r\in (b_{2i}, b_{2i+ 1}] \\
\alpha_{2i+ 1} b_{2i+ 2}^{\epsilon_{2i+ 1}} r^{1- \epsilon_{2i+ 1}} \quad \quad  r\in (b_{2i+ 1}, b_{2i+ 2}] \\
\end{array} \right.
\end{equation}
where for $i= 0, 1, 2, \cdots$, we have
\begin{assu}\label{Assu 4}
{
\begin{align}
&\alpha_{2i+ 2}> \alpha_{2i}> 0\ , \quad \alpha_{2i+ 1}> \alpha_{2i+ 3}> 0  \label{8.4.0} \\ 
&\lim_{i\rightarrow \infty}\alpha_{2i}= 1\ , \quad \alpha_0\geq \frac{99}{100} \ , \quad \lim_{i\rightarrow \infty}\alpha_{2i+ 1}= 0\ , \quad \alpha_1\leq \frac{1}{100} \label{8.4}
\end{align}
}
\end{assu}
$\epsilon_i$, $b_i$ are to be determined later, and satisfies
\begin{assu}\label{Assu 5}
{\begin{align}
1< b_0< b_1< b_2< \cdots \ , \quad \lim_{i\rightarrow \infty}b_i= \infty \label{8.6} \\
1> \epsilon_0> \epsilon_1> \epsilon_2> \cdots \ , \quad \lim_{i\rightarrow \infty}\epsilon_i= 0 \label{8.7}
\end{align}
}
\end{assu}

We also have the following equations for $i= 0, 1, 2, \cdots $,
\begin{assu}\label{Assu 6}
{
\begin{align}
\frac{\alpha_{2i}}{\alpha_{2i+ 1}}&= \Big(\frac{1- \epsilon_{2i+ 1}}{1+ \epsilon_{2i}} \Big) \Big(\frac{b_{2i+ 2}}{b_{2i+ 1}}\Big)^{\epsilon_{2i+ 1}} \label{8.8} \\
\frac{\alpha_{2i+ 2}}{\alpha_{2i+ 1}}&= \Big(\frac{1- \epsilon_{2i+ 1}}{1+ \epsilon_{2i+ 2}} \Big) \Big(\frac{b_{2i+ 3}}{b_{2i+ 2}}\Big)^{\epsilon_{2i+ 2}}  \label{8.9}
\end{align}
}
\end{assu}
which implies that for $i\geq 1$,
\begin{align}
\bar{h}'(b_{i})= \lim_{r\rightarrow b_{i}^{+}}\bar{h}'(r)= \lim_{r\rightarrow b_{i}^{-}}\bar{h}'(r) \label{1st der 2}
\end{align}

In the rest part of $\mathbf{Step\ (I)}$, we will prove $Rc(\bar{f}, \bar{h})> 0$ on $(b_0, \infty)$ except the points $b_i$, $i= 1, 2, \cdots$.

$(\mathbf{i})$. We firstly consider the interval $(b_0, b_1]$: 

If we assume that
\begin{assu}\label{Assu 7}
{$\displaystyle b_0^{\eta_1}\geq 7 $
}
\end{assu}
then 
\begin{align}
Rc|_{k_1}(\bar{f}, \bar{h})\geq \frac{2}{r^2}\big[\beta_0^{-2}b_1^{2\omega_0}r^{2\eta_1}- (1- \eta_1)(3+ 2\epsilon_0- \eta_1)\big]> 0 \label{8.10}
\end{align}

If we assume that
\begin{assu}\label{Assu 8}
{$\displaystyle \alpha_0 b_0^{\eta_1}> \Big(\frac{b_1}{b_0}\Big)^{\epsilon_0} $
}
\end{assu}
then $\bar{h}(r)> \bar{f}(r)$. And if we further assume that 
\begin{assu}\label{Assu 9}
{$\displaystyle \epsilon_{2i}< \frac{1}{2}(\alpha_{2i}^{-1}- 1)\ , \quad i= 0, 1, 2, \cdots $
}
\end{assu}
we have
\begin{align}
Rc|_{k_2}(\bar{f}, \bar{h})&> \frac{6}{\bar{h}^2}- \frac{\epsilon_0(1+ \epsilon_0)}{r^2}- \frac{3(1+ \epsilon_0)^2}{r^2}- \frac{3(1- \eta_1)(1+ \epsilon_0)}{r^2} \nonumber \\
& > \frac{6}{r^2}\big[\alpha_0^{-2}- (1+ 2\epsilon_0)^2\big]\geq 0 \label{8.11}
\end{align}

If we assume 
\begin{assu}\label{Assu 10}
{$\displaystyle \epsilon_0< \frac{1}{4}\eta_1(1- \eta_1)b_0^{-\eta_1}b_1^{-\omega_0} $
}
\end{assu}
then 
\begin{align}
Rc(\vec{n}, \vec{n})= \frac{3\eta_1 (1- \eta_1)}{r^2}- 4\frac{\epsilon_0(1+ \epsilon_0)}{r^2}> 0 \label{8.12}
\end{align}

From (\ref{8.10}), (\ref{8.11}) and (\ref{8.12}), we get $Rc(\bar{f}, \bar{h})> 0$ on $(b_0, b_1)$.

$(\mathbf{ii})$. Next we consider the interval $(b_{2i}, b_{2i+ 1}]$, $i\geq 1$. 

From Assumption \ref{Assu 7} and (\ref{8.1}), it is easy to get 
\begin{align}
Rc|_{k_1}> 0 \label{8.13}
\end{align}

If we assume that 
\begin{assu}\label{Assu 11}
{$\displaystyle b_{2i}^{\eta_1}> \alpha_{2i- 1}^{- 1} \Big(\frac{1+ \epsilon_{2i}}{1- \epsilon_{2i- 1}}\Big) \ , \quad i= 1, 2, \cdots $
}
\end{assu}
then $\bar{h}(r)> \bar{f}(r)$, from it and Assumption \ref{Assu 9}, we get
\begin{align}
Rc|_{k_2}> \frac{6}{\bar{h}^2}- \frac{\epsilon_{2i}(1+ \epsilon_{2i})}{r^2}- \frac{3(1+ \epsilon_{2i})^2}{r^2}- \frac{3(1- \eta_1)(1+ \epsilon_{2i})}{r^2}> 0 \label{8.14}
\end{align}

Similarly, from Assumption \ref{Assu 10}, we get 
\begin{align}
Rc(\vec{n}, \vec{n})> 0 \label{8.15}
\end{align}

From (\ref{8.13}), (\ref{8.14}) and (\ref{8.15}), we get that $Rc(\bar{f}, \bar{h})> 0$ on $(b_{2i}, b_{2i+ 1})$, where $i\geq 1$.

$(\mathbf{iii})$. Finally, we consider the interval $(b_{2i+1}, b_{2i+ 2}]$, $i\geq 0$.

From Assumption \ref{Assu 7}, (\ref{f1}) and (\ref{h1}), it is easy to get
\begin{align} 
Rc|_{k_1}\geq \frac{2}{r^2}\Big[\big(\frac{1- \eta_2}{1- \eta_1}\big)^2\beta_{2i}^{-2}b_{2i+ 1}^{2\eta_1+ 2\omega_{2i}}- 3\Big]> 0 \label{8.16}
\end{align}

From Assumption \ref{Assu 7}, it is easy to get that $\bar{h}(r)> \bar{f}(r)$. From (\ref{8.1.0}) and Assumption \ref{Assu 9}, we get
\begin{align}
Rc|_{k_2}> \frac{6}{\bar{h}^2}- \frac{3(1- \epsilon_{2i+ 1})^2}{r^2}- 3\frac{(1- \eta_2)(1- \epsilon_{2i+ 1})}{r^2}> 0 \label{8.17}
\end{align} 

From $\bar{f}''< 0$ and $\bar{h}''< 0$, it is easy to get
\begin{align}
Rc(\vec{n}, \vec{n})> 0 \label{8.18}
\end{align}

From (\ref{8.16}), (\ref{8.17}) and (\ref{8.18}), we get that $Rc(\bar{f}, \bar{h})> 0$ on $(b_{2i+ 1}, b_{2i+ 2})$, where $i\geq 0$.

From all the above, we get $Rc(\bar{f}, \bar{h})> 0$ on $(b_0, \infty)$ excepts the points $b_i$, where $i\geq 1$.

Note $\bar{f}$ has jump discontinuities at the points $b_j$, $j= 1, 2, \cdots$. For $i= 0, 1, 2, \cdots$, we have
\begin{align}
\tau_{2i+ 1}&\doteqdot \bar{f}(b_{2i+ 1})- \lim_{r\rightarrow b_{2i+ 1}^{+}}\bar{f}(r)= -\Big(\frac{\eta_2- \eta_1}{1- \eta_1}\Big)\beta_{2i+ 1}b_{2i+ 2}^{\omega_{2i+ 1}} b_{2i+ 1}^{1- \eta_2} \label{tau1} \\
\tau_{2i+ 2}&\doteqdot \bar{f}(b_{2i+ 2})- \lim_{r\rightarrow b_{2i+ 2}^{+}}\bar{f}(r)= \Big(\frac{\eta_2- \eta_1}{1- \eta_1}\Big)\beta_{2i+ 1} b_{2i+ 2}^{1- \eta_2+ \omega_{2i+ 1}}\label{tau2}
\end{align}

Similarly, $\bar{h}$ has jump discontinuities at the points $b_j$, $j= 1, 2, \cdots$. For $i= 0, 1, 2, \cdots$, we have
\begin{align}
\delta_{2i+ 1}&\doteqdot \bar{h}(b_{2i+ 1})- \lim_{r\rightarrow b_{2i+ 1}^{+}}\bar{h}(r)= -\alpha_{2i} b_{2i+ 1} \Big(\frac{\epsilon_{2i+ 1}+ \epsilon_{2i}}{1- \epsilon_{2i+ 1}}\Big) \label{delta1} \\
\delta_{2i+ 2}&\doteqdot \bar{h}(b_{2i+ 2})- \lim_{r\rightarrow b_{2i+ 2}^{+}}\bar{h}(r)= \alpha_{2i+ 1} b_{2i+ 2} \Big(\frac{\epsilon_{2i+ 2}+ \epsilon_{2i+ 1}}{1+ \epsilon_{2i+ 2}}\Big)\label{delta2}
\end{align}

\subsection{Step $(\mathbf{II})$}\label{subsection 8.2}
$\newline$ We construct $\tilde{f}$, $\tilde{h}$ on the interval $[0, b_0]$ in this step.

Define $\tilde{f}$, $\tilde{h}$ on $[0, b_0]$ as the following:
\begin{equation}\label{f2}
\tilde{f}(r)\doteqdot \left\{
\begin{array}{rl}
r \quad \quad r\in [0, \frac{b_0}{2}] \\
r- C_1\big(r- \frac{b_0}{2}\big)^2 \quad \quad  r\in (\frac{b_{0}}{2}, b_{0}] \\
\end{array} \right.
\end{equation}
where $C_1= \frac{1}{b_0}\Big[1- \beta_0 b_1^{-\omega_0}b_0^{-\eta_1}(1- \eta_1)\Big]> 0$.
\begin{equation}\label{h2}
\tilde{h}(r)\doteqdot \left\{
\begin{array}{rl}
r \quad \quad r\in [0, \frac{b_0}{2}] \\
r- C_2\big(r- \frac{b_0}{2}\big)^2 \quad \quad  r\in (\frac{b_{0}}{2}, b_{0}] \\
\end{array} \right.
\end{equation}
where $C_2= \frac{1}{b_0}\big[1- \alpha_0(1+ \epsilon_0)(\frac{b_0}{b_1})^{\epsilon_0}\big]> 0$.

Then $\tilde{f}(0)=\tilde{h}(0)= 0$, $\tilde{f}'(0)= \tilde{h}'(0)= 1$,
\begin{align}
\tilde{f}'(b_0)= \beta_0 b_1^{-\omega_0} b_0^{-\eta_1} (1- \eta_1) \ , \quad \tilde{h}'(b_0)= \alpha_0 (1+ \epsilon_0)(\frac{b_0}{b_1})^{\epsilon_0} \nonumber
\end{align}
On $(\frac{b_0}{2}, b_0)$, we have $\tilde{f}''(r)= -2C_1< 0$ and $\tilde{h}''(r)= -2C_2< 0$, hence 
\begin{align}
(1- \eta)b_0^{-\eta}\leq \tilde{f}'(r)\leq 1\ , \quad \alpha_0(1+ \epsilon_0)(\frac{b_0}{b_1})^{\epsilon_0}\leq \tilde{h}'(r)\leq 1 \label{8.19}
\end{align}

It is easy to see $\displaystyle \tilde{f}'(b_0)= \lim_{r\rightarrow b_0^{+}} \bar{f}'(r)$ and $\displaystyle \tilde{h}'(b_0)= \lim_{r\rightarrow b_0^{+}} \bar{h}'(r)$.

In the rest part of $\mathbf{Step\ (II)}$, we will show that $Rc(\tilde{f}, \tilde{h})> 0$ on $(\frac{b_0}{2}, b_0)$. 

It is obvious that $Rc(\tilde{f}, \tilde{h})= 0$ on $[0, \frac{b_0}{2})$.

Also it is easy to get that $Rc(\vec{n}, \vec{n})> 0$ from $\tilde{f}''< 0$ and $\tilde{h}''< 0$ on $(\frac{b_0}{2}, b_0)$.

Next if we assume that 
\begin{assu}\label{Assu 13}
{$\displaystyle \alpha_0 (1+ \epsilon_0) \big(\frac{b_0}{b_1}\big)^{\epsilon_0} = \frac{2}{3}+ \frac{1}{3}\beta_0 b_1^{-\omega_0} b_0^{-\eta_1}(1- \eta_1) $
}
\end{assu}
then $C_2= \frac{1}{3}C_1$, $\tilde{h}'(r)> \tilde{f}'(r)$ on $[\frac{b_0}{2}, b_0]$.
Hence $\tilde{h}(r)\geq \tilde{f}(r)$.

If we further have 
\begin{align}
\Big(\frac{\tilde{f}}{\tilde{h}}\Big)^3\geq \tilde{f}'\tilde{h}' \label{f/h}
\end{align}
then it is easy to get $Rc|_{k_1}> 0$ and $Rc|_{k_2}> 0$ on $(\frac{b_0}{2}, b_0)$. 

To show $\big(\frac{\tilde{f}}{\tilde{h}}\big)^3\geq \tilde{f}'\tilde{h}'$ on $[\frac{b_0}{2}, b_0]$, we consider the function
\[\varphi_1(r)\doteqdot \tilde{f}^3- \tilde{f}'\tilde{h}'\tilde{h}^3\]
Note $\varphi_1(\frac{b_0}{2})= 0$, we only need to show $\varphi_1'(r)\geq 0$ on $[\frac{b_0}{2}, b_0]$. It is easy to get
\begin{align}
\varphi_1'(r)\geq 3\tilde{f}'[\tilde{f}- \tilde{h}'h][\tilde{f}+ \tilde{h}'\tilde{h}] \nonumber
\end{align}
Hence we just need to show that $\tilde{f}- \tilde{h}'\tilde{h}\geq 0$ on $[\frac{b_0}{2}, b_0]$. Define 
\begin{align}
\varphi_2(r)\doteqdot \tilde{f}- \tilde{h}'\tilde{h} \nonumber 
\end{align}
Observe that $\varphi_2(\frac{b_0}{2})= 0$, the problem reduces to show that 
\begin{align}
\varphi_2'(r)= \tilde{f}'- (\tilde{h}')^2- \tilde{h}\tilde{h}'' \geq 0 \ , \quad r\in [\frac{b_0}{2}, b_0]
\end{align}
Let $\varphi_3(r)\doteqdot \tilde{f}'- (\tilde{h}')^2- \tilde{h}\tilde{h}''$, then 
\begin{align}
\varphi_3'(r)= 6C_2\tilde{h}'- 2C_1\leq 0 \label{8.20}
\end{align}
Now using Assumption \ref{Assu 13}, which is equivalent to $C_2= \frac{1}{3}C_1$, it is easy to get 
\begin{align}
\varphi_3(b_0)= b_0C_2(1- \frac{1}{2}b_0C_1)> 0 \label{8.21}
\end{align}
From (\ref{8.20}) and (\ref{8.21}), we get $\varphi_3(r)\geq 0$. Hence (\ref{f/h}) is obtained, we are done.

\subsection{Step $(\mathbf{III})$}\label{subsection 8.3}
$\newline$ By adjusting the values of the functions $\bar{f}$, $\bar{h}$, by suitable constants on each interval $(b_i, b_{i+ 1}]$, we can remove the jump discontinuities, thereby obtaining $C^1$ functions $\hat{f}$, $\hat{h}$ by gluing $\tilde{f}$, $\tilde{h}$ with $\bar{f}$, $\bar{h}$.

The functions $\hat{f}$, $\hat{h}$ may not have the second derivatives at the points $b_{i}$.

Now we define
\begin{equation}\label{f3}
\hat{f}(r)\doteqdot \left\{
\begin{array}{rl}
&\tilde{f}(r) \quad \quad \quad \quad  \quad  \quad r\in [0, b_0] \\
&\bar{f}(r)+ \sum_{l= 0}^{k} \tau_l \quad \quad   r\in (b_k, b_{k+ 1}]\ , \ k= 0, 1, 2, \cdots \\
\end{array} \right.
\end{equation}
where $\displaystyle \tau_0\doteqdot \tilde{f}(b_0)- \lim_{r\rightarrow b_0^{+}} \bar{f}(r)= \frac{b_0}{4}\Big[3- (3+ \eta_1)\beta_0 b_1^{-\omega_0} b_0^{-\eta_1}\Big]$, when $l\geq 1$, $\tau_l$ is defined in (\ref{tau1}) and (\ref{tau2}). From Assumption \ref{Assu 7} we can get that $\tau_0\in (0, \frac{3}{4}b_0)$, and it is also easy to check that $\hat{f}$ is of class $C^1$ on $[0, \infty)$

Similarly, we define 
\begin{equation}\label{h3}
\hat{h}(r)\doteqdot \left\{
\begin{array}{rl}
&\tilde{h}(r) \quad \quad \quad \quad  \quad r\in [0, b_0] \\
&\bar{h}(r)+ \sum_{l= 0}^k \delta_l \quad \quad   r\in (b_k, b_{k+ 1}]\ , k= 0, 1, 2, \cdots \\
\end{array} \right.
\end{equation}
where $\displaystyle \delta_0\doteqdot \tilde{h}(b_0)- \lim_{r\rightarrow b_0^{+}} \bar{h}(r)$, when $l\geq 1$, $\delta_l$ is defined in (\ref{delta1}) and (\ref{delta2}).
It is easy to get that $\delta_0\in (0, \frac{3}{4}b_0]$. And it is also easy to check that $\hat{h}$ is of class $C^1$ on $[0, \infty)$.

If we assume 
\begin{assu}\label{Assu 14}
$\displaystyle b_{2k+1}^{1- \frac{1}{2}(\eta_2+ \eta_1)}\leq b_{2k+ 2}^{1- \eta_2}\ , \quad \forall  k\geq 0$.
\end{assu}
We have the following claim about $\tau_i$:
\begin{claim}\label{claim tau}
{\begin{align}
|\tau_0|&\leq \Big(\frac{2b_0}{b_1^{1-\eta_1- \omega_0}}\Big) \min_{r\in (b_j, b_{j+ 1}]}  \bar{f}(r)\ , \quad \quad j\geq 1 \label{tau inq 1.1} \\
|\tau_i|&\leq \Big(\frac{\eta_2- \eta_1}{1- \eta_2}\Big) \min_{r\in (b_i, b_{i+ 1}]} \bar{f}(r) \ , \quad \quad i\geq 1\label{tau inq 1.2} \\
|\tau_i|&\leq \big(b_i^{-\frac{1}{2}(\eta_2- \eta_1)}\big) \min_{r\in (b_j, b_{j+ 1}]} \bar{f}(r) \ , \quad \quad j> i\geq 1 \label{tau inq 1.3}
\end{align}
}
\end{claim}

\pf
{(\ref{tau inq 1.1}) follows directly from the definition of $\tau_0$, (\ref{f1}), (\ref{8.1.3}) and (\ref{8.1.4}).

(\ref{tau inq 1.2}) follows from (\ref{tau1}), (\ref{tau2}), (\ref{f1}) and (\ref{8.1.4}).

There are five cases for (\ref{tau inq 1.3}), in the rest of the proof, $k\geq 0$.

$(1)$. When $i= 2k+ 1$, $j= 2k+ 2$, $k\geq 0$, we have
\begin{align}
\frac{|\tau_i|}{\min_{r\in (b_j, b_{j+ 1}]}\bar{f}(r)}= \Big(\frac{\eta_2- \eta_1}{1- \eta_2}\Big)\cdot \Big(\frac{b_{2k+ 1}}{b_{2k+ 2}}\Big)^{1- \eta_2} \nonumber
\end{align}
using Assumption \ref{Assu 14}, (\ref{tau inq 1.3}) is obtained in this case.

$(2)$. When $i= 2k+ 1$, $j= 2\tilde{k}$, $\tilde{k}> k+ 1$, we have
\begin{align}
\frac{|\tau_i|}{\min_{r\in(b_j, b_{j+ 1}]}\bar{f}(r)}&= \frac{\eta_2- \eta_1}{1- \eta_1}\cdot \frac{\beta_{2k+ 1}}{\beta_{2\tilde{k}- 2}}\cdot \frac{b_{2k+ 2}^{\omega_{2k+ 1}} b_{2k+ 1}^{1- \eta_2}}{b_{2\tilde{k}}^{1- \eta_2} b_{2\tilde{k}- 1}^{\eta_2- \eta_1- \omega_{2\tilde{k}- 2}}} \nonumber \\
&\leq b_{2k+ 2}^{\eta_1- \eta_2+ \omega_{2k+ 1}+ \omega_{2\tilde{k}- 1}} \leq b_{2k+ 1}^{-\frac{1}{2}(\eta_2- \eta_1)} \nonumber
\end{align}
Then (\ref{tau inq 1.3}) is obtained in this case.

$(3)$. When $i= 2k+ 1$, $j= 2\tilde{k}+ 1$, $\tilde{k}> k$, we have
\begin{align}
\frac{|\tau_i|}{\min_{r\in(b_j, b_{j+ 1}]}\bar{f}(r)}&= \frac{1- \eta_2}{1- \eta_1}\cdot\frac{\eta_2- \eta_1}{1- \eta_1}\cdot \frac{\beta_{2k+ 1}}{\beta_{2\tilde{k}}}\cdot \frac{b_{2k+ 2}^{\omega_{2k+ 1}} b_{2k+ 1}^{1- \eta_2}}{b_{2\tilde{k}+ 1}^{1- \eta_1- \omega_{2\tilde{k}}}} \nonumber \\
&\leq b_{2k+ 1}^{-\frac{1}{2}(\eta_2- \eta_1)} \nonumber
\end{align}
Hence (\ref{tau inq 1.3}) holds in this case.

$(4)$. When $i= 2k+ 2$, $j= 2\tilde{k}$, $\tilde{k}> k+ 1$, we have
\begin{align}
\frac{|\tau_i|}{\min_{r\in(b_j, b_{j+ 1}]}\bar{f}(r)}&= \frac{\eta_2- \eta_1}{1- \eta_1}\cdot \frac{\beta_{2k+ 1}}{\beta_{2\tilde{k}- 2}}\cdot \frac{b_{2k+ 2}^{1- \eta_2+ \omega_{2k+ 1}}}{b_{2\tilde{k}}^{1- \eta_2} b_{2\tilde{k}- 1}^{\eta_2- \eta_1- \omega_{2\tilde{k}- 2}}} \nonumber \\
&\leq b_{2k+ 2}^{-\frac{1}{2}(\eta_2- \eta_1)} \nonumber
\end{align}
(\ref{tau inq 1.3}) is got here.

$(5)$. When $i= 2k+ 2$, $j= 2\tilde{k}+ 1$, $\tilde{k}> k$, we have
\begin{align}
\frac{|\tau_i|}{\min_{r\in(b_j, b_{j+ 1}]}\bar{f}(r)}&= \frac{1- \eta_2}{1- \eta_1}\cdot\frac{\eta_2- \eta_1}{1- \eta_1}\cdot \frac{\beta_{2k+ 1}}{\beta_{2\tilde{k}}}\cdot \frac{b_{2k+ 2}^{1- \eta_2+ \omega_{2k+ 1}}}{b_{2\tilde{k}+ 1}^{1- \eta_1- \omega_{2\tilde{k}}}} \nonumber \\
&\leq b_{2k+ 2}^{-\frac{1}{2}(\eta_2- \eta_1)} \nonumber
\end{align}

This completes our proof of (\ref{tau inq 1.3}).
}
\qed

Similarly, We have the following claim about $\delta_i$:
\begin{claim}\label{claim delta}
{\begin{align}
|\delta_0|&\leq 3\frac{b_0}{b_1}\min_{r\in (b_j, b_{j+ 1}]}\bar{h}(r)\ , \quad j\geq 1 \label{delta inq 1.1} \\
|\delta_i|&\leq 4\epsilon_{i- 1}\min_{r\in (b_j, b_{j+ 1}]}\bar{h}(r) \ , \quad 1\leq i\leq j \label{delta inq 1.2}  
\end{align}
}
\end{claim}

\pf
{For $i\geq 1$, we can get the following estimate:
\begin{align}
\frac{\delta_0}{\min_{r\in (b_{2i}, b_{2i+ 1}]}\bar{h}(r)}&\leq \frac{3}{4}\cdot \frac{b_0}{\alpha_{2i}b_{2i+ 1}^{-\epsilon_{2i}} b_{2i}^{1+ \epsilon_{2i}}} \nonumber \\
&= \frac{3}{4} \frac{1}{\alpha_{2i- 1}}\cdot \frac{1+ \epsilon_{2i}}{1- \epsilon_{2i- 1}}\cdot \frac{b_0}{b_{2i- 1}}\cdot \Big(\frac{\alpha_{2i- 1}}{\alpha_{2i- 2}}\cdot \frac{1- \epsilon_{2i- 1}}{1+ \epsilon_{2i- 2}}\Big)^{\frac{1}{\epsilon_{2i- 1}}} \nonumber \\
&\leq \frac{3}{2}\frac{1}{\alpha_0}\frac{b_0}{b_1}\leq 3\frac{b_0}{b_1} \nonumber
\end{align}

Similarly, we can get that for $i\geq 0$,
\begin{align}
\frac{\delta_0}{\min_{r\in (b_{2i+ 1}, b_{2i+ 2}]}\bar{h}(r)}&\leq \frac{3}{4}\cdot \frac{b_0}{\alpha_{2i+ 1}b_{2i+ 2}^{\epsilon_{2i+ 1}} b_{2i+ 1}^{1- \epsilon_{2i+ 1}}} = \frac{3}{4} \frac{b_0}{\alpha_{2i}}\cdot \frac{1- \epsilon_{2i+ 1}}{(1+ \epsilon_{2i})b_{2i+ 1}}\nonumber \\
&\leq \frac{3}{2}\frac{1}{\alpha_0}\frac{b_0}{b_1}\leq 3\frac{b_0}{b_1} \nonumber
\end{align}

By the above two inequalities, we obtain (\ref{delta inq 1.1}).

For $k\geq 1$, 
\begin{align}
\min_{r\in (b_{2k}, b_{2k+ 1}]}\bar{h}(r)&\geq \alpha_{2k}b_{2k+ 1}^{-\epsilon_{2k}} b_{2k}^{1+ \epsilon_{2k}}= \alpha_{2k- 1}b_{2k}\Big(\frac{1+ \epsilon_{2k}}{1- \epsilon_{2k- 1}}\Big) \nonumber\\
&= \alpha_{2k- 2} b_{2k}\big(\frac{b_{2k- 1}}{b_{2k}}\big)^{\epsilon_{2k- 1}} \cdot \Big(\frac{1+ \epsilon_{2k- 2}}{1- \epsilon_{2k- 1}}\Big) \cdot \Big(\frac{1+ \epsilon_{2k}}{1- \epsilon_{2k- 1}}\Big) \nonumber \\
&\geq \alpha_{2k- 2}b_{2k- 1} \nonumber
\end{align}

When $1\leq i< k$, 
\begin{align}
|\delta_{2i}|&= \alpha_{2i- 1}b_{2i}\cdot \frac{\epsilon_{2i- 1}+ \epsilon_{2i}}{1+ \epsilon_{2i}}\leq 2\epsilon_{2i- 1}\alpha_{2i- 1} b_{2i} \nonumber \\
&\leq 2\epsilon_{2i- 1}\alpha_{2k- 2} b_{2k- 1} \leq 2\epsilon_{2i- 1}\min_{r\in (b_{2k}, b_{2k+ 1}]}\bar{h}(r) \nonumber
\end{align}
And 
\begin{align}
|\delta_{2k}|&= \alpha_{2k- 1}b_{2k}\cdot \frac{\epsilon_{2k}+ \epsilon_{2k- 1}}{1+ \epsilon_{2k}} \nonumber \\
&\leq \min_{r\in (b_{2k}, b_{2k+ 1}]}\bar{h}(r) \cdot \frac{1- \epsilon_{2k- 1}}{1+ \epsilon_{2k}} \cdot \frac{\epsilon_{2k}+ \epsilon_{2k- 1}}{1+ \epsilon_{2k}} \nonumber \\
&\leq 2\epsilon_{2k- 1} \min_{r\in (b_{2k}, b_{2k+ 1}]}\bar{h}(r) \nonumber 
\end{align}

When $0\leq i< k$, 
\begin{align}
|\delta_{2i+ 1}|= \alpha_{2i}b_{2i+ 1}\cdot \frac{\epsilon_{2i+ 1}+ \epsilon_{2i}}{1- \epsilon_{2i+ 1}} \leq 4\epsilon_{2i} \min_{r\in (b_{2k}, b_{2k+ 1}]}\bar{h}(r) \nonumber 
\end{align}

From all the above, we get 
\begin{align}
|\delta_i|\leq 4\epsilon_{i- 1}\min_{r\in (b_{2k}, b_{2k+ 1}]}\bar{h}(r)\ , \quad 1\leq i\leq 2k \label{8.23}
\end{align}

For $k\geq 0$, 
\begin{align}
\min_{r\in (b_{2k+ 1}, b_{2k+ 2}]}\bar{h}(r)&\geq \alpha_{2k+ 1}b_{2k+ 1}\Big(\frac{b_{2k+ 2}}{b_{2k+ 1}}\Big)^{\epsilon_{2k+ 1}}= \alpha_{2k}b_{2k+ 1}\Big(\frac{1+ \epsilon_{2k}}{1- \epsilon_{2k+ 1}}\Big) \nonumber \\
&\geq \alpha_{2k}b_{2k+ 1} \nonumber
\end{align}

When $0\leq i\leq k$, 
\begin{align}
|\delta_{2i+ 1}|\leq 4\epsilon_{2i}\alpha_{2k}b_{2k+ 1}\leq 4\epsilon_{2i}\min_{r\in (b_{2k+ 1}, b_{2k+ 2}]}\bar{h}(r)\nonumber
\end{align}

When $1\leq i\leq k$, 
\begin{align}
|\delta_{2i}|\leq 2\epsilon_{2i- 1}\alpha_{2k}b_{2k+ 1}\leq 2\epsilon_{2i- 1} \min_{r\in (b_{2k+ 1}, b_{2k+ 2}]}\bar{h}(r) \nonumber
\end{align}

Hence we obtain that
\begin{align}
|\delta_i|\leq 4\epsilon_{i- 1} \min_{r\in (b_{2k+ 1}, b_{2k+ 2}]}\bar{h}(r)\ , \quad 1\leq i\leq 2k+ 1\label{8.24}
\end{align}

From (\ref{8.23}) and (\ref{8.24}), we (\ref{delta inq 1.2}).
}
\qed

We will assume
\begin{assu}\label{Assu 15}
$\displaystyle \sum_{l= 0}^{\infty} \epsilon_l\doteqdot \delta< 1$, $\quad \displaystyle \sum_{l= 1}^{\infty} b_l^{-\frac{1}{2}(\eta_2- \eta_1)}\doteqdot \tau< 1$
\end{assu}
where $\delta$ and $\tau$ are positive constants to be determined later.

We define $\displaystyle\zeta_k= \sum_{l= 0}^{k} \tau_l$, $\displaystyle\xi_k= \sum_{l= 0}^k \delta_l$, then 
\begin{align}
\hat{f}|_{(b_k, b_{k+ 1}]}= \bar{f}+ \zeta_k\ , \quad \hat{h}|_{(b_k, b_{k+ 1}]}= \bar{h}+ \xi_k \nonumber
\end{align}
Note that we have $Rc(\hat{f}, \hat{h})\geq 0$ on $[0, b_0)$ from $(\mathbf{II})$. In the rest part of $(\mathbf{III})$, we will prove $Rc(\hat{f}, \hat{h})> 0$ on $(b_0, \infty)$ except at points $b_j$, $j=1, 2, \cdots$.

$(\mathbf{i})$. We firstly consider the interval $(b_0, b_1]$.

on $(b_0, b_1]$, from Assumption \ref{Assu 7}, we get 
\begin{align}
Rc|_{k_1}(\hat{f}, \hat{h})&\geq \frac{2\big[r^2- \bar{f}^2(1- \eta_1)^2- 2\bar{f}(\bar{f}+ \tau_0) \big]}{(\bar{f}+ \tau_0)^2 r^2} \nonumber \\
&\geq \frac{2\bar{f}^2\Big[\big(\delta_0^{-1}b_1^{\omega_0}r^{\eta_1}\big)^2- (1- \eta_1)^2- 2\big(1+ \frac{3}{4}\delta_0^{-1}b_1^{\omega_0}b_0^{\eta_1}\big)\Big]}{(\bar{f}+ \tau_0)^2 r^2} \nonumber\\
&\geq \frac{2\bar{f}^2\Big[\delta_0^{-1}b_1^{\omega_0}\big(b_0^{2\eta_1}- 1- 4b_0^{\eta_1}\big)\Big]}{(\bar{f}+ \tau_0)^2 r^2}> 0 \nonumber
\end{align}

From Assumption \ref{Assu 10}, we obtain that
\begin{align}
Rc(\vec{n}, \vec{n})&> \frac{3\eta_1(1- \eta_1)\bar{f}}{\big(\bar{f}+ \tau_0\big)r^2}- \frac{4\epsilon_0(\epsilon_0+ 1)}{r^2} \nonumber \\
&\geq \frac{2}{r^2}\Big[\frac{3}{2}\cdot\frac{(1- \eta_1)\eta_1}{1+ b_0^{\eta_1}b_1^{\omega_0}}- 2\epsilon_0(\epsilon_0+ 1) \Big]> 0 \nonumber
\end{align}

We assume that for $i\geq 0$,
\begin{assu}\label{Assu 17}
{$\displaystyle \epsilon_{2i}< \omega_{2i}$
}
\end{assu}
then $\hat{h}'> \hat{f}'$ on $(b_0, b_1]$. Combining with $\hat{h}(b_0)> \hat{f}(b_0)$, we get that $\hat{h}> \hat{f}$ on $(b_0, b_1]$. 

From (\ref{8.1.0}), combining with Assumption \ref{Assu 9}, we can get
\begin{align}
Rc|_{k_2}(\hat{f}, \hat{h})&\geq \frac{6}{\hat{h}^2}- \frac{\hat{h}''}{\hat{h}}- \frac{3(\hat{h}')^2}{\hat{h}^2}- \frac{3\hat{f}'\hat{h}'}{\hat{f}\hat{h}} \nonumber \\
&= \frac{1}{\hat{h}^2}\Big\{6- \big[\frac{\bar{h}}{r}+ \frac{\delta_0}{r}\big](1+ \epsilon_0)\epsilon_0\frac{\bar{h}}{r}- 3\big[(1+ \epsilon_0)\frac{\bar{h}}{r}\big]^2 \nonumber \\
&\quad \quad - 3(1- \eta_1)\frac{\bar{f}}{\bar{f}+ \tau_0}(1+ \epsilon_0)\frac{\bar{h}}{r}\big(\frac{\bar{h}}{r}+ \frac{\delta_0}{r} \big) \Big\} \nonumber \\
&\geq \frac{1}{\hat{h}^2} \Big[6- \alpha_0(1+ \epsilon_0)\epsilon_0(\frac{3}{4}+ \alpha_0)- 3\big(\alpha_0(1+ \epsilon_0)\big)^2 \nonumber \\
&\quad \quad - 3(1- \eta_1)\alpha_0(1+ \epsilon_0)(\alpha_0+ \frac{3}{4}) \Big] \nonumber \\
&\geq \frac{1}{\hat{h}^2} [6- 2\epsilon_0- 3- 6(1-\eta_1)]> 0 \nonumber
\end{align}

So, we proved that $Rc(\hat{f}, \hat{h})> 0$ on $(b_0, b_1)$.

$(\mathbf{ii})$. Next we consider the interval $(b_{2i}, b_{2i+ 1}]$, $i\geq 1$.

We assume that
\begin{assu}\label{Assu 18}
{$\displaystyle \frac{2b_0}{b_1^{1- \eta_1- \omega_0}}+ \frac{\eta_2- \eta_1}{1- \eta_2}+ \tau< \eta_1^{3}\ , \quad b_1^{2\eta_1}> 2+ 20\Big(1- 3\frac{b_0}{b_1}- 4\delta \Big)^{-1}$
}
\end{assu}
then combining with Claim \ref{claim tau} and Claim \ref{claim delta}, we get
\begin{align}
&Rc|_{k_1}(\hat{f}, \hat{h})\geq \frac{2\Big[1- (1- \eta_1)^2\big(\frac{\bar{f}}{r}\big)^2\Big]}{\big(\bar{f}+ \zeta_{2i}\big)^2}- \frac{4(1- \eta_1)(1+ \epsilon_{2i})\big(\frac{\bar{f}}{r}\big)\big(\frac{\bar{h}}{r}\big)}{\big(\bar{f}+ \zeta_{2i}\big) \big(\bar{h}+ \xi_{2i}\big)} \label{8.25}\\
&= \frac{2\bar{f}^2}{\big(\bar{f}+ \zeta_{2i}\big)^2 r^2} \Big[\frac{r^2}{\bar{f}^2}- (1- \eta_1)^2- 2(1- \eta_1)(1+ \epsilon_{2i})\frac{\bar{h}}{\bar{h}+ \xi_{2i}}\big(1+ \frac{\zeta_{2i}}{\bar{f}}\big)\Big] \nonumber \\
&\geq \frac{2\bar{f}^2}{\big(\bar{f}+ \zeta_{2i}\big)^2 r^2} \Big[b_2^{2\eta_1}- 1- \frac{2\Big(1+ \frac{2b_0}{b_1^{1- \eta_1- \omega_0}}+ \frac{\eta_2- \eta_1}{1- \eta_2}+ \tau \Big)}{1- \frac{3b_0}{b_1}- 4\delta}\Big] \nonumber \\
&> 0 \nonumber
\end{align}

If we further assume that
\begin{assu}\label{Assu 19}
{$\displaystyle \epsilon_0\leq \frac{1}{10}\eta_1(1- \eta_1)\big(1- 3\frac{b_0}{b_1}- 4\delta \big)$
}
\end{assu}
then we get
\begin{align}
Rc(\vec{n}, \vec{n})&= \frac{3\eta_1(1- \eta_1)\bar{f}}{(\bar{f}+ \zeta_{2i})r^2}- \frac{4\epsilon_{2i}(1+ \epsilon_{2i})\bar{h}}{(\bar{h}+ \xi_{2i})r^2} \nonumber \\
&\geq \frac{3\eta_1(1- \eta_1)}{\Big(1+ \frac{2b_0}{b_1^{1- \eta_1- \omega_0}}+ \frac{\eta_2- \eta_1}{1- \eta_2}+ \tau\Big)r^2}- \frac{4\epsilon_{2i}(1+ \epsilon_{2i})}{\Big(1- \frac{3b_0}{b_1}- 4\delta \Big)r^2} \nonumber \\
&\geq \frac{3}{r^2\Big(1- 3\frac{b_0}{b_1}- 4\delta \Big)} \Big[\frac{\eta_1(1- \eta_1)\Big(1- 3\frac{b_0}{b_1}- 4\delta \Big)}{\Big(1+ \frac{2b_0}{b_1^{1- \eta_1- \omega_0}}+ \frac{\eta_2- \eta_1}{1- \eta_2}+ \tau\Big)}- 2\epsilon_0\Big] \nonumber \\
&> 0 \nonumber 
\end{align}

We assume that 
\begin{assu}\label{Assu 20}
{$b_{1}^{\eta_1}\geq \frac{100}{1- 4\delta}$
}
\end{assu}

From Assumptions \ref{Assu 17} and \ref{Assu 20}, we get $\bar{h}\geq \alpha_0 b_{2i}^{\eta_1} \bar{f}$ on $(b_{2i}, b_{2i+ 1}]$. Also note that the following holds:
\begin{align}
\hat{f}= \bar{f}+ \zeta_{2i}\leq 5\bar{f}\ , \quad \hat{h}= \bar{h}+ \xi_{2i}\geq (1- 4\delta)\bar{h} \nonumber 
\end{align}
the above three inequalities imply that $\hat{h}\geq \hat{f}$ on $(b_{2i}, b_{2i+ 1}]$.

We further assume that
\begin{assu}\label{Assu 21}
{$\displaystyle \frac{3b_0}{b_1}+ 4\delta\leq \eta_1$
}
\end{assu}
Observe that $\epsilon_{0}< \frac{\eta_1^2}{1+ \eta_1}$, then using Assumption \ref{Assu 9}, we have
\begin{align}
Rc|_{k_2}(\hat{f}, \hat{h})&\geq \frac{1}{\hat{h}^2}\Big[6- \hat{h}\hat{h}''- 3(\hat{h}')^2- 3\frac{\hat{f}'}{\hat{f}}\hat{h}'\hat{h} \Big] \nonumber \\
&= \frac{1}{\hat{h}^2}\Big\{6- \frac{\bar{h}^2}{r^2}\epsilon_{2i}(1+ \epsilon_{2i})\big(\frac{\hat{h}}{\bar{h}}\big)- 3\big[(1+ \epsilon_{2i})\frac{\bar{h}}{r}\big]^2 \nonumber \\
&\quad\quad - 3\frac{\bar{f}(1- \eta_1)(1+ \epsilon_{2i})}{\bar{f}+ \zeta_{2i}}\cdot \frac{\bar{h}(\bar{h}+ \xi_{2i})}{r^2}\Big\} \nonumber \\
&\geq \frac{1}{\hat{h}^2}\Big\{6- \epsilon_{2i}(1+ \epsilon_{2i})\alpha_{2i}^2\big(1+ 3\frac{b_0}{b_1}+ 4\delta\big)- 3\big[\alpha_{2i}(1+ \epsilon_{2i})\big]^2 \nonumber \\
&\quad \quad 3\frac{(1- \eta_1)(1+ \epsilon_{2i})\big(1+ 3\frac{b_0}{b_1}+ 4\delta\big)}{1- \frac{2b_0}{b_1^{1- \eta_1- \omega_0}}- \frac{\eta_2- \eta_1}{1- \eta_2}- \tau} \alpha_{2i}^2 \Big\} \nonumber \\
&\geq \frac{1}{\hat{h}^2} \Big[3- \frac{\eta_1^2}{1+ \eta_1}(1+ \eta_1)- 3\frac{1- \eta_1}{1- \eta_1^3}(1+ \eta_1) \Big] > 0 \nonumber
\end{align}

From all the above, we proved that $Rc(\hat{f}, \hat{h})> 0$ on $(b_{2i}, b_{2i+ 1})$, where $i\geq 1$.

$(\mathbf{iii})$. Finally we consider the interval $(b_{2i+ 1}, b_{2i+ 2})$, where $i\geq 0$.

From Assumption \ref{Assu 18}, similarly as (\ref{8.25}), we get
\begin{align}
Rc|_{k_1}&\geq \frac{2\bar{f}^2}{\big(\bar{f}+ \zeta_{2i+ 1}\big)^2 r^2}\Big[\Big(\frac{r^{\eta_2}}{\beta_{2i+ 1}b_{2i+ 2}^{\omega_{2i+ 1}}}\Big)^{2}- 1 \nonumber \\
&\quad \quad \quad - \frac{2}{1- \frac{3b_0}{b_1}- 4\delta}\Big(1+ \frac{2b_0}{b_1^{1- \eta_1- \omega_0}}+ \frac{\eta_2- \eta_1}{1- \eta_2}+ \tau \Big)\Big] \nonumber \\
&\geq \frac{2\bar{f}^2}{\big(\bar{f}+ \zeta_{2i+ 1}\big)^2 r^2}\Big[\Big(\frac{1- \eta_2}{\beta_{2i}(1- \eta_1)}\Big)^{2}b_{2i+ 1}^{2\eta_1}- 1- \frac{4}{1- \frac{3b_0}{b_1}- 4\delta}\Big] \nonumber \\
&\geq \frac{\bar{f}^2}{\big(\bar{f}+ \zeta_{2i+ 1}\big)^2 r^2}\Big[b_{2i+ 1}^{2\eta_1}- 2-  8\Big(1- \frac{3b_0}{b_1}- 4\delta\Big)^{-1}\Big] > 0 \nonumber
\end{align}

And $Rc(\vec{n}, \vec{n})> 0$ is trivial by $\hat{f}''< 0$ and $\hat{h}''< 0$.

It is easy to see that we also have
\begin{align}
\hat{f}= \bar{f}+ \zeta_{2i+ 1}\leq 5\bar{f}\ , \quad \hat{h}= \bar{h}+ \xi_{2i+ 1}\geq (1- 4\delta)\bar{h} \label{8.26}
\end{align}

Using (\ref{8.1.3}) and (\ref{8.8}), we have
\begin{align}
\bar{h}\geq \frac{\alpha_{2i}}{\beta_{2i}}\cdot \frac{1+ \epsilon_{2i}}{1- \epsilon_{2i+ 1}}\cdot \frac{1- \eta_2}{1- \eta_1}b_{2i+ 1}^{\eta_1+ \omega_{2i}} \bar{f}\geq \frac{\alpha_0}{2}b_{2i+ 1}^{\eta_1}\bar{f} \label{8.27}
\end{align}
from (\ref{8.26}), (\ref{8.27}) and Assumption \ref{Assu 20}, we can get $\hat{h}> \hat{f}$ on $(b_{2i+ 1}, b_{2i+ 2})$, where $i\geq 0$.

From Assumption \ref{Assu 18} and Assumption \ref{Assu 21}, we can get
\begin{align}
Rc|_{k_2}(\hat{f}, \hat{h})&\geq \frac{1}{\hat{h}^2}\Big[6- 3\Big(\frac{\bar{h}}{r}\Big)^2(1- \epsilon_{2i+ 1})^2 \nonumber \\
&\quad \quad \quad - 3\frac{\bar{h}^2}{r^2}(1- \epsilon_{2i+ 1})(1- \eta_2)\cdot \frac{\Big(1+ \frac{\xi_{2i+ 1}}{\bar{h}}\Big)}{\Big(1+ \frac{\zeta_{2i+ 1}}{\bar{f}}\Big)}\Big] \nonumber \\
&\geq \frac{1}{\hat{h}^2}\Big\{6- 3\big[\alpha_{2i}(1+ \epsilon_{2i})\big]^2\Big[1+ \frac{(1+ \eta_1)(1- \eta_2)}{(1- \eta_1^3)(1- \epsilon_{2i+ 1})}\Big]\Big\} \nonumber \\
&\geq \frac{1}{\hat{h}^2}\Big\{6- 3\Big[1+ \frac{1+ \eta_1}{1+ \eta_1+ \eta_1^2}\Big]\Big\}> 0 \nonumber
\end{align}
in the last inequality, we used the inequality $\displaystyle\frac{1- \eta_2}{1- \eta_1}< 1- \epsilon_0$ from (\ref{8.1.0}).

From all the above, we get $Rc(\hat{f}, \hat{h})> 0$ on $(b_0, \infty)$ excepts at points $b_i$, $i= 1, 2, \cdots$.

\subsection{Step $(\mathbf{IV})$}\label{subsection 8.4}
$\newline$ Finally, we can remove the jump discontinuities in the functions, $\hat{f}''$, $\hat{h}''$, by modifying them by linear interpolation, in arbitrarily small neighborhoods of the points, $\{b_i\}_{i= 0}^{\infty}$. Call the resulting functions $f''$, $h''$, and let the corresponding functions, $f$, $h$, be obtained by integration with respect to $r$, subject to the conditions, $f(0)= h(0)= 0$, $f'(0)= h'(0)= 1$. The modification in the second derivatives can be performed on intervals whose size decreases rapidly enough to ensure the nonnegative property of $Rc|_{k_1}(f, h)$, $Rc|_{k_2}(f, h)$ and $Rc(\vec{n}, \vec{n})(f, h)$ on $[0, \infty)$.


For $(M^8, g)$, $M^8$ is diffeomorphic to $\mathbb{R}^8$, $g= dr^2+ f^2k_1+ h^2 k_2$, define two sequences $\{t_i\}_{i= 0}^{\infty}$, $\{\tilde{t}_i\}_{i= 0}^{\infty}$ as the following:
\begin{align}
t_i= \Big(b_{2i+ 1}^{1- \epsilon_{2i}}\Big)^2\ , \quad \quad \tilde{t}_i= \Big(b_{2i+ 2}^{1- \epsilon_{2i+ 1}}\Big)^2\ , \quad \quad i= 0, 1, 2, \cdots \nonumber
\end{align}
And define the scaling metrics $g_i\doteqdot t_i^{-1}g$ and $\tilde{g}_i\doteqdot \tilde{t}_i^{-1}g$, we also assume that
\begin{assu}\label{Assu 22}
{$\displaystyle\lim_{i\rightarrow \infty}b_i^{\epsilon_i}= 1\ , \quad \lim_{i\rightarrow \infty}b_{i+ 1}^{\epsilon_i^2}= 1\ , \quad \lim_{i\rightarrow \infty}b_{i+ 1}^{\epsilon_i}= \infty$
}
\end{assu}

It is not hard to check that we can find sequences $\{b_i\}$, $\{\alpha_i\}$, $\{\beta_i\}$, $\{\omega_i\}$, $\{\epsilon_i\}$ and $\eta_1$, $\eta_2$ satisfying the Assumptions $1-20$. Hence we get  
\begin{align}
\big(M^8, g_i, y, \nu_i\big)\stackrel{d_{GH}}{\longrightarrow}\big(M_{\infty}, \rho_{\infty}, y_{\infty}, \nu_{\infty}\big) \label{conv 1}
\end{align}
and define $\nu_i(A)\doteqdot t_i^{\frac{n}{2}} V(\sqrt{t_i})^{-1} \mu_{i}(A)$, where $\mu_i$ is the volume element determined by metric $g_i$. 

$M_{\infty}$ is diffeomorphic to $\mathbb{R}^5$ with metric $\rho_{\infty}= dr^2+ \frac{1}{4}r^2g^{\mathbb{S}^4}$, and 
\begin{align}
\nu_{\infty}\big(B_{\infty}(r)\big)= r^{8- 3\eta_1} \nonumber
\end{align}

On the other side, we have 
\begin{align}
\big(M^8, \tilde{g}_i, y, \tilde{\nu}_i\big)\stackrel{d_{GH}}{\longrightarrow}\big(\tilde{M}_{\infty}, \tilde{\rho}_{\infty}, y_{\infty}, \tilde{\nu}_{\infty}\big) \label{conv 2}
\end{align}
where $\tilde{\nu}_i(A)\doteqdot \tilde{t}_i^{\frac{n}{2}} V(\sqrt{\tilde{t}_i})^{-1}\tilde{\mu}_i(A)$, and $\tilde{u}_i$ is the volume element determined by $\tilde{g}_i$. 

$\tilde{M}_{\infty}$ is diffeomorphic to $\overline{\mathbb{R}^{+}}$ with metric $\tilde{\rho}_{\infty}= dr^2$ and 
\begin{align}
\tilde{\nu}_{\infty}\big(B_{\infty}(r)\big)= r^{8- 3\eta_2} \nonumber
\end{align} 

From the proof of Theorem \ref{thm 11.0}, we can get that for rotational symmetric functions on $(M_{\infty}, \rho_{\infty}, y_{\infty}, \nu_{\infty})$ and $(\tilde{M}_{\infty}, \tilde{\rho}_{\infty}, y_{\infty}, \tilde{\nu}_{\infty})$ respectively, 
\begin{align}
\Delta_{(\rho_{\infty}, \nu_{\infty})}= \frac{\partial^2}{\partial r^2}+ \frac{7- 3\eta_1}{r}\cdot \frac{\partial}{\partial r}\ , \quad \quad \Delta_{(\tilde{\rho}_{\infty}, \tilde{\nu}_{\infty})}= \frac{\partial^2}{\partial r^2}+ \frac{7- 3\eta_2}{r}\cdot \frac{\partial}{\partial r} \nonumber
\end{align}

Then it is not hard to get 
\begin{align}
H_{\infty}(x_{\infty}, y_{\infty}, t)&= C_{H}\cdot t^{-\frac{1}{2}(8- 3\eta_1)} \exp{\Big(-\frac{d_{\rho_{\infty}}(x_{\infty}, y_{\infty})^2}{4t} \Big)} \label{HK 1} \\
\tilde{H}_{\infty}(x_{\infty}, y_{\infty}, t)&= C_{\tilde{H}}\cdot t^{-\frac{1}{2}(8- 3\eta_2)} \exp{\Big(-\frac{d_{\tilde{\rho}_{\infty}}(x_{\infty}, y_{\infty})^2}{4t} \Big)} \label{HK 2}
\end{align}
where $\displaystyle C_H= \Big(\int_0^{\infty} e^{-\frac{u^2}{4}} u^{7- 3\eta_1} du\Big)^{-1}$ and $\displaystyle C_{\tilde{H}}= \Big(\int_0^{\infty} e^{-\frac{u^2}{4}} u^{7- 3\eta_2} du\Big)^{-1}$, which follows from 
\begin{align}
\int_{M_{\infty}} H_{\infty} d\nu_{\infty}= 1\ , \quad \quad \int_{\tilde{M}_{\infty}} \tilde{H}_{\infty} d\tilde{\nu}_{\infty}= 1 \nonumber 
\end{align}

From (\ref{3.1}) and Proposition \ref{prop 3.14}, we get
\begin{align}
\lim_{i\rightarrow \infty} V(\sqrt{t_i}) H(x, y, t_i)&= H_{\infty}(y_{\infty}, y_{\infty}, 1)= C_H \nonumber \\
\lim_{i\rightarrow \infty} V(\sqrt{\tilde{t}_i}) H(x, y, \tilde{t}_i)&= \tilde{H}_{\infty}(y_{\infty}, y_{\infty}, 1)= C_{\tilde{H}} \nonumber
\end{align}

But from $\eta_1< \eta_2$, it is easy to see that $C_H< C_{\tilde{H}}$. Hence
\begin{align}
\lim_{i\rightarrow \infty} V(\sqrt{t_i}) H(x, y, t_i)< \lim_{i\rightarrow \infty} V(\sqrt{\tilde{t}_i}) H(x, y, \tilde{t}_i) \nonumber 
\end{align}

This answers one open question raised in \cite{Li} negatively. That is, without maximal volume growth assumption, $\displaystyle \lim_{t\rightarrow \infty} V(\sqrt{t}) H(x, y, t)$ does not generally exist.

\appendix
\section{Rellich-type Compactness theorem}\label{App 1}
Similar with the Rellich-Kondrakov Theorem for Sobolev spaces on a fixed domain, we have Rellich-type Compactness Theorem in the Gromov-Hausdorff sense, which was used in the proof of Theorem \ref{thm 4.2}. In this appendix we will give a complete proof of Rellich-type Compactness Theorem. 

We firstly state some background knowledge needed for the proof.
\begin{definition}[Measure approximation, \cite{KS0}]\label{def ma}
{Let $M_i$ and $M_{\infty}$ be measure spaces.  A net $\{\varphi_i: M_i\supset \mathcal{D}(\varphi_i)\rightarrow M_{\infty}\}$ of maps is called a \textbf{measure approximation} if the following are satisfied:
\begin{itemize}
\item Each $\varphi_i$ is a measurable map from a Borel subset $\mathcal{D}(\varphi_i)$ of $M_i$ to $M_{\infty}$.
\item The push-forward by $\varphi_i$ of the measure on $M_i$ weakly-* converges to the measure on $M_{\infty}$, i.e., for any $f\in C_c(M_{\infty})$, 
\begin{align}
\lim_{i\rightarrow \infty} \int_{\mathcal{D}(\varphi_i)} f\circ \varphi_i d\nu_{i}= \int_{M_{\infty}} f d\nu_{\infty} \label{ma 1}
\end{align}
where $C_c(M_{\infty})$ is the set of continuous functions on $M_{\infty}$ with compact support.
\end{itemize}
}
\end{definition}

As in \cite{Fukaya} (also see \cite{KS0}), there is another definition of measured Gromov-Hausdorff convergence as the following.

\begin{definition}[Measured Gromov-Hausdorff convergence]\label{def MGH1}
{If $\nu_i$, $\nu_{\infty}$ are Borel regular measures on $M_i^n$, $M_{\infty}$, we say that $(M_i^n, y_i, \rho_i, \nu_i)$ converges to $(M_{\infty}, y_{\infty}, \rho_{\infty}, \nu_{\infty})$ \textbf{in the measured Gromov-Hausdorff sense}, if 
there exists a measure approximation $\{\varphi_i: M_i\rightarrow M_{\infty}\}$,  such that each $\varphi_i$ is an $\epsilon_i$-Gromov-Hausdorff approximation for some $\epsilon_i\rightarrow 0$, and $\varphi_i(y_i)= y_{\infty}$.
}
\end{definition}

\begin{remark}[Fukaya's definition VS definition of Cheeger \& Colding]\label{rem Fu vs CC}
{$\newline$ If $(M_i^n, y_i, \rho_i, \nu_i)$ converges to $(M_{\infty}, y_{\infty}, \rho_{\infty}, \nu_{\infty})$ in the measured Gromov-Hausdorff sense, from the above definition, we have
\begin{itemize}
\item $(M_i^n, y_i, \rho_i)\stackrel{d_{GH}}{\longrightarrow} (M_{\infty}, y_{\infty}, \rho_{\infty})$. 
\item In addition, for any $x_i\rightarrow x_{\infty}$, ($x_i\in M_i^n$, $x_{\infty}\in M_{\infty}$), $r> 0$, we have 
\[\nu_{i}\Big(B_{i}(x_i, r)\Big)\rightarrow \nu_{\infty}\Big(B_{\infty} (x_{\infty}, r)\Big)\]
where $(M_{\infty}, \rho_{\infty})$ is a length space with length metric $\rho_{\infty}$, and  
\begin{align}
B_i(x_{i}, r)= \{z\in M_i^n|\  d_{\rho_i}(z, x_i)\leq r\}\ , \quad B_{\infty}(x_{\infty}, r) = \{z\in M_{\infty}|\  d_{\rho_{\infty}}(z, x_{\infty})\leq r\} \nonumber
\end{align}  
\end{itemize}

The above two items were used to define the measured Gromov-Hausdorff convergence in \cite{Cheeger}(also see Definition \ref{def MGH}). Hence the definition of the measured Gromov-Hausdorff convergence we chose (following \cite{Fukaya}), implies the measured Gromov-Hausdorff convergence discussed in Cheeger and Colding's work.

However, from Proposition $2.2$ in \cite{KS0}, in fact, the definition of the measured Gromov-Hausdorff convergence in the Definition \ref{def MGH} is equivalent to the one used by Cheeger and Colding.
}
\end{remark}

In most parts of the paper, we used the definition of the measured Gromov-Hausdorff convergence by Cheeger and Colding as in Definition \ref{def MGH}. However, to prove the following Rellich-type compactness result in the Gromov-Hausdorff sense, we will use the definition of Fukaya in the Definition \ref{def MGH1}. 

And as in \cite{KS}, we define $L^p$ convergence in Gromov-Hausdorff topology in the following.
\begin{definition}[$L^p$ Convergence in G-H topology]\label{def Lp}
{Assume that $\{f_i\}_{i= 1}^{\infty}$ are functions on $M_i^n$, $f_{\infty}$ is a function on $M_{\infty}$, we say $f_i\rightarrow f_{\infty}$ in $L^p$ sense on $U\subset M_{\infty}$, if there exists $f_{\infty}^{(j)}\in C_c(U)$, such that 
\begin{align}
\lim_{j\rightarrow \infty} \int_U |f_{\infty}^{(j)}- f_{\infty}|^p d\nu_{\infty}= 0 \ , \quad \lim_{j\rightarrow \infty} \varlimsup_{i\rightarrow \infty} \int_{U_i} |f_i- f_{\infty}^{(j)}\circ\varphi_i|^2 d\nu_i= 0 \label{Lp 1}
\end{align}
where $\varphi_i: U_i\rightarrow U$ is a measure approximation and an $\epsilon_i$-Gromov-Hausdorff approximation for some $\epsilon_i\rightarrow 0$.
}
\end{definition}

\begin{theorem}[Rellich-type Compactness Theorem]\label{thm 2.3}
{Assume 
\begin{align}
B_i(x_i, r)\subset (M_i^n, y_i, \rho_i, \nu_i)\ , \quad \quad B_{\infty}(x_{\infty}, r)\subset (M_{\infty}, y_{\infty}, \rho_{\infty}, \nu_{\infty}) \nonumber
\end{align} 
and $\displaystyle B_i(x_i, r)\stackrel{d_{GH}}{\longrightarrow} B_{\infty}(x_{\infty}, r) $ in the measured Gromov-Hausdorff sense, $u_i$ is a function on $M_i^n$, and for some fixed constant $N> 0$, 
\begin{align}
\int_{B_i(x_i, r)} \Big[|u_i|^2+ |\nabla u_i|^2 \Big]d\nu_i\leq N \label{1.2.1}
\end{align}
Then there exists a subsequence of $\{u_i\}$ such that $u_i\rightarrow u_{\infty}$ in $L^2$ sense on any $K_{\infty}\subset\subset \mathring{B}_{\infty}(x_{\infty}, r)$, where $\mathring{B}_{\infty}(x_{\infty}, r)$ denotes the interior of  $B_{\infty}(x_{\infty}, r)$.
}
\end{theorem}

\begin{remark}\label{rem apd1}
{The proof of the above theorem was sketched in \cite{Ding}. Following closely the argument in \cite{KS} (see Theorem $4.15$ there), also compare \cite{CMAnn}, we give a detailed proof here. 
}
\end{remark}

\pf
{For $K_{\infty}\subset\subset \mathring{B}_{\infty}(x_{\infty}, r)$, assume $d_{\infty}(K_{\infty}, \partial B_{\infty})= 100 r_0> 0$. Then there exists $i_0> 0$, for $i> i_0$, $d_{\rho_i}\big(\phi_i(K_{\infty}), \partial B_i\big)= 10r_0> 0$. 

Define $K_i= \phi_i(K_{\infty})\subset B_i(x_i, r)$. Take a sequence of numbers $r_j \searrow 0$, $j= 1, 2, \cdots$, and $r_j\leq 0$. Let $\{B_i(z_{jk}^i, r_j)\}_{k= 1}^{N_j^i}$ be a maximal set of disjoint balls with radius $r_j$, centers $z_{jk}^i$ in $K_i$.

First, by the volume comparison theorem, 
\begin{align}
\nu_i\big(B_i(z_{jk}^i, r_j)\big)&\geq \nu_i\big(B_i(z_{jk}^i, r_j+ 2r)\big)\cdot \Big(\frac{r_j}{r_j+ 2r}\Big)^{n} \nonumber
&\geq C(r_j, r, n)\nu_i\big(B_i(x_i, r)\big)  \nonumber
\end{align}

Note 
\begin{align}
\sum_{k= 1}^{N_j^i} \nu_i\big(B_i(z_{jk}^i, r_j)\big)\leq \nu_i\big(B_i(x_i, r)\big) \nonumber
\end{align}
therefore
\begin{align}
N_j^i\leq C(r_j, r, n) \nonumber
\end{align}

It follows from maximality that double the balls covers $K_i$. We now get $N_j^i$ disjoint subsets $S_{j1}^i, S_{j2}^i \cdots, S_{j N_j^i}^{i}$ which cover $K_i$, where 
\begin{align}
S_{jk}^i= B_i(z_{jk}^i, 2r_j)\backslash \Big(\cup_{l= 1}^{k- 1} B_i(z_{jl}^i, 2r_j)\Big) \nonumber
\end{align}

We define a step function $\bar{u}_j^i: K_i\rightarrow \mathbb{R}$ by $\bar{u}_j^i= \bar{u}_{jk}^i$ on each $S_{jk}^i$, where 
\begin{align}
\bar{u}_{jk}^i= \frac{1}{\nu_i\Big(B_{i}\big(z_{jk}^i, 2r_j \big)\Big)}\int_{B_{i}\big(z_{jk}^i, 2r_j \big)} u_i d\nu_i \nonumber
\end{align}

Let $\eta(y)$ be the number of $k$, such that $y\in B_{i}\big(z_{jk}^i, 4r_j \big)$ and let $\bar{C}_i= \max_{y\in B_i(x_i, r)} \eta(y)$. 

If $y\in \cap_{m= 1}^{\eta(y)} B_{i}\big(z_{jk}^i, 4r_j\big)$, it follows that $B_i(y, 5r_j)$ contains all of the balls 
\[B_{i}\big(z_{j1}^i, r_j \big), B_{i}\big(z_{j2}^i, r_j \big), \cdots, B_{i}\big(z_{j\eta(y)}^i, r_j\big)\] 
Since these are disjoint, 
\begin{align}
\sum_{m= 1}^{\eta(y)} \nu_i\Big(B_{i}\big(z_{jm}^i, r_j \big)\Big)\leq \nu_i\Big(B_i(y, 5r_j)\Big) \label{2.3.1}
\end{align}

Also for each $m= 1, 2, \cdots, \eta(y)$, the doubling condition together with the triangle inequality yields
\begin{align}
\nu_i\Big(B_i(y, 5r_j)\Big)\leq \nu_i\Big(B_{i}\big(z_{jm}^i, 9r_j \big)\Big)\leq 9^n \nu_i\Big(B_{i}\big(z_{jm}^i, r_j\big)\Big) \label{2.3.2}
\end{align}

Combining (\ref{2.3.1}) and (\ref{2.3.2}), we see that $\eta(y)\leq 9^n= C(n)$, hence $\bar{C}_i\leq C(n)$.

We have the following claim.
\begin{claim}\label{claim 2.3.1}
{\begin{align}
\lim_{j\rightarrow \infty}\varlimsup_{i\rightarrow \infty} \int_{K_i} |u_i- \bar{u}_j^i|^2 d\nu_i= 0 \label{2.3.3}
\end{align}
}
\end{claim}

\pf
{\begin{align}
\int_{K_i} |u_i- \bar{u}_j^i|^2&= \sum_{k= 1}^{N_j^i} \int_{S_{jk}^i} |u_i- \bar{u}_{jk}^i|^2 \nonumber \\
&\leq \sum_{k= 1}^{N_j^i} \int_{B_i(z_{jk}^i, 2r_j)} |u_i- \bar{u}_{jk}^i|^2 \nonumber \\
&\leq \sum_{k= 1}^{N_j^i} C(n) (2r_j)^2 \int_{B_i(z_{jk}^i, 4r_j)} |\nabla u_i|^2 \nonumber \\
&\leq \bar{C}_i C(n) r_j^2 \int_{B_i(x_i, r)} |\nabla u_i|^2\leq C(n, N)r_j^2 \nonumber
\end{align}

The conclusion follows from it, and $\lim_{j\rightarrow \infty}r_j= 0$.
}
\qed

It follows from the Cauchy-Schwarz inequality, together with the doubling condition that
\begin{align}
\bar{u}_{jk}^i&\leq \frac{1}{\sqrt{\nu_i\Big(B_{i}\big(z_{jk}^{i}, 2r_j\big)\Big)}} \Big(\int_{B_{i}\big(z_{jk}^{i}, 2r_j\big)} u_i^2\Big)^{\frac{1}{2}} \nonumber \\
&\leq \frac{N}{\sqrt{\nu_i\Big(B_{i}\big(x_i, r\big)\Big)}}\cdot \Big(\frac{r_j+ r}{r_j}\Big)^{\frac{n}{2}} \nonumber \\
&\leq N\cdot C\Big(n, r_0, r, \nu_{\infty}\big(B_{\infty}(x, r)\big) \Big) \label{2.3.4}
\end{align}
note that the bound on the right side is independent of $i, k$. Hence for fixed $j, k$, $\{\bar{u}_{jk}^i\}_{i= 1}^{\infty}$ has a convergent subsequence.

There is a measure approximation $\varphi_i: B_i(x_i, r)\rightarrow B_{\infty}(x_{\infty}, r)$, such that each $\varphi_i$ is an $\epsilon_i$-approximation for some $\epsilon_i\searrow 0+$. There is a subsequence of $\{i\}$ depending on $j$, denoted as $\mathcal{I}_j$, such that for every $k= 1, 2, \cdots, N_j^i$,
\begin{align}
z_{jk}\doteqdot \lim_{i\rightarrow \infty}\varphi_i (z_{jk}^i)\ , \quad N_j\doteqdot \lim_{i\rightarrow \infty}N_J^i\ , \quad \bar{u}_{jk}\doteqdot \lim_{i\rightarrow \infty} \bar{u}_{jk}^i \nonumber
\end{align}
all the above limits exist, where $i\in \mathcal{I}_j$.

By (\ref{2.3.4}), replacing $\mathcal{I}_j$ with a subset of $\mathcal{I}_j$, also denoted as $\mathcal{I}_j$, we can assume that $N_j= N_j^i$ for all $i\in \mathcal{I}_j$. We may assume that $\mathcal{I}_{j+ 1}\subset \mathcal{I}_j$ for every $j$.

Therefore, by a diagonal argument, we find a common cofinal subnet of all $\mathcal{I}_j$, and denote it by $\mathcal{I}$. Set 
\begin{align}
S_{jk}\doteqdot B_{\infty}(z_{jk}, 2r_j)\backslash \cup_{l= 1}^{k- 1} B_{\infty}(z_{jl}, 2r_j)\ , \quad \quad 1\leq k\leq N_j  \nonumber
\end{align}

Define  
\begin{equation}\nonumber
\xi[x, a, b](y)= \left\{
\begin{array}{rl}
1 &\quad if \ \rho_{\infty}(x, y)\leq a \\
\frac{b- \rho_{\infty}(x, y)}{b- a} &\quad if \ a< \rho_{\infty}(x, y)< b \\
0 &\quad if \ \rho_{\infty}(x, y)\geq b 
\end{array} \right.
\end{equation}
We see that $\xi[x, a, b]$ is a Lipschitz function with Lipschitz constant $\frac{1}{b- a}$.

For any $\epsilon> 0$, $y\in K_{\infty}$, we define 
\begin{align}
\zeta_{S_{jk}}^{\epsilon} (y)= \xi [z_{jk}, r_j- 2\epsilon, r_j- \epsilon](y)\cdot \prod_{l= 1}^{k- 1}\Big\{1- \xi[z_{jl}, r_j- 2\epsilon, r_j- \epsilon]\Big\} \nonumber
\end{align}

It is easy to check that 
\begin{align}
\lim_{\epsilon\rightarrow 0+} |\zeta_{S_{jk}}^{\epsilon}- \chi_{S_{jk}}|_{L^2(K_{\infty})}= 0\ , \quad \lim_{\epsilon\rightarrow 0+}\lim_{i\rightarrow \infty} |\zeta_{S_{jk}}^{\epsilon}\circ \varphi_i- \chi_{S_{jk}^i}|_{L^2(K_i)}= 0 \nonumber
\end{align}
for $i\in \mathcal{I}$ and any $j= 1, 2, \cdots$, $k= 1, 2, \cdots, N_j$.

For $\bar{u}_{jk}= \lim_{i\rightarrow \infty} \bar{u}_{jk}^i$, we define two functions by
\begin{align}
\bar{u}_j(x)= \sum_{k= 1}^{N_j} \chi_{S_{jk}}(x)\bar{u}_{jk}\ , \quad \tilde{u}_j^{\epsilon}(x)= \sum_{k= 1}^{N_j} \zeta_{S_{jk}}^{\epsilon}(x)\bar{u}_{jk} \nonumber
\end{align}

Then 
\begin{align}
&\lim_{\epsilon\rightarrow \infty}\lim_{i\rightarrow \infty} |\tilde{u}_j^{\epsilon}- \bar{u}_j^i|_{L^2(K_i)} \nonumber \\
&\quad \leq  \lim_{\epsilon\rightarrow \infty}\lim_{i\rightarrow \infty}  \sum_{k= 1}^{N_j} \Big[|\bar{u}_{jk}|\cdot |\zeta_{S_{jk}}^{\epsilon}\circ \varphi_i- \chi_{S_{jk}^i}|_{L^2(K_i)}+ \nu_i(K_i)|\bar{u}_{jk}- u_{jk}^i|\Big]= 0 \nonumber
\end{align}
that is $\lim_{\epsilon\rightarrow \infty}\lim_{i\rightarrow \infty} |\tilde{u}_j^{\epsilon}- \bar{u}_j^i|_{L^2(K_i)} =0$.

Hence
\begin{align}
|\bar{u}_j- \bar{u}_{j'}|_{L^2}&\leq \lim_{\epsilon\rightarrow \infty} \Big(|\bar{u}_j- \tilde{u}_j^{\epsilon}|_{L^2}+ |\bar{u}_{j'}- \tilde{u}_{j'}^{\epsilon}|_{L^2}+ |\tilde{u}_j^{\epsilon}- \tilde{u}_{j'}^{\epsilon}|_{L^2}\Big) \nonumber \\
&\leq \lim_{\epsilon\rightarrow \infty}\lim_{i\rightarrow \infty} |\tilde{u}_j^{\epsilon}\circ \varphi_i- \tilde{u}_{j'}^{\epsilon}\circ \varphi_i|_{L^2} \nonumber \\
&\leq \lim_{i\rightarrow \infty} |\bar{u}_j^i- \bar{u}_{j'}^i|_{L^2} \nonumber \\
&\leq \lim_{i\rightarrow \infty} |\bar{u}_j^i- u_i|_{L^2}+ \lim_{i\rightarrow \infty} |\bar{u}_{j'}^i- u_i|_{L^2} \nonumber
\end{align}

From Claim \ref{claim 2.3.1}, we get that $\{\bar{u}_j\}$ is a Cauchy sequence in $L^2(K_{\infty})$, then set $u_{\infty}\doteqdot \lim_{j\rightarrow \infty} \bar{u}_j\in L^2(K_{\infty})$. From the above argument, it is easy to see that  $u_i\rightarrow u_{\infty}$ in $L^2$ sense on $K_{\infty}$, this completes the proof of Theorem \ref{thm 2.3}.
}
\qed

\section*{Acknowledgements}
Part of the work was done when the author visited University of California, Irvine, and he thanks the Department of Mathematics at UCI for their hospitality. The author would like to thank Peter Li and Laurent Saloff-Coste for their comments and suggestions, Shing-Tung Yau for his encouragement. He is deeply grateful to Jiaping Wang for his interest, continuous encouragement and support. He is also indebted to Yu Ding, Jozef Dodziuk, Renjie Feng, Jiaxin Hu, Jeff Streets, Luen-Fai Tam for their comments, Zhiqin Lu for the discussion, Yat Sun Poon and Li-Sheng Tseng for their support.

\end{document}